\documentclass[10pt]{article}
\usepackage{amsfonts,color}
\usepackage{amssymb}
\usepackage{footnote}
\usepackage{mathtools}
\usepackage{amsthm,wasysym}
\usepackage[english]{babel}
\usepackage{bbm}
\usepackage{colonequals}
\usepackage{setspace}
\usepackage{titlesec}
\usepackage{float}

\usepackage{graphicx}
\usepackage{amsfonts}
\usepackage{hyperref}
\usepackage{subcaption}
\usepackage{amsmath}
\usepackage{nicefrac}
\usepackage{gensymb}
\usepackage{enumerate}
\usepackage[nameinlink,capitalize]{cleveref}
\usepackage{chngcntr}
\usepackage{multirow}

\usepackage{mdframed}

\usepackage{cite}
\usepackage[nottoc,numbib]{tocbibind}

\usepackage{pgf,tikz,pgfplots}
\pgfplotsset{compat=1.14}
\usepackage{mathrsfs}
\usetikzlibrary{arrows}
\usetikzlibrary{arrows.meta}

\newcommand{\vb}[1]{\mathbf{#1}}
\newcommand{\bm}[1]{\boldsymbol{#1}}

\DeclareMathOperator{\range}{\mathrm{range}}

\DeclareMathOperator{\sym}{\mathrm{sym}}

\DeclareMathOperator{\vol}{\mathrm{vol}}
\DeclareMathOperator{\dev}{\mathrm{dev}}
\DeclareMathOperator{\skw}{\mathrm{skw}}
\DeclareMathOperator{\Anti}{\mathrm{Anti}}
\DeclareMathOperator{\axl}{\mathrm{axl}}
\DeclareMathOperator{\tr}{\mathrm{tr}}
\DeclareMathOperator{\Tr}{\mathrm{Tr}}

\newcommand{\one}{\bm{\mathbbm{1}}}
\newcommand{\con}[2]{\langle {#1} , \, {#2} \rangle}
\newcommand{\norm}[1]{\| {#1} \|}
\newcommand\subsetsim{\mathrel{%
  \ooalign{\raise0.2ex\hbox{$\subset$}\cr\hidewidth\raise-0.8ex\hbox{\scalebox{0.9}{$\sim$}}\hidewidth\cr}}}

\newcommand{\dd}{\mathrm{d}}
\newcommand{\D}{\mathrm{D}}
\DeclareMathOperator{\di}{\mathrm{div}}
\DeclareMathOperator{\Di}{\mathrm{Div}}
\DeclareMathOperator{\rot}{\mathrm{rot}}
\DeclareMathOperator{\Rot}{\mathrm{Rot}}
\DeclareMathOperator{\curl}{\mathrm{curl}}
\DeclareMathOperator{\Curl}{\mathrm{Curl}}
\DeclareMathOperator{\inc}{\mathrm{Inc}}
\DeclareMathOperator{\hess}{\mathrm{hess}}
\DeclareMathOperator{\airy}{\mathrm{airy}}

\newcommand{\so}{\mathfrak{so}}
\newcommand{\Sym}{\mathrm{Sym}}

\newcommand{\AW}{\mathcal{AW}}
\newcommand{\PS}{\mathcal{PS}}
\newcommand{\HHJ}{\mathcal{HHJ}}
\newcommand{\HZ}{\mathcal{HZ}}

\newcommand{\Leg}{\mathcal{L}}

\newcommand{\Po}{\mathit{P}}
\newcommand{\Le}{{\mathit{L}^2}}

\newcommand{\Hone}{\mathit{H}^1}
\newcommand{\Honez}{\mathit{H}_0^1}
\newcommand{\HsD}[1]{\mathit{H}^\mathrm{sym}(\mathrm{Div}{#1})}

\newcommand{\Hc}[1]{\mathit{H}(\mathrm{curl}{#1})}

\newcommand{\HdD}[1]{\mathit{H}(\mathrm{div}\mathrm{Div}{#1})}

\newcommand{\body}{V}
\newcommand{\surf}{A}
\newcommand{\curv}{s}

\newcommand{\RM}{\mathit{RM}}

\newcommand{\R}{\mathbb{R}}

\newcommand{\X}{\mathit{X}}
\newcommand{\Y}{\mathit{Y}}

\newcommand{\C}{\mathit{C}}
\renewcommand{\H}{\mathit{H}}

\newcommand{\Cm}{\mathbb{C}}
\newcommand{\J}{\mathbb{J}}
\newcommand{\A}{\mathbb{A}}

\newcommand{\wrt}{\text{w.r.t.}}

\allowdisplaybreaks[1]

\makeindex

\newtheoremstyle{break}
{\topsep}{\topsep}%
{\itshape}{}%
{\bfseries}{}%
{\newline}{}%
\theoremstyle{break}

\newtheorem{theorem}{Theorem}
\newtheorem{lemma}{Lemma}
\newtheorem{corollary}{Corollary}
\newtheorem{remark}{Remark}
\newtheorem{observation}{Observation}

\newmdtheoremenv{definition}{Definition}

\counterwithin{theorem}{section}
\counterwithin{lemma}{section}
\counterwithin{corollary}{section}
\counterwithin{remark}{section}
\counterwithin{observation}{section}
\counterwithin{definition}{section}

\counterwithin{figure}{section}
\counterwithin{equation}{section}

\makeatletter
\let\@fnsymbol\@arabic
\makeatother

\crefname{Problem}{Problem.}{Problem.}

\title{Symmetric unisolvent equations for linear elasticity purely in stresses}
\author{\normalsize{Adam Sky}\thanks{Corresponding author: Adam Sky, Institute of Computational Engineering and Sciences, Department of Engineering, Faculty of Science, Technology and Medicine, University of Luxembourg, 6 Avenue de la Fonte, L-4362 Esch-sur-Alzette, Luxembourg, email: adam.sky@uni.lu}
    \quad 
    \normalsize{and} \quad
	\normalsize{Andreas Zilian}\thanks{Andreas Zilian, Institute of Computational Engineering and Sciences, Department of Engineering, Faculty of Science, Technology and Medicine, University of Luxembourg, 6 Avenue de la Fonte, L-4362 Esch-sur-Alzette, Luxembourg, email: andreas.zilian@uni.lu}
}

\begin{document}

\maketitle

\begin{abstract}
In this work we introduce novel stress-only formulations of linear elasticity with special attention to their approximate solution using weighted residual methods. We present four sets of boundary value problems for a pure stress formulation of three-dimensional solids, and in two dimensions for plane stress and plane strain. The associated governing equations are derived by modifications and combinations of the Beltrami-Michell equations and the Navier-Cauchy equations.   
The corresponding variational forms of dimension $d \in \{2,3\}$ allow to approximate the stress tensor directly, without any presupposed potential stress functions, and are shown to be well-posed in $\Hone \otimes \Sym(d)$ in the framework of functional analysis via the Lax-Milgram theorem, making their finite element implementation using $\C^0$-continuous elements straightforward. Further, in the finite element setting we provide a treatment for constant and piece-wise constant body forces via distributions. The operators and differential identities in this work are provided in modern tensor notation and rely on exact sequences, making the resulting equations and differential relations directly comprehensible. Finally, numerical benchmarks for convergence as well as spectral analysis are used to test the limits and identify viable use-cases of the equations.   
\\
\vspace*{0.25cm}
\\
{\bf{Key words:}}  Beltrami-Michell equations, \and stress formulation, \and linear elasticity, \and numerical analysis, \and finite element method.
\\
\end{abstract}

\section{Introduction}
The classical formulation of linear elasticity is per the Navier-Cauchy equations, where the displacement field $\vb{u}$ is the main unknown. This formulation has been implemented countless times in finite element frameworks due to its simplicity, its straight-forward incorporation into the variational framework and the corresponding proofs of robustness in the compressible regime. However, as is well-known, the formulation is unstable for incompressible materials $\nu = 0.5$, since the Saint Venant-Kirchhoff material tensor tends to infinity $\Cm \to \infty$ in its components related to the trace of the strain tensor $\tr \bm{\varepsilon}$. A remedy for the latter phenomenon was introduced via the Hellinger-Reissner principle, which reformulated linear elasticity as a mixed problem, making use of the compliance tensor $\A$. The compliance tensor, being the inverse of the material tensor $\A = \Cm^{-1}$, remains finite for incompressible materials. Thus, a two-field variational formulation of linear elasticity in both displacements and stresses emerged, where the displacement field is set to be discontinuous $\vb{u} \in [\Le]^d$ and the stress tensor is symmetric and normal-continuous $\bm{\sigma} \in \HsD{}$. Subsequently, various finite element approaches were introduced for the mixed formulation \cite{Boon2021,Stenberg1988,arnold_mixed_2007}. However, the consistent solution relying on conforming finite elements based on commuting Hilbert space complexes \cite{arnold_complexes_2021,SkyNovel,CRMECA_2023__351_S1_A8_0} remained impractical, since the elements required for $\HsD{}$-discretisations did not allow for a reference-to-physical mapping of their base functions or were simply computationally too expensive. In fact, an efficient mapping procedure in two-dimensions was only recently introduced in \cite{SkyReissner} and implemented in the open source software for finite element computations NGSolve \cite{Sch2014,Sch1997} for the Hu-Zhang \cite{hu_family_2014,hu_family_2015,hu_finite_2016} element $\HZ \subset \HsD{}$, which satisfies the commutative property in the elasticity complex \cite{pauly_elasticity_2022,pauly_hilbert_2022,chen,Christiansen2023}. We mention that an alternative mapping by combination of mapped base functions was recently introduced in \cite{aznaran_transformations_2021} for the Arnold-Winther element $\AW \subset \HsD{}$ \cite{arnold_finite_2008,arnold_mixed_2002}, being the first conforming element for $\HsD{}$ in the elasticity complex. An alternative novel method called Tangential-Displacement-Normal-Normal-Stresses (TDNNS) \cite{pechstein_tdnns_2017,pechstein_analysis_2018,pechstein_anisotropic_2012,PECHSTEINTDNNS} alleviated the need for conforming $\HsD{}$-elements by applying the Hellan-Herrmann-Johnson principle, for which the two-dimensional Hellan-Herrmann-Johnson \cite{neunteufel_hellanherrmannjohnson_2019,arnold_hellan--herrmann--johnson_2020} elements $\HHJ \subsetsim \HdD{}$ were available and the three-dimensional Pechstein-Sch\"oberl elements $\PS \subsetsim \HdD{}$ were introduced. This innovative approach relied on an unconventional interpretation of the contraction between the divergence of the stress tensor and the displacement field as the dual product $\con{\Di \bm{\sigma}}{\vb{u}}_{\H^{-1}(\di) \times \Hc{}}$, where the stress field is defined as $\bm{\sigma} \in \HdD{}$, such that its divergence is in $\H^{-1}(\di)$ which is dual to $\Hc{}$, in which the displacement lives. The newly introduced finite elements for the stress field $\bm{\sigma}$ in this formulation are symmetric normal-normal continuous tensor fields, such that the contra-variant Piola transformation can be readily applied in order to map base functions from the reference to the physical element. For the displacement field $\vb{u} \in \Hc{}$ the formulation makes use of the classical N\'ed\'elec elements $\mathcal{N} \subset \Hc{}$ \cite{sky_hybrid_2021,SKY2022115298,SkyHigher,SkyOn,sky_polytopal_2022}.             

In a mixed formulation, the relation of the stress field to the displacement is given by the compliance tensor $\sym \D \vb{u} = \A \bm{\sigma}$. This relation can also be understood as the requirement that $\bm{\sigma}$ is to be related through the compliance tensor to some compatible tensor field $\A\bm{\sigma} \in \ker(\inc)$. This interpretation gives rise to the Beltrami-Michell equations \cite{CHANDRASEKHARAIAH1994363,Patnaik2000,Patnaik2006}, which replace the compliance relation with the compatibility relation $\inc(\A \bm{\sigma} ) = 0$, also known as the Saint-Venant compatibility condition. A core characteristic of these equations is that they are formulated purely in stresses without any additional fields. Thus, they represent the stress analogue of the displacement formulation, and proof of the equivalence of the two formulations can be found in \cite{Hackl1988} per the uniqueness of the solution.
Despite this desirable feature, the equations have found only limited numerical usage, specifically using finite differences \cite{WANG2022109752,10.2118/182595-MS,10.2118/193931-MS} or recently via Fourier series \cite{ANDRADE2024105309}, and in fact, application of the theory is almost exclusively analytical. Possibly, this is due to the fact that the equations are ill-suited for standard finite element computations. Indeed, in \cite{CIARLET2004307} the authors prove well-posedness of the variational formulation of the Beltrami-Michell equations for a discontinuous and compatible stress field $\bm{\sigma} \in [\Le \otimes \Sym(3)] \cap \ker(\inc)$ with $\inc\bm{\sigma} \in \H^{-2} \otimes \Sym(d)$, for the which the construction of conforming finite subspaces is not trivial. As a curiosity, we mention here that one other special case of a pure stress formulation can be derived specifically for the eigenvalue problem of linear elasticity \cite{Inzunza}, not yielding the Beltrami-Michell equations.

The Beltrami-Michell equations pose several difficulties for a straight-forward implementation in a standard finite element framework. Firstly, the three-dimensional equations are asymmetric, making the use of efficient solvers for symmetric positive definite problems inapplicable. Secondly, as we show in this work, for a mixed Dirichlet-Neumann boundary the associated variational form is unstable. Thirdly, the two-dimensional Beltrami-Michell equations for plane stress and plane strain characterise only the mean stress in the plane, making them incomplete in the sense that one cannot solve for the full stress tensor, but rather only its trace $\tr \bm{\sigma}$. We note that for the three-dimensional case the problem of asymmetry in the strong form was first noticed and alleviated in \cite{Pobedrya}. In spite of the non-existent finite element framework for the Beltrami-Michell equations, the latter have been employed successfully in multiple applications \cite{doi:10.2514/3.13034,KUCHER20053611,doi:10.2514/6.2003-2000}. Often, not the Beltrami-Michell equations themselves but rather their special-case descendants, the Beltrami, Maxwell, Morera and Airy \cite{FAAL20071133,Pommaret,JIANG2021220,Alshaya2021} stress equations, were used in order to analytically investigate the flow of stresses in a medium. Otherwise, the Beltrami-Michell equations in their original form have been used in the theory of mixtures \cite{MUTI2015140}, in order to describe coupled thermal-hydraulic-mechanical processes \cite{WANG2022109752,10.2118/182595-MS,10.2118/193931-MS}, and in the prediction of stresses in a reservoir \cite{ANDRADE2024105309}.

In this work we derive the Beltrami-Michell equations in modern tensor notation while relying on the concept of exact sequences for our argumentation, thus making the involved operators and differential relations clear. We symmetrise the three-dimensional equations, and combine a weaker version of planar static equilibrium with the two-dimensional Beltrami-Michell equations to introduce a novel full planar field problem for stresses in two-dimensions. By variational calculus we introduce the complete boundary value problems and their variational forms, for which we subsequently prove well-posedness on the Hilbert space $\Honez \otimes \Sym(d)$ via the Lax-Milgram theorem and several original lemmas. Finally, we investigate the numerical stability of the formulations via convergence estimates and spectral analysis, which leads us to present novel stabilised forms of the equations for mixed Dirichlet-Neumann boundary conditions. The novel equations introduced in this paper are summarised in \cref{ap:eqs}. 

The numerical benchmarks in this work are computed using the open source software NGSolve\footnote{www.ngsolve.org}, such that we rely on three- and two-dimensional conforming elements $\Leg^p(\body) \subset \Hone(\body)$, and $\Leg^p(\surf) \subset \Hone(\surf)$, constructed via hierarchical Legendre polynomials \cite{Zaglmayr2006}.

\subsection{Notation}
The following notation is used throughout this work.
Exception to these rules are made clear in the precise context.
\begin{itemize}
    \item Vectors are defined as bold lower-case letters $\vb{v}, \, \bm{\xi}$
    \item Second order tensors are denoted with bold capital letters $\bm{T}$
    \item Fourth-order tensors are designated by the blackboard-bold format $\mathbb{A}$
    \item We denote the Cartesian basis as $\{\vb{e}_1, \, \vb{e}_2, \, \vb{e}_3\}$
    \item Summation over indices follows the standard rule of repeating indices. Latin indices represent summation over the full dimension, whereas Greek indices define summation over the co-dimension.  
    \item The angle-brackets are used to define scalar products of arbitrary dimensions $\con{\vb{a}}{\vb{b}} = a_i b_i$, $\con{\bm{A}}{\bm{B}} = A_{ij}B_{ij}$
    \item The matrix product is used to indicate all partial-contractions between a higher-order and a lower-order tensor $\bm{A}\vb{v} = A_{ij} v_j \vb{e}_i$, $\mathbb{A}\bm{B} = A_{ijkl}B_{kl}\vb{e}_i \otimes \vb{e}_j$
    \item The second-order identity tensor is defined via $\one$, such that $\one \vb{v} = \vb{v}$. Analogously, the fourth-order identity tensor $\J$ yields $\J \bm{T} = \bm{T}$ 
    \item Subsequently, we define various differential operators based on the Nabla-operator $\nabla = \partial_i \vb{e}_i$, which is defined with respect to the dimension of the domain
    \item Volumes, surfaces and curves of the physical domain are identified via $\body$, $\surf$ and $\curv$, respectively.
    \item Tangent and normal vectors on the physical domain are designated by $\vb{t}$ and $\vb{n}$, respectively. 
\end{itemize}
We define the constant space of symmetric second order tensors as
\begin{align}
    \Sym(d) = \{ \bm{T} \in \R^{d \times d} \; | \; \bm{T} = \bm{T}^{T} \} \, .
\end{align}
Its counterpart is the space of skew-symmetric tensors
\begin{align}
    \so(d) = \{ \bm{T} \in \R^{d \times d} \; | \; \bm{T} = -\bm{T}^{T} \} \, .
\end{align}
Further, for our variational formulations we introduce the following Hilbert spaces and their respective norms 
\begin{subequations}
    \begin{align}
\Le(\body) &= \{ u : \body \to \mathbb{R} \; | \; \| u \|_\Le < \infty  \} \, , && \|u\|_\Le^2 = \int_\body u^2 \, \dd \body \, , \\
    \Hone(\body) &= \{ u \in \Le(\body) \; | \; \nabla u \in [\Le(\body)]^d \} \, , & \norm{u}^2_{\Hone} &= \norm{u}^2_\Le + \norm{\nabla u}^2_\Le \, . 
\end{align}
\end{subequations}
Hilbert spaces with vanishing traces are marked with a zero-subscript, for example $\Honez(\surf)$. 
Scalar products pertaining to the Hilbert spaces are indicated by a subscript on the angle-brackets
\begin{align}
    \con{u}{v}_{\Le} = \int_\body \con{u}{v} \, \dd \body \, ,
\end{align}
where the domain is clear from context.

In the following we define various differential operators derived from the $\nabla$-operator. 
\begin{itemize}
    \item The left-gradient is given via $\nabla$, such that $\nabla \lambda = \nabla \otimes \lambda$
    \item The right-gradient is defied for vectors and higher order tensors via $\D$, such that $\D \vb{v} = \vb{v} \otimes \nabla$
    \item The Hessian is given as the composite $\hess \lambda = \D \nabla \lambda$
    \item We define the vectorial divergence as $\di \vb{v} = \con{\nabla}{\vb{v}}$
    \item The tensor divergence is given by $\Di \bm{T} = \bm{T} \nabla$, implying a single contraction and acting row-wise 
    \item The Laplacian is given via $\Delta \lambda = \di \nabla \lambda$ for scalars and as $\Delta \bm{T} = \Di \D \bm{T}$ for higher order tensors
    \item The vectorial curl operator reads $\curl \vb{v} = \nabla \times \vb{v}$
    \item For tensors the operator is given by $\Curl \bm{T} = -\bm{T} \times \nabla$, acting row-wise
    \item The composite incompatibility operator reads $\inc \bm{T} = \curl \Curl \bm{T} = - \nabla \times \bm{T} \times \nabla$
    \item In two dimensions the vectorial $\curl$-operator induces the scalar operator $\nabla^\perp \lambda = \bm{R} \nabla \lambda$ and vectorial operator $\rot \vb{v} = \di(\bm{R}\vb{v})$, with $\bm{R} = \vb{e}_1 \otimes \vb{e}_2 - \vb{e}_2 \otimes \vb{e}_1$
    \item Analogously, the tensorial $\Curl$-operator induces the vectorial operator $\D^\perp \vb{v} = (\D \vb{v}) \bm{R}^T$ and the tensorial operator $\Rot \bm{T} = \Di(\bm{T}\bm{R}^T)$, acting row-wise
    \item Consequently, the composite incompatibility operator reduces to $\rot \Rot \bm{T}$ for tensors and the Airy-operator $\airy \lambda = \D^\perp \nabla^\perp \lambda$ for scalars 
\end{itemize}

Lastly, in this work we employ the standard algebraic operators
\begin{align}
    &\sym \bm{T} = \dfrac{1}{2}(\bm{T} + \bm{T}^T) \, , 
    && \skw \bm{T} = \dfrac{1}{2}(\bm{T} - \bm{T}^T) \, , 
    &&\tr \bm{T} = \con{\bm{T}}{\one} \, , \notag \\
    &\vol \bm{T} = \dfrac{1}{d} (\tr \bm{T}) \one \, ,
    &&\dev \bm{T} = \bm{T} - \vol \bm{T} \, , && \one,\bm{T} \in \R^{d \times d} \, ,
\end{align}
defining the symmetry, skew-symmetry, trace, volumetric and deviatoric operators on tensors. We generalise the trace operator to higher order tensors with the notation
\begin{align}
    &\Tr \bm{T} = T_{ijk} \vb{e}_i \con{\vb{e}_j \otimes \vb{e}_k}{\one} = T_{ijj} \vb{e}_i \in \R \, , && \bm{T} \in \R^{d \times d \times d} \, .  
\end{align}
For skew-symmetric tensors we also introduce the operators
\begin{align}
    \vb{a} \times \vb{v} = (\Anti \vb{a}) \vb{v} = [\axl (\Anti \vb{a})] \times \vb{v} \in \R \, ,   
\end{align}
where $(\Anti \vb{a}) \in \so(3) \subset \R^{3 \times 3}$ maps the vector $\vb{a} \in \R^3$ to its corresponding skew-symmetric tensor for an equivalent form of the cross product given through a single contraction, and $\axl(\Anti \vb{a}) = \vb{a} \in \R^3$ is the inverse operator extracting the axial vector.    

\section{Linear elasticity}
In this section we briefly recap the Navier-Cauchy equations of linear elasticity. These equations are subsequently used in the derivation of the Beltrami-Michell equations.

\subsection{Static equilibrium}
The energy functional of linear elasticity with the linear strain measure $\bm{\varepsilon} = \sym \D \vb{u}$ is given by 
\begin{align}
    &I(\vb{u}) = \int_\body \dfrac{1}{2} \langle \sym\D\vb{u}, \, \Cm \sym\D\vb{u} \rangle  - \langle \vb{u} ,\, \vb{f} \rangle \, \dd \body - \int_{\surf_N} \con{\vb{u}}{\bm{\mu}} \, \dd \surf \quad \to \quad \min \quad \wrt \quad \vb{u} \, , 
\end{align}
where $\vb{u} : \overline{\body} \subset \R^3 \to \R^3$ is the displacement field, $\Cm \in \R^{3 \times 3 \times 3 \times 3}$ is the tensor of material constants, $\vb{f} : \overline{\body} \subset \R^3 \to \R^3$ represents the body forces, and $\bm{\mu} : \surf_N \subset \R^2 \to \R^3$ are the tractions on the Neumann boundary $\surf_N$, see \cref{fig:domain}.
\begin{figure}
		\centering
		\definecolor{asl}{rgb}{0.4980392156862745,0.,1.}
		\definecolor{asb}{rgb}{0.,0.4,0.6}
		\begin{tikzpicture}[line cap=round,line join=round,>=triangle 45,x=1.0cm,y=1.0cm]
			\clip(-0.5,-0.5) rectangle (13,4.5);
			
			\fill [asb, opacity=0.1] plot [smooth cycle] coordinates {(1,3) (3,4) (7, 2) (10,3) (12,1) (10,0) (5,1) (2,1)};
			
			\begin{scope}
				\clip(5,-0.5) rectangle (12.5,1.5);
				\draw [asl, dashed] plot [smooth cycle] coordinates {(1,3) (3,4) (7, 2) (10,3) (12,1) (10,0) (5,1) (2,1)};
			\end{scope}
		    \begin{scope}
		    	\clip(5,1.5) rectangle (12.5,4.5);
		    	\draw [asl, dashed] plot [smooth cycle] coordinates {(1,3) (3,4) (7, 2) (10,3) (12,1) (10,0) (5,1) (2,1)};
		    \end{scope}
			\begin{scope}
				\clip(0,-0.5) rectangle (5,4.5);
				\draw [asb] plot [smooth cycle] coordinates {(1,3) (3,4) (7, 2) (10,3) (12,1) (10,0) (5,1) (2,1)};
			\end{scope}
			
			\draw [-to,color=black,line width=1.pt] (0,0) -- (1,0);
			\draw [-to,color=black,line width=1.pt] (0,0) -- (0,1);
			\draw (1,0) node[color=black,anchor=west] {$x$};
			\draw (0,1) node[color=black,anchor=south] {$y$};
			
			\draw [-to,color=asb,line width=1.pt] (10.,1.2) -- (10.3,0.9);
			\draw [-to,color=asb,line width=1.pt] (9.6,1.2) -- (9.3,0.9);
			\draw [-to,color=asb,line width=1.pt] (9.8,1.2) -- (9.8,0.7);
			
			\draw [-to,color=asl,line width=1.pt] (11.35,2) -- (12,2.55);
			\draw (12,2.55) node[color=asl,anchor=south west] {$\vb{n}$};
			\draw (9.8,1.2) node[color=asb,anchor=south] {$\vb{f}$};
            \draw [to-,color=asb,line width=1.pt] (11.35,2) -- (12.5,1.5);
			\draw (12.5,1.65) node[color=asb,anchor=north west] {$\bm{\mu}$};
			
			\draw (6.5,1.15) node[color=asb,anchor=south] {$\body$};
			\draw (4.3,3.6) node[color=asb,anchor=west] {$\surf_D$};
			\draw (10.6,2.95) node[color=asl,anchor=south] {$\surf_N$};
		\end{tikzpicture}
		\caption{The domain $\body \subset \R^3$ with Dirichlet $\surf_D\subset \R^2$ and Neumann $\surf_N\subset \R^2$ boundaries under internal body forces $\vb{f}$ and tractions $\bm{\mu}$ with the outer normal vector $\vb{n}$.}
		\label{fig:domain}
	\end{figure}
For a linear isotropic Saint Venant-Kirchhoff material the  elasticity tensor reads
\begin{align}
    &\bm{\sigma} = \Cm \sym\D\vb{u} \, , && \Cm =  \dfrac{E}{1 + \nu} \left [\J + \dfrac{\nu}{1-2\nu} \one \otimes \one \right ]  \, , 
\end{align}
where $\{E, \nu \}$ are Young's modulus and the Poisson ratio, $\one : \R^3 \to \R^3$ is the second order identity tensor, and $\J: \R^{3\times 3} \to \R^{3\times 3}$ is the fourth order identity tensor. Variation of the energy functional with respect to the displacement field yields
\begin{align}
    & \int_\body \con{ \sym\D\delta\vb{u}}{\Cm \sym\D\vb{u}}\, \dd \body = \int_\body \con{\delta \vb{u}}{ \vb{f}} \, \dd \body + \int_{\surf_N} \con{\delta \vb{u}}{\bm{\mu}} \, \dd \surf \qquad \forall \, \delta \vb{u} \in [\C^\infty(\overline{\body})]^3 \, . 
\end{align}
Applying partial integration leads to 
\begin{align}
    & \int_{\partial \body} \con{\delta \vb{u}}{(\Cm \sym\D\vb{u}) \,\vb{n}} \, \dd \surf  - \int_\body \con{\delta \vb{u}}{ \Di (\Cm \sym \D \vb{u}) } \, \dd \body = \int_\body \con{\delta \vb{u}}{ \vb{f}} \, \dd \body + \int_{\surf_N} \con{\delta \vb{u}}{\bm{\mu}} \, \dd \surf \quad \forall \, \delta \vb{u} \in [\C^\infty(\overline{\body})]^3 \, , 
\end{align}
from which we extract the boundary value problem of linear elasticity in three dimensions  
\begin{subequations}
    \begin{align}
        -\Di (\Cm \sym \D \vb{u}) &= \vb{f} && \text{in} && \body \, , \\
        (\Cm \sym \D \vb{u}) \, \vb{n} &= \bm{\mu} && \text{on} && \surf_N \, , \\
        \vb{u} &= \widetilde{\vb{u}} && \text{on} && \surf_D \, ,
    \end{align}
\end{subequations}
by splitting the boundary between Dirichlet and Neumann $\partial \body = \surf_D \cup \surf_N$ with $\surf_D \cap \surf_N = \emptyset$.

\subsection{Planar static equilibrium}
In the planar form of linear elasticity it is assumed that the displacement vector is a function solely of the $x-y$-plane $\vb{u} = \vb{u}(x,y)$ with relevant components only in the $x$- and $y$-directions $\vb{u}:\overline{\surf} \subset \R^2 \to \R^2$. Under the assumption of plane stress $\sigma_{i3} = \sigma_{3i} = 0$, the material tensor is adapted to 
\begin{align}
    &\bm{\sigma} = \Cm \sym \D \vb{u} \, , && \Cm = \dfrac{E}{1-\nu^2}[(1-\nu)\J + \nu  \one\otimes \one] \in \R^{2 \times 2\times 2\times 2} \, ,
\end{align}
where now $\J: \R^{2 \times 2\times 2\times 2} \to \R^{2 \times 2\times 2\times 2}$ and $\one: \R^{2 \times 2} \to \R^{2 \times 2}$.
In the case of plane strain $\varepsilon_{i3} = \varepsilon_{3i} = 0$ with $\bm{\varepsilon} = \sym \D \vb{u}$, the material tensor has the same form as for the three-dimensional case, reduced to its planar components $\Cm \in \R^{2\times 2\times 2\times 2}$. Since in both cases the quadratic form $\con{\sym\D\vb{u}}{\bm{\sigma}} = \con{\sym\D\vb{u}}{\Cm \sym \D \vb{u}}$ does not produce energy for the out-of-plane components, the tensor fields of the strain and stress can be reduced to $\sym\D\vb{u} : \overline{\surf} \subset \R^2 \to \Sym(2)$ and $\bm{\sigma} : \overline{\surf} \subset \R^2 \to \Sym(2)$, where $\Sym(2)$ is the space of constant symmetric second order tensors in $\R^{2 \times 2}$. Further, the form of the boundary value problem remains the same as in three dimensions
\begin{subequations}
    \begin{align}
        -\Di (\Cm \sym \D \vb{u}) &= \vb{f} && \text{in} && \surf \, , \\
        (\Cm \sym \D \vb{u}) \, \vb{n} &= \bm{\mu} && \text{on} && \curv_N \, , \\
        \vb{u} &= \widetilde{\vb{u}} && \text{on} && \curv_D \, ,
    \end{align}
\end{subequations}
with an adjustment of the domain and boundary to $\partial \surf = \curv_D \cup \curv_N$, and the body forces and tractions to $\vb{f} : \overline{\surf} \subset \R^2 \to \R^2$ and $\bm{\mu} : \curv_N \subset \R \to \R^2$.

\subsection{The Beltrami-Michell equations}
The Beltrami-Michell equations reformulate isotropic linear elasticity purely in stresses.
In order to derive the Beltrami-Michell equations we introduce the compliance tensor for the inverse stress-strain relation
\begin{align}
     &\sym \D \vb{u} = \A \bm{\sigma} \, , &&  \A = \Cm^{-1} = \dfrac{1}{E} [ (1+\nu) \J - \nu  \one \otimes \one] \, .
\end{align}
With the strong form of equilibrium and constitutive relation in strain-form 
\begin{subequations}
    \begin{align}
        -\Di \bm{\sigma} &= \vb{f} && \text{in} && \body \, , \\
        \A \bm{\sigma} - \sym \D \vb{u} &= 0  && \text{in} && \body \, , 
    \end{align}
\end{subequations}
the two-field problem of displacement and stresses as per the Hellinger-Reissner principle \cite{SkyReissner} is obtained. In order to construct a problem purely in $\bm{\sigma}$, one starts from the Saint-Venant compatibility condition
\begin{align}
    \inc \bm{\varepsilon} = \curl \Curl \bm{\varepsilon} =  0 \qquad \forall \, \bm{\varepsilon} \in \sym \D [\C^\infty(\body)]^3 \, ,
\end{align}
which asserts that the strain tensor $\bm{\varepsilon}$ is the symmetric part of the gradient of some vector field on contractible domains by the exact sequence property of the elasticity complex \cite{pauly_elasticity_2022,pauly_hilbert_2022,chen,Christiansen2023}, see \cref{fig:elasticity}.
\begin{figure}
		\centering
		\begin{tikzpicture}[line cap=round,line join=round,>=triangle 45,x=1.0cm,y=1.0cm]
				\clip(0,-0.25) rectangle (16,0.5);
                \draw (0,0) node[anchor=west] {$\RM(\body) \hookrightarrow [\C^\infty(\body)]^3$};
                \draw [-Triangle,line width=.5pt] (3.5,0.) -- (5,0.);
                \draw (0.7+3.5,0) node[anchor=south] {$\sym \D$};
                \draw (5,0) node[anchor=west] {$\C^\infty(\body)\otimes\Sym(3)$};
                \draw [-Triangle,line width=.5pt] (8,0.) -- (9.5,0.);
                \draw (0.7+8,0) node[anchor=south] {$\inc$};
                \draw (9.5,0) node[anchor=west] {$\C^\infty(\body)\otimes\Sym(3)$};
                \draw [-Triangle,line width=.5pt] (12.5,0.) -- (14,0);
                \draw (0.7+12.5,0) node[anchor=south] {$\Di$};
                \draw (14,0) node[anchor=west] {$[\C^\infty(\body)]^3$};
\end{tikzpicture}
		\caption{The three-dimensional elasticity sequence. The space of rigid body motions is the kernel of the symmetrised gradient operator $\ker(\sym\D) = \RM(\body)$, and the range of the latter is the kernel of the incompatibility operator $\sym \D [\C^\infty(\body)]^3 = \ker(\inc)$. 
        The range of the incompatibility operator is exactly the kernel of the divergence operator $\inc[\C^\infty(\body)\otimes \Sym(3)] = \ker(\Di)$ for symmetric tensors, and the range of the divergence operator is a surjection onto the last space in the sequence $\Di[\C^\infty(\body)\otimes \Sym(3)] = [\C^\infty(\body)]^3$.}
		\label{fig:elasticity}
	\end{figure}
Note that the compatibility condition does not identify the underlying vector field uniquely since $\inc (\sym \D \vb{u}) = \inc [\sym \D (\vb{u} + \vb{v})] = 0$ for any pair $\{\vb{u}, \vb{v}\} \in [\C^\infty(\body)]^3$ \cite{ANDRIANOV2022123}. Consequently, the equation of static equilibrium is paramount to determine a unique solution field in the domain $\body$. The field equations are now given by
\begin{subequations}
    \begin{align}
        -\Di \bm{\sigma} &= \vb{f} && \text{in} && \body \, , \label{eq:equi} \\
        \inc (\A \bm{\sigma})&= 0  && \text{in} && \body \, . 
    \end{align}
\end{subequations}
The three-dimensional compliance relation $\A \bm{\sigma}$ can be expanded into
\begin{align}
    \A \bm{\sigma} = \dfrac{1}{E} [(1 + \nu) \bm{\sigma} - (\nu \tr \bm{\sigma}) \one] \, .
\end{align}
Since the incompatibility equation is zero on the right-hand-side, for a medium composed of a single isotropic material we can multiply the equation with Young's modulus $E$ and divide it by $(1+\nu)$ before inserting the compliance relation 
\begin{align}
    \dfrac{E}{1+\nu}\inc (\A \bm{\sigma})&=  \inc \bm{\sigma} - \dfrac{\nu}{1+\nu} \inc[(\tr \bm{\sigma}) \one] = 0 \, . \label{eq:Ainc}
\end{align}
The application of the incompatibility operator \cite{Marguerre} on a symmetric tensor can be expressed as 
\begin{align}
    \inc \bm{\sigma} = 2\sym (\D \Di \bm{\sigma}) - \Delta \bm{\sigma} - \hess (\tr \bm{\sigma}) + [ \Delta \tr \bm{\sigma}  - \di \Di \bm{\sigma}] \one \, ,
    \label{eq:incid}
\end{align}
as per the Schaefer-Kr\"oner formula \cite{Schaefer,Kroner1954}.
For a volumetric tensor such as $(\tr \bm{\sigma}) \one$ the result simplifies \cite{Lewintan2021} to 
\begin{align}
    \inc [ (\tr \bm{\sigma}) \one ] = (\Delta \tr \bm{\sigma}) \one -  \hess (\tr \bm{\sigma})  \, .
    \label{eq:hess}
\end{align}
Inserting the latter two identities into \cref{eq:Ainc} yields
\begin{align}
    -\Delta \bm{\sigma} - \dfrac{1}{1+\nu} \hess( \tr \bm{\sigma} ) + \dfrac{1}{1+\nu} (\Delta \tr \bm{\sigma}) \one - (\di \Di \bm{\sigma}) \one + 2 \sym \D \Di \bm{\sigma} = 0 \, ,
    \label{eq:laplacesig}
\end{align}
where we can now directly insert the static equilibrium equation \cref{eq:equi} to find 
\begin{align}
    -\Delta \bm{\sigma} - \dfrac{1}{1+\nu} \hess( \tr \bm{\sigma} ) + \dfrac{1}{1+\nu} (\Delta \tr \bm{\sigma}) \one  = 2 \sym \D \vb{f} - (\di \vb{f}) \one \, . \label{eq:almost-beltrami}
\end{align}
The equation can be further simplified by observing that the first invariant of the vanishing field $\inc (\A \bm{\sigma}) = 0$ must also vanish $\tr [\inc (\A \bm{\sigma})] = 0$.
With 
\begin{subequations}
    \begin{align}
    \tr (\inc \bm{\sigma}) &= \Delta \tr \bm{\sigma} - \di \Di \bm{\sigma} \, , \\
    \tr (\inc [ (\tr \bm{\sigma}) \one ])  & = 2 \Delta \tr \bm{\sigma} \, ,
\end{align}
\end{subequations}
the first invariant leads to
\begin{align}
    (1-\nu) \Delta \tr \bm{\sigma} - (1+\nu) \di \Di \bm{\sigma} = 0 \, . \label{eq:invone}
\end{align}
With the latter identity at hand we can eliminate $\Delta \tr \bm{\sigma}$ in \cref{eq:almost-beltrami} and retrieve the classical form of the Beltrami-Michell field equation
\begin{align}
    -\Delta \bm{\sigma} - \dfrac{1}{1+\nu} \hess (\tr \bm{\sigma}) &= 2 \sym \D \vb{f} + \dfrac{\nu}{1-\nu} (\di \vb{f}) \one  && \text{in} && \body \, . 
    \label{eq:beltrami}
\end{align}
Clearly, the Beltrami-Michell equations introduce the problem of finding meaningful stress fields for certain force fields $\vb{f}$. In fact, we can characterise cases for which $\bm{F}(\vb{f}) = 2 \sym \D \vb{f} + \nu/(1-\nu) (\di \vb{f}) \one$ vanishes, leaving no right-hand side $\bm{F}(\vb{f}) = 0$. 
\begin{remark}[Vanishing right-hand side] \label{le:vanishing-rhs}
    Let the body forces $\vb{f}$ be in the polynomial space of rigid body motions, there holds
    \begin{align}
        &\bm{F}(\vb{f}) = 0 \qquad \forall\,\vb{f} \in \RM(\body) = \R^3 \oplus \R^3 \times \vb{x} \, , && \vb{x} = x \vb{e}_1 + y \vb{e}_2 + z \vb{e}_3 \, ,  
    \end{align}
    asserting that the right-hand side of the Beltrami-Michell equation vanishes. This is clear since the overlap of the kernel of the divergence operator $\ker(\di) = \range(\curl)$ and the symmetrised gradient operator $\ker(\sym\D) = \RM(\body)$ is the space of rigid body motions $\ker(\di) \cap \ker(\sym \D) = \RM(\body)$. The respective characterisations are given by the exact de Rham sequence \cite{PaulyDeRham,Demkowicz2000} in \cref{fig:deRham} 
    \begin{figure}
		\centering
		\begin{tikzpicture}[line cap=round,line join=round,>=triangle 45,x=1.0cm,y=1.0cm]
				\clip(1,-0.25) rectangle (13,0.5);
                \draw (3.5,0) node[anchor=east] {$\R \hookrightarrow \C^\infty(\body)$};
                \draw [-Triangle,line width=.5pt] (3.5,0.) -- (5,0.);
                \draw (0.7+3.5,0) node[anchor=south] {$\nabla$};
                \draw (5,0) node[anchor=west] {$[\C^\infty(\body)]^3$};
                \draw [-Triangle,line width=.5pt] (6.8,0.) -- (8.3,0.);
                \draw (0.7+6.8,0) node[anchor=south] {$\curl$};
                \draw (8.3,0) node[anchor=west] {$[\C^\infty(\body)]^3$};
                \draw [-Triangle,line width=.5pt] (10.1,0.) -- (11.6,0);
                \draw (0.7+10.1,0) node[anchor=south] {$\di$};
                \draw (11.6,0) node[anchor=west] {$\C^\infty(\body)$};
\end{tikzpicture}
		\caption{The three-dimensional de Rham sequence. The space of constants is the kernel of the gradient operator $\R = \ker(\nabla)$, and the range of the latter is the kernel of the curl operator $\nabla \C^\infty(\body) = \ker(\curl)$. 
        The range of the curl operator is exactly the kernel of the divergence operator $\curl[\C^\infty(\body)]^3 = \ker(\di)$, of which the range is a surjection onto the last space $\di[\C^\infty(\body)]^3 = \C^\infty(\body)$.}
		\label{fig:deRham}
	\end{figure}
    and the elasticity sequence in \cref{fig:elasticity}.
\end{remark}

\subsection{The Beltrami stress function}
Under the assumption of vanishing body forces in the domain $\vb{f} = 0$, the stress field $\bm{\sigma}$ must clearly belong to the null space $\ker(\Di)$, thus satisfying static equilibrium a priori $\Di \bm{\sigma} = 0$.
By the exact elasticity sequence there holds
\begin{align}
    \forall \, \bm{\sigma} \in [\C^\infty(\body)\otimes\Sym(3)] \cap \ker(\Di) \qquad \exists \, \bm{\Sigma} \in \C^\infty(\body)\otimes\Sym(3): \qquad \bm{\sigma} = \inc \bm{\Sigma} \, ,
\end{align}
on contractible domains as per \cref{fig:elasticity}.
Consequently, we can insert the tensorial potential $\inc\bm{\Sigma}$ into the Beltrami-Michell equations to find 
\begin{align}
    -\Delta [\Delta \bm{\Sigma} - 2  \sym \D \di \bm{\Sigma} - ( \Delta \tr \bm{\Sigma} -  \di \Di \bm{\Sigma}) \one ] - \dfrac{1}{1+\nu} \hess(\nu\Delta \tr \bm{\Sigma} + \di \Di \bm{\Sigma}) = 0 \, ,
\end{align}
which is a fourth-order partial differential equation with $\bm{\Sigma}$ being the Beltrami stress function.
\begin{remark}[Specialised stress functions]
    We note that the Maxwell stress function $\bm{\Sigma} =  \Sigma_1 \vb{e}_1 \otimes \vb{e}_1 + \Sigma_2 \vb{e}_2 \otimes \vb{e}_2 + \Sigma_3 \vb{e}_3 \otimes \vb{e}_3$ and the Morera stress function $\bm{\Sigma} = \Sigma_{12}\vb{e}_1 \otimes \vb{e}_2 + \Sigma_{13}\vb{e}_1 \otimes \vb{e}_3 + \Sigma_{23}\vb{e}_2 \otimes \vb{e}_3$ are specialised forms of the general Beltrami stress function \cite{Pommaret}, limiting the number of independent terms in the Beltrami stress tensor.
\end{remark}

\subsection{The planar Beltrami-Michell equations}
In two dimensions the compatibility condition is given by the operator $\rot\Rot \bm{\varepsilon} = 0$, and the equilibrium equation is simply adapted to a two-dimensional plane 
\begin{align}
    -\Di \bm{\sigma} &= \vb{f} && \text{in} && \surf \, , 
    \label{eq:equi-2D}
\end{align}
where $\bm{\sigma}:\surf \to \R^{2 \times 2}$ and $\vb{f}:\surf \to \R^2$.
The out-of-plane strain can be retrieved with the classical relation $\varepsilon_{33} = -(\nu/E)\tr \bm{\sigma}$. The case of plane stress follows analogously to the three-dimensional case since the compliance relation is the same, adjusted to two-dimensions $\A  = (1/E) [(1+ \nu) \mathbb{J} - \nu \one \otimes \one ]\in \R^{2 \times 2 \times 2 \times 2}$. 
For a symmetric tensor there holds
\begin{align}
    \rot\Rot \bm{\varepsilon} = \Delta \tr \bm{\varepsilon} - \di \Di \bm{\varepsilon} = 0 \qquad \forall \, \bm{\varepsilon} \in \sym \D [\C^\infty(\surf)]^2 \, ,
\end{align}
as per the two-dimensional elasticity complex, \cref{fig:elast}.
\begin{figure}
		\centering
		\begin{tikzpicture}[line cap=round,line join=round,>=triangle 45,x=1.0cm,y=1.0cm]
				\clip(0,-0.25) rectangle (11.5,0.5);
                \draw (0,0) node[anchor=west] {$\RM(\surf) \hookrightarrow [\C^\infty(\surf)]^2$};
                \draw [-Triangle,line width=.5pt] (3.5,0.) -- (5,0.);
                \draw (0.7+3.5,0) node[anchor=south] {$\sym \D$};
                \draw (5,0) node[anchor=west] {$\C^\infty(\surf)\otimes\Sym(2)$};
                \draw [-Triangle,line width=.5pt] (8,0.) -- (9.5,0.);
                \draw (0.7+8,0) node[anchor=south] {$\rot \Rot$};
                \draw (9.5,0) node[anchor=west] {$\C^\infty(\surf)$};
\end{tikzpicture}
        \vspace{0.1cm}
        \begin{tikzpicture}[line cap=round,line join=round,>=triangle 45,x=1.0cm,y=1.0cm]
				\clip(0,-0.25) rectangle (11.5,0.5);
                \draw (3.5,0) node[anchor=east] {$\Po^1(\surf) \hookrightarrow \C^\infty(\surf)$};
                \draw [-Triangle,line width=.5pt] (3.5,0.) -- (5,0.);
                \draw (0.7+3.5,0) node[anchor=south] {$\airy$};
                \draw (5,0) node[anchor=west] {$\C^\infty(\surf)\otimes\Sym(2)$};
                \draw [-Triangle,line width=.5pt] (8,0.) -- (9.5,0.);
                \draw (0.7+8,0) node[anchor=south] {$\Di$};
                \draw (9.5,0) node[anchor=west] {$[\C^\infty(\surf)]^2$};
\end{tikzpicture}
		\caption{Two planar elasticity sequences derived from the reduction of the incompatibility operator. For two-dimensional tensors, the incompatibility  operator turns into the $\rot \Rot$-operator and for scalars into the airy operator.
        In the first sequence the two-dimensional space of rigid body motions is the kernel of the $\sym\D$-operator $\RM(\surf) = \ker(\sym \D)$, of which the range is the kernel of the $\rot \Rot$-operator $\sym \D [\C^\infty(\surf)]^2 = \ker(\rot \Rot)$. Finally, the $\rot\Rot$-operator yields a surjection onto the last space in the sequence $\rot \Rot[\C^\infty(\surf) \otimes \Sym(2)] = \C^\infty(\surf)$. In the second sequence the kernel of the $\airy$-operator is given by the linear polynomial space $\Po^1(\surf) = \ker(\airy)$, and the $\airy$-operator maps the kernel the divergence operator for symmetric tensors $\airy \C^\infty(\surf) = \ker(\Di)$. Lastly, the divergence yields a surjection onto the last space in the sequence $\Di [\C^\infty(\surf) \otimes \Sym(2)] = [\C^\infty(\surf)]^2$. 
        }
		\label{fig:elast}
	\end{figure}
Applying the latter to $\A \bm{\sigma}$ yields
\begin{align}
    \dfrac{E}{1+\nu} \rot\Rot (\A \bm{\sigma}) = \dfrac{1}{1+\nu} \Delta \tr \bm{\sigma} - \di \Di \bm{\sigma} = 0 \, , 
\end{align}
since $\rot\Rot[(\tr \bm{\sigma}) \one] = \Delta \tr \bm{\sigma}$. Now inserting \cref{eq:equi-2D} results in the two-dimensional field equations 
\begin{align}
    -\dfrac{1}{1+\nu} \Delta \tr \bm{\sigma} &= \di \vb{f} && \text{in} && \surf \, .
    \label{eq:planestress}
\end{align}

In the case of plane strain one assumes zero out-of-plane strains, such that the compliance equation reads
\begin{align}
    &\bm{\varepsilon} = \A \bm{\sigma} = \dfrac{1+\nu}{E}[\bm{\sigma} - (\nu \tr \bm{\sigma})\one] \, , &&
    \A = \dfrac{1+\nu}{E}[\J - \nu \one \otimes \one] \, , 
\end{align}
for the in-plane components of the stain tensor.
The out-of-plane stress can be recovered via $\sigma_{33} = \lambda \tr (\sym \D \vb{u}) = \lambda \di \vb{u}\, ,$ with the Lam\'e constant $\lambda = E\nu /[(1+\nu)(1-2\nu)]$, but does not produce energy.
Consequently, we find for the compatibility
\begin{align}
    \dfrac{E}{1+\nu} \rot\Rot (\A \bm{\sigma}) = (1-\nu)\Delta \tr \bm{\sigma} - \di \Di \bm{\sigma} = 0 \, ,
\end{align}
which by incorporating static equilibrium results in the field equation
\begin{align}
    -(1-\nu) \Delta \tr \bm{\sigma} &= \di \vb{f} && \text{in} && \surf \, .
    \label{eq:planestrain}
\end{align}

\begin{remark}[Vanishing planar right-hand side]
    For both plane stress and plane strain, the right-hand side of the equations is given by the divergence of the body forces $\di \vb{f}$. Consequently, solenoidal force fields $\vb{f} = \nabla^\perp f$ are not captured by the equations since $\ker(\di) = \range(\nabla^\perp)$ in two dimensions by the exact de Rham complex as per \cref{fig:deRham2D}.
\end{remark}
\begin{figure}
		\centering
		\begin{tikzpicture}[line cap=round,line join=round,>=triangle 45,x=1.0cm,y=1.0cm]
				\clip(1,-0.25) rectangle (10,0.5);
                \draw (3.5,0) node[anchor=east] {$\R \hookrightarrow \C^\infty(\surf)$};
                \draw [-Triangle,line width=.5pt] (3.5,0.) -- (5,0.);
                \draw (0.7+3.5,0) node[anchor=south] {$\nabla^\perp$};
                \draw (5,0) node[anchor=west] {$[\C^\infty(\surf)]^2$};
                \draw [-Triangle,line width=.5pt] (6.8,0.) -- (8.3,0.);
                \draw (0.7+6.8,0) node[anchor=south] {$\di$};
                \draw (8.3,0) node[anchor=west] {$\C^\infty(\surf)$};
\end{tikzpicture}
		\caption{One of two possible two-dimensional de Rham sequences derived by the reduction of the curl operator to two dimensions. The space of constants is the kernel of the rotated gradient operator $\R = \ker(\nabla^\perp)$, and the range of the latter is the kernel of the divergence operator $\nabla^\perp\C^\infty(\surf) = \ker(\di)$. The range of the divergence operator is a surjection onto the last space $\di[\C^\infty(\surf)]^2 =  \C^\infty(\surf)$.}
		\label{fig:deRham2D}
	\end{figure}

\subsection{The Airy function}
In the special case where in the two-dimensional Beltrami-Michell equations no body forces occur in the domain $\surf$, the field equations simplify to
\begin{subequations}
    \begin{align}
        \Di \bm{\sigma} &= 0 && \text{in} && \surf \, , \\
        \Delta \tr \bm{\sigma} &= 0 && \text{in} && \surf \, .
    \end{align}
\end{subequations}
By the exact two-dimensional elasticity complex \cite{SkyReissner,arnold_mixed_2002} there holds 
\begin{align}
    \exists \, \varsigma \in \C^\infty(\surf): \qquad \bm{\sigma} = \airy \varsigma \qquad \forall \,  \bm{\sigma} \in [\C^\infty(\surf)\otimes\Sym(2)] \cap \ker(\Di)   \, ,  
\end{align}
on contractible domains as per \cref{fig:elast}, where the airy operator reads
\begin{align}
    \airy \varsigma = \D^\perp \nabla^\perp \varsigma \, .
\end{align}
Consequently, we can insert the scalar potential $\varsigma$ into the two-dimensional Beltrami-Michell equations to find the biharmonic field equation 
\begin{align}
    \Delta \Delta \varsigma &= 0 && \text{in} && \surf \, ,
\end{align}
since $\tr (\airy \varsigma) = \Delta \varsigma$. In this context $\varsigma$ is called the Airy-stress function.

\section{Symmetrised stress formulation of isotropic linear elasticity}
The Beltrami-Michell equation in its classical form in \cref{eq:beltrami} does not lend itself to the construction of a symmetric bilinear form. The reason lies in the occurrence of the Hessian of the trace of the stress tensor $\hess (\tr \bm{\sigma})$, which does not yield a symmetric bilinear form after testing and partial integration. Fortunately, the field equation of static equilibrium can be reformulated into
\begin{align}
    -\dfrac{1}{1+\nu}(\di\Di \bm{\sigma}) \one = \dfrac{1}{1+\nu} (\di \vb{f}) \one \, ,
\end{align}
by taking the divergence of both sides and multiplying with the constant tensor $[1/(1+\nu)]\one$. Now, adding this zero-sum term to \cref{eq:beltrami} yields
\begin{align}
    -\Delta \bm{\sigma} - \dfrac{1}{1+\nu} [\hess( \tr \bm{\sigma} ) + (\di \Di \bm{\sigma}) \one ] =  2 \sym \D \vb{f} + \dfrac{1+\nu^2}{1-\nu^2} (\di \vb{f}) \one \, ,
    \label{eq:symBeltrami}
\end{align}
which does lend itself to a symmetric bilinear form. We apply the test function $\bm{\tau} \in \C^\infty(\overline{\body}) \otimes \Sym(3)$ to find
\begin{align}
    &\int_\body \con{\D \bm{\tau}}{ \D \bm{\sigma}} + \dfrac{1}{1+\nu} [\con{\Di \bm{\tau}}{\nabla \tr \bm{\sigma}} + \con{\nabla \tr \bm{\tau}}{\Di \bm{\sigma}}]  \, \dd \body = \int_\body 2\con{\bm{\tau}}{\sym\D \vb{f}} + \dfrac{1+\nu^2}{1-\nu^2} \con{\tr \bm{\tau}}{ \di \vb{f}} \, \dd \body \notag \\
    &\qquad + \int_{\partial \body} \con{\bm{\tau}}{(\D \bm{\sigma}) \vb{n}} + \dfrac{1}{1+\nu}[\con{\bm{\tau}}{(\nabla \tr \bm{\sigma}) \otimes \vb{n}} +  (\tr\bm{\tau})\con{\Di \bm{\sigma}}{\vb{n}} ]  \, \dd \surf \qquad \forall \, \bm{\tau} \in \C^\infty(\overline{\body})\otimes\Sym(3) \, ,   
\end{align}
by partial integration. Thus, we extract the Neumann boundary term
\begin{align}
    \bm{\kappa} = (\D \bm{\sigma}) \vb{n} + \dfrac{1}{1+\nu}[{(\nabla \tr \bm{\sigma}) \otimes \vb{n}} +  \con{\Di \bm{\sigma}}{\vb{n}} \one ] \, ,
\end{align}
by splitting the boundary between Dirichlet and Neumann $\partial \body = \surf_D \cup \surf_N$, such that $\surf_D \cap \surf_N = \emptyset$.
\begin{definition}[The pure stress boundary value problem for solids I]
\label{def:bvp3d}
    The field equations and boundary conditions of the novel three-dimensional stress boundary value problem read
    \begin{subequations}
        \begin{align}
            -\Delta \bm{\sigma} - \dfrac{1}{1+\nu} [\hess( \tr \bm{\sigma} ) + (\di \Di \bm{\sigma}) \one ]  &= 2 \sym \D \vb{f}   + \dfrac{1+\nu^2}{1-\nu^2} (\di \vb{f}) \one && \text{in} && \body \, ,  \\
            (\D \bm{\sigma})\vb{n} + \dfrac{1}{1+\nu} [\nabla \tr \bm{\sigma} \otimes \vb{n} - \con{\vb{f}}{\vb{n}}\one]  &= \bm{\kappa} && \text{on} && \surf_N \, , \\
         \bm{\sigma} &= \widetilde{\bm{\sigma}} && \text{on} && \surf_D  \, .
        \end{align}
    \end{subequations}
    where $\bm{\kappa}$ is a mixed measure of equilibrium and compatibility on the Neumann boundary where we applied $\Di \bm{\sigma} = - \vb{f}$.  
\end{definition}
We observe that continuous constant force fields $\vb{f} \in \R^3$ over the domain are captured by the Dirichlet and partially by the Neumann boundary conditions. We can now state the variational form.
\begin{definition}[Variational form of the pure stress problem I]
    \label{def:var3D}
    The weak formulation of the boundary value problem of linear elasticity written purely in stresses reads
    \begin{align}
    &\int_\body \con{\D \bm{\tau}}{ \D \bm{\sigma}} + \dfrac{1}{1+\nu} [\con{\Di \bm{\tau}}{\nabla \tr \bm{\sigma}} + \con{\nabla \tr \bm{\tau}}{\Di \bm{\sigma}}]  \, \dd \body = \int_{\surf_N} \con{\bm{\tau}}{\bm{\kappa}} \, \dd \surf \label{eq:weak3D} \\
    &\qquad + \int_\body 2\con{\bm{\tau}}{\sym\D \vb{f}} + \dfrac{1+\nu^2}{1-\nu^2} \con{\tr \bm{\tau}}{ \di \vb{f}} \, \dd \body \qquad \forall \, \bm{\tau} \in \C_{\surf_D}^\infty(\body) \otimes \Sym(3) \, , \notag
\end{align}
and yields a fully symmetric left-hand side. Here $\bm{\tau}$ is assumed to be compatible with the Dirichlet boundary.
\end{definition}
Clearly, it can be derived from the variation of a functional.
\begin{definition}[Variational functional in pure stresses I]
\label{def:energy3D}
    The weak formulation of the new boundary value problem can directly constructed as the variation of the functional
    \begin{align}
    I(\bm{\sigma}) &= \int_\body \dfrac{1}{2} \norm{\D \bm{\sigma}}^2 + \dfrac{1}{1+\nu} \con{\Di \bm{\sigma}}{\nabla \tr \bm{\sigma}}   -  2 \con{\bm{\sigma}}{\sym \D \vb{f}} - \dfrac{1+\nu^2}{1-\nu^2} \con{\tr \bm{\sigma}}{ \di \vb{f}}  \, \dd \body - \int_{\surf_N} \con{\bm{\sigma}}{\bm{\kappa}} \, \dd \surf  \, , 
\end{align}
    with respect to the stress tensor $\bm{\sigma}$.
\end{definition}

\subsection{Existence and uniqueness}
From the weak form in \cref{eq:weak3D} we extract the bilinear and linear forms
\begin{subequations}
\label{eq:biforms}
    \begin{align}
        a(\bm{\tau},\bm{\sigma}) &= \int_\body \con{\D \bm{\tau}}{ \D \bm{\sigma}} + \chi [\con{\Di \bm{\tau}}{\nabla \tr \bm{\sigma}} + \con{\nabla \tr \bm{\tau}}{\Di \bm{\sigma}}]  \, \dd \body \, , \\
        l(\bm{\tau}) &= \int_\body 2\con{\bm{\tau}}{\sym\D \vb{f}} + \dfrac{1+\nu^2}{1-\nu^2} \con{\tr \bm{\tau}}{ \di \vb{f}}  \, \dd \body \, ,
\end{align}
\end{subequations}
with the material parameter $\chi = 1/(1+\nu) \in [2/3, 1]$, given by the Poisson ratio $\nu \in [0, 1/2]$. In order to show existence and uniqueness we first derive some preliminary results. 
\begin{lemma}[Upper bound of $\norm{\D \bm{\sigma}}_\Le$]
\label{le:dS}
    Let $\bm{\sigma} \in \Honez(\body) \otimes \Sym(3)$ for all $\epsilon > 0$ there holds
    \begin{align}
        \norm{\D \bm{\sigma}}^2_{\Le} \leq  \max \left \{ 2 + \epsilon, 1 + \dfrac{1}{\epsilon} \right \} ( \norm{\Di \bm{\sigma} }^2_{\Le} +  \norm{\nabla \tr \bm{\sigma}}^2_{\Le} + \norm{\sym \Curl \bm{\sigma}}^2_{\Le} ) \qquad \forall \, \bm{\sigma} \in \Honez(\body) \otimes \Sym(3) \, .
    \end{align}
\end{lemma}
\begin{proof}
    Starting with \cref{eq:incid}, we reformulate the differential identity to
\begin{align}
     \Delta \bm{\sigma} = 2\sym (\D \Di \bm{\sigma}) - \hess (\tr \bm{\sigma}) + [ \Delta \tr \bm{\sigma}  - \di \Di \bm{\sigma}] \one - \inc \bm{\sigma} \, .
\end{align}
Now, testing the equation with $-\bm{\sigma} \in \Honez(\body)\otimes \Sym(3)$ and integrating by parts yields
\begin{align}
    \norm{\D \bm{\sigma}}^2_{\Le} &= 2 \norm{\Di \bm{\sigma} }^2_{\Le} - 2 \con{\Di \bm{\sigma}}{\nabla \tr \bm{\sigma}}_{\Le} + \norm{\nabla \tr \bm{\sigma}}^2_{\Le} + \con{\bm{\sigma}}{\inc \bm{\sigma}}_{\Le} \qquad \forall \, \bm{\sigma} \in \Honez(\body) \otimes \Sym(3) \, ,
    \label{eq:gradsig}
\end{align}
in the distributional sense.
For the incompatibility of the stress tensor $\inc \bm{\sigma} = \curl \Curl \bm{\sigma}$ we observe
\begin{align}
    \con{\bm{\sigma}}{\curl \Curl \bm{\sigma}}_{\Le} = \con{\curl \bm{\sigma}}{\Curl \bm{\sigma}}_{\Le} = \con{(\Curl \bm{\sigma})^T}{\Curl \bm{\sigma}}_{\Le} = \norm{\sym \Curl \bm{\sigma}}^2_{\Le} - \norm{\skw \Curl \bm{\sigma}}^2_{\Le} \, ,
\end{align}
where we algebraically decomposed $(\Curl \bm{\sigma})^T$ and $\Curl \bm{\sigma}$ into their symmetric and skew-symmetric parts. 
Using the latter we find the upper bound
\begin{align}
    \norm{\D \bm{\sigma}}^2_{\Le} &= 2 \norm{\Di \bm{\sigma} }^2_{\Le} - 2 \con{\Di \bm{\sigma}}{\nabla \tr \bm{\sigma}}_{\Le} + \norm{\nabla \tr \bm{\sigma}}^2_{\Le} + \norm{\sym \Curl \bm{\sigma}}^2_{\Le} - \norm{\skw \Curl \bm{\sigma}}^2_{\Le} 
    \notag \\
    &\overset{Y}{\leq} (2 + \epsilon) \norm{\Di \bm{\sigma} }^2_{\Le} + \left (1 + \dfrac{1}{\epsilon} \right ) \norm{\nabla \tr \bm{\sigma}}^2_{\Le} + \norm{\sym \Curl \bm{\sigma}}^2_{\Le} - \norm{\skw \Curl \bm{\sigma}}^2_{\Le}
     \\
    &\overset{\phantom{Y}}{\leq} \max \left \{ 2 + \epsilon, 1 + \dfrac{1}{\epsilon} \right \} ( \norm{\Di \bm{\sigma} }^2_{\Le} +  \norm{\nabla \tr \bm{\sigma}}^2_{\Le} + \norm{\sym \Curl \bm{\sigma}}^2_{\Le} ) \qquad \forall \, \bm{\sigma} \in \Honez(\body) \otimes \Sym(3) \, , \notag
\end{align}
where we applied Young's inequality\footnote{Young's inequality: $\con{x}{y}\leq \frac{\norm{x}^2}{2\epsilon}+\frac{\epsilon\norm{y}^2}{2} \quad \forall\, \epsilon > 0 \, .$} and dropped the negative term $- \norm{\skw \Curl \bm{\sigma}}^2_{\Le}$. 
\end{proof}
The properties of the $\sym \Curl$-operator along with its corresponding bounds can be found in \cite{Lewintan2021,LewintanInc,LewintanInc2,NEFF20151267,Neff2012,Gmeineder1,Gmeineder2,SkyNovel}.

We also require an upper bound for the skew symmetric curl of the stress tensor.
\begin{lemma}[Upper bound of $\norm{\skw \Curl \bm{\sigma}}$]
\label{le:skwCurls}
    Let $\bm{\sigma} \in \Honez(\body) \times \Sym(3)$, then there holds the upper bound
    \begin{align}
        \norm{\skw \Curl \bm{\sigma}}_{\Le}^2 \leq
        \dfrac{1 + \epsilon}{2} \norm{\Di \bm{\sigma}}_{\Le}^2 +  \dfrac{1 + \epsilon}{2\epsilon}\norm{\nabla \tr \bm{\sigma}}_{\Le}^2  \qquad \forall \, \bm{\sigma} \in \Hone(\body) \otimes \Sym(3) \, ,
    \end{align}
    for all constant $\epsilon > 0$.
\end{lemma} 
\begin{proof}
    The norm of the skew-symmetric part of the Curl of the stress tensor can be reformulated into
\begin{align}
    \norm{\skw \Curl \bm{\sigma}}_{\Le} &= \sup_{\bm{\tau} \in \Le \otimes \R^{3 \times 3}} \dfrac{\con{\bm{\tau}}{\skw \Curl \bm{\sigma}}_\Le}{\norm{\bm{\tau}}_\Le} \notag \\
    &= \sup_{\bm{\tau} \in \Le \otimes \R^{3 \times 3}} \dfrac{\con{ \Curl( \skw \bm{\tau})}{\bm{\sigma}}_{\H^{-1}\times \Hone}}{\norm{\bm{\tau}}_\Le}
    \notag \\
    &= \sup_{\bm{\tau} \in \Le \otimes \R^{3 \times 3}} \dfrac{\con{ (\di \axl \skw \bm{\tau}) \one - (\D \axl \skw \bm{\tau} )^T }{\bm{\sigma}}_{\H^{-1}\times \Hone}}{\norm{\bm{\tau}}_\Le}
    \notag \\
    &= \sup_{\bm{\tau} \in \Le \otimes \R^{3 \times 3}} \dfrac{\con{ (\di \axl \skw \bm{\tau}) \one - \sym(\D \axl \skw \bm{\tau} ) }{\bm{\sigma}}_{\H^{-1}\times \Hone}}{\norm{\bm{\tau}}_\Le}
    \notag \\
    &= \sup_{\bm{\tau} \in \Le \otimes \R^{3 \times 3}} \dfrac{\con{ \axl \skw \bm{\tau} }{\Di \bm{\sigma} - \nabla \tr \bm{\sigma}}_\Le}{\norm{\bm{\tau}}_\Le}
    \notag \\
    &= \sup_{\bm{\tau} \in \Le \otimes \R^{3 \times 3}} \dfrac{\con{ \bm{\tau} }{  \Anti( \Di \bm{\sigma} - \nabla \tr \bm{\sigma} )}_\Le}{ 2\norm{\bm{\tau}}_\Le }
    \notag \\
    &= \dfrac{1}{2}\norm{\Anti( \Di \bm{\sigma} - \nabla \tr \bm{\sigma} )}_\Le 
    \notag \\
    &= \dfrac{1}{\sqrt{2}} \norm{ \Di \bm{\sigma} - \nabla \tr \bm{\sigma} }_\Le \qquad \forall \, \bm{\sigma} \in \Honez(\body) \otimes \Sym(3) \, , 
\end{align}
where $\H^{-1}(\body)$ is the dual space of $\Honez(\body)$ and $\con{\cdot}{\cdot}_{\H^{-1} \times \Hone}$ is the dual product. In the derivation we used \cite{Lewintan2021} the algebraic identity $\Curl \bm{A} = (\di \axl \bm{A}) \one - (\D \axl \bm{A})^T$ and the relation $\con{\axl \bm{A}}{\axl \bm{A}} = (1/2) \con{\bm{A}}{\bm{A}}$ for some skew-symmetric tensor $\bm{A}$, as well as the symmetry of $\bm{\sigma}$. 
In fact, by the estimate we derive the following algebraic identity
\begin{align}
    \skw\Curl\bm{\sigma} = \dfrac{1}{2} \Anti( \Di \bm{\sigma} - \nabla \tr \bm{\sigma} ) \qquad \forall \,  \bm{\sigma} \in \C^\infty(\overline{\body}) \otimes \Sym(3) \, , 
\end{align}
which holds in general for symmetric tensors $\bm{\sigma} = \bm{\sigma}^T$. 
Employing the latter we find
\begin{align}
    \norm{\skw \Curl \bm{\sigma}}_{\Le}^2 &= \dfrac{1}{2}\norm{\Di \bm{\sigma} - \nabla \tr \bm{\sigma}}_\Le^2  
    \notag \\
    & \overset{T}{\leq} \dfrac{1}{2} (\norm{\Di \bm{\sigma}}_{\Le} + \norm{\nabla \tr \bm{\sigma}}_{\Le})^2
    \notag \\
    & \overset{\phantom{T}}{=} \dfrac{1}{2} (\norm{\Di \bm{\sigma}}_{\Le}^2 + 2\norm{\Di \bm{\sigma}}_{\Le}\norm{\nabla \tr \bm{\sigma}}_{\Le} + \norm{\nabla \tr \bm{\sigma}}_{\Le}^2)
    \notag \\
    & \overset{Y}{\leq} \dfrac{1 + \epsilon}{2} \norm{\Di \bm{\sigma}}_{\Le}^2 +  \dfrac{1 + \epsilon}{2\epsilon}\norm{\nabla \tr \bm{\sigma}}_{\Le}^2  \qquad \forall \, \bm{\sigma} \in \Hone(\body) \otimes \Sym(3) \, , 
\end{align}
where we used the triangle\footnote{Triangle inequality: $ \norm{x + y} \leq \norm{x} + \norm{y} \, .$} inequality and Young's inequality.
\end{proof}
With the latter bounds at hand we can show well-posedness of the three-dimensional problem.
\begin{theorem}[Existence and uniqueness in three dimensions]
\label{th:wellposed3D}
    Let the Poisson ratio be $\nu \in [0,1/2]$ such that $\chi = 1/(1+\nu) \in [2/3, 1]$, then the problem
    \begin{align}
        a(\bm{\tau},\bm{\sigma}) = l(\bm{\tau}) \qquad \forall \, \bm{\tau} \in \X_0(\body) \, , 
    \end{align}
    where the bilinear and linear forms are according to \cref{eq:biforms}, and the space $\X_0(\body)$ is given by
    \begin{align}
        &\X_0(\body) = \Honez(\body) \otimes \Sym(3) \, , && \norm{\bm{\sigma}}_{\X}^2 = \norm{\bm{\sigma}}_{\Hone}^2 = \norm{\bm{\sigma}}^2_{\Le} + \norm{\D\bm{\sigma}}^2_{\Le} \, ,
    \end{align}
    has a unique solution $\bm{\sigma} \in \X_0(\body)$ for every right-hand side with the stability estimate 
    \begin{align}
        \norm{\bm{\sigma}}_{\X} \leq \dfrac{1}{\beta} \norm{l}_{\X'} \, ,
    \end{align}
    where $\beta = \beta(\nu) > 0$ is the coercivity constant and $\norm{\cdot}_{\X'}$ is the dual norm with respect to $\X_0(\body)$. 
\end{theorem}
\begin{proof}
    The proof follows by the Lax-Milgram theorem, where we neglect Neumann boundary conditions. The continuity of the linear form is obvious. For the bilinear form we first observe that there holds
    \begin{align}
        \norm{\nabla \tr \bm{\sigma}}_{\Le} &= \norm{\Tr (\nabla \otimes \bm{\sigma})}_{\Le} \overset{CS}{\leq} \norm{\Tr}_*\norm{\nabla \otimes \bm{\sigma}}_{\Le} = \sqrt{3}\norm{\D \bm{\sigma}}_{\Le} \, ,  \notag \\ 
        \norm{\Di \bm{\sigma}}_{\Le} &= \norm{\Tr \D \bm{\sigma}}_{\Le} \overset{CS}{\leq} \norm{\Tr}_*\norm{\D \bm{\sigma}}_{\Le} = \sqrt{3}\norm{\D \bm{\sigma}}_{\Le} \, , 
        \label{eq:divtrdSig}
    \end{align}
    which follows from the Cauchy-Schwarz inequality\footnote{Cauchy-Schwarz inequality: $\con{x}{y}\leq \norm{x}\norm{y} \, .$} with the operator norm $\norm{\Tr}_* = \norm{\one} = \sqrt{3}$. Thus, we find
    \begin{align}
        a(\bm{\tau},\bm{\sigma}) &=  \con{\D \bm{\tau}}{ \D \bm{\sigma}}_\Le + \chi [\con{\Di \bm{\tau}}{\nabla \tr \bm{\sigma}}_\Le + \con{\nabla \tr \bm{\tau}}{\Di \bm{\sigma}}_\Le]  \notag \\
        & \overset{CS}{\leq} \norm{\D \bm{\tau}}_\Le \norm{\D \bm{\sigma}}_\Le + \chi (\norm{\Di \bm{\tau}}_\Le \norm{\nabla \tr \bm{\sigma}}_\Le +  \norm{\nabla \tr \bm{\tau}}_\Le \norm{\Di \bm{\sigma}}_\Le)   \notag \\
        & \overset{*}{\leq} \norm{\D \bm{\tau}}_\Le \norm{\D \bm{\sigma}}_\Le + 6 \chi \norm{\D \bm{\tau}}_\Le \norm{\D \bm{\sigma}}_\Le  \notag \\
        & \overset{\phantom{T}}{\leq} (1 + 6\chi ) \norm{\bm{\tau}}_{\X} \norm{\bm{\sigma}}_{\X} \qquad \forall \, \bm{\tau},\bm{\sigma} \in \X(\body) \, , 
    \end{align}
    where we applied the Cauchy-Schwarz inequality and the operator bounds. In order to show coercivity we reformulate \cref{eq:gradsig} into
    \begin{align}
         2 \con{\Di \bm{\sigma}}{\nabla \tr \bm{\sigma}}_{\Le} &= 2 \norm{\Di \bm{\sigma} }^2_{\Le} - \norm{\D \bm{\sigma}}^2_{\Le} + \norm{\nabla \tr \bm{\sigma}}^2_{\Le} + \norm{\sym \Curl \bm{\sigma}}_{\Le}^2 - \norm{\skw \Curl \bm{\sigma}}_{\Le}^2 \qquad \forall \, \bm{\sigma} \in \X_0(\body) \, .
    \end{align}
    Now, for the bilinear form we find
    \begin{align}
        a(\bm{\sigma},\bm{\sigma}) &= \norm{\D \bm{\sigma}}_\Le^2 + 2\chi \con{\Di \bm{\sigma}}{\nabla \tr \bm{\sigma}}_\Le  \notag \\
        &\overset{\phantom{T}}{=} (1 - \chi) \norm{\D \bm{\sigma}}_\Le^2 + 2\chi\norm{\Di \bm{\sigma}}_\Le^2 + \chi \norm{\nabla \tr \bm{\sigma}}_{\Le}^2 + \chi\norm{\sym \Curl \bm{\sigma}}_{\Le}^2 - \chi\norm{\skw \Curl \bm{\sigma}}_{\Le}^2 
        \notag \\
        &\overset{*}{\geq} (1 - \chi) \norm{\D \bm{\sigma}}_\Le^2 +   \dfrac{(3 -  \epsilon)\chi}{2} \norm{\Di \bm{\sigma}}_\Le^2 + \dfrac{(\epsilon-1)\chi}{2\epsilon} \norm{\nabla \tr \bm{\sigma}}_{\Le}^2 + \chi\norm{\sym \Curl \bm{\sigma}}_{\Le}^2  
        \notag \\
        &\overset{\phantom{T}}{\geq} (1 - \chi) \norm{\D \bm{\sigma}}_\Le^2 + \chi \min \left \{  \dfrac{3 -  \epsilon}{2}, 1, \dfrac{\epsilon-1}{2\epsilon} \right \}  (\norm{\Di \bm{\sigma}}_\Le^2  + \norm{\sym \Curl \bm{\sigma}}_{\Le}^2 +  \norm{\nabla \tr \bm{\sigma}}_{\Le}^2 ) 
        \notag \\
        &\overset{\star}{\geq} (1 - \chi) \norm{\D \bm{\sigma}}_\Le^2 + \dfrac{\chi}{3} \min \left \{  \dfrac{3 -  \epsilon}{2}, 1, \dfrac{\epsilon-1}{2\epsilon} \right \}  \norm{\D \bm{\sigma}}_\Le^2     
        \notag \\
        &\overset{PF}{\geq} \dfrac{1}{1+c_F^2}  \left [ (1 - \chi ) + \dfrac{\chi}{3} \min \left \{  \dfrac{3 -  \epsilon}{2}, 1, \dfrac{\epsilon-1}{2\epsilon} \right \} \right ] \norm{\bm{\sigma}}_{\X}^2  \qquad \forall \, \bm{\sigma} \in \X_0(\body) \, ,
    \end{align}
    where we used the upper bound of $\norm{\skw \Curl \bm{\sigma}}_\Le^2$ from \cref{le:skwCurls} in $*$, the upper bound of $\norm{\D \bm{\sigma}}_\Le^2$ from \cref{le:dS} in $\star$ with $\epsilon = 1$, and finally applied the 
    Poincar\'e-Friedrich inequality\footnote{Poincar\'e-Friedrich inequality: $\exists \, c_F > 0 : \quad \norm{x}_\Le \leq c_F\norm{\nabla x}_\Le \quad \forall \, x \in \Honez(\body) \, ,$ such that $c_F$ depends only on the domain $\body \subseteq \R^3$ and its boundary $\partial \body$.}.
    The inequality holds for the complete range of $\chi \in [2/3, 1]$ with $1 \leq \epsilon \leq 3$, since even for $\chi = 1$ we can set $1 < \epsilon < 3$, such that the gradient of the stress tensor $\norm{\D \bm{\sigma}}_\Le^2$ does not vanish.  
\end{proof}
\begin{corollary}[Well-posedness of the three-dimensional Beltrami-Michell equations]
    From the proof of \cref{th:wellposed3D} we automatically get that the variational formulation of the original Beltrami-Michell equations in \cref{eq:beltrami} is also well-posed for $\bm{\sigma} \in \X_0(\body)$. The bilinear form is simply
    \begin{align}
        a(\bm{\tau},\bm{\sigma}) &= \int_\body \con{\D \bm{\tau}}{ \D \bm{\sigma}} + \chi \con{\Di \bm{\tau}}{\nabla \tr \bm{\sigma}}  \, \dd \body \, ,
    \end{align}
    such that continuity and coercivity follow the same lines of our proof simply with $\chi$ instead of $2\chi$. However, this bilinear form is non-symmetric, making computational approaches to its solution far less efficient. 
\end{corollary}

\section{Novel stress formulation of planar isotropic linear elasticity}
From \cref{eq:planestress} and \cref{eq:planestrain} it is clear that the planar Beltrami-Michell equations characterise only the mean stress $\sigma_m = (1/2) \tr \bm{\sigma}$, related to the first invariant of the stress tensor. However, we can derive from static equilibrium
\begin{align}
    -\psi\sym(\D \Di \bm{\sigma}) = \psi\sym\D\vb{f} \, ,
\end{align}
with some constant $\psi > 0$.
Now, multiplying either \cref{eq:planestress} or \cref{eq:planestrain} with the two-dimensional identity tensor $\one$ and adding the symmetrised gradient of static equilibrium to it yields the field equations
\begin{subequations}
\label{eq:planarmod}
    \begin{align}
      -\psi\sym(\D \Di \bm{\sigma}) -  \dfrac{1}{1+\nu}(\Delta \tr \bm{\sigma}) \one &= \psi\sym\D \vb{f} + (\di \vb{f}) \one \, , \\
      -\psi\sym(\D \Di \bm{\sigma}) - [(1-\nu) \Delta \tr \bm{\sigma}] \one &= \psi\sym\D \vb{f} + (\di \vb{f}) \one \, .  
\end{align}
\end{subequations}
The right-hand side now clearly vanishes for body forces $\vb{f}$ in the space of rigid body motions
\begin{align}
    \RM(\surf) &= \R^2 \oplus \bm{R}\vb{x} = \ker(\sym\D) \cap \ker(\di) \subset \ker(\di) \, ,  && \bm{R} = \begin{bmatrix}
        0 & 1 \\
        -1 & 0 
    \end{bmatrix} \, , && \vb{x} = x \vb{e}_1 + x \vb{e}_2 \, ,
\end{align}
analogously to the three-dimensional case.
By respectively defining the constant $\chi = 1/(1+\nu)$ or $\chi = 1-\nu$ for plane stress or plane strain, we can denote either problem with the single equation
\begin{align}
    -\psi \sym (\D \Di \bm{\sigma}) - (\chi \Delta \tr \bm{\sigma} ) \one = \psi \sym \D \vb{f} + (\di \vb{f}) \one \, .
    \label{eq:planarStrong}
\end{align}
We apply the test function $\bm{\tau} \in \C^\infty(\overline{\surf}) \otimes \Sym(2)$ and integrate by parts to find
\begin{align}
    &\int_\surf \psi \con{\Di \bm{\tau}}{ \Di \bm{\sigma} } + \chi \con{\nabla \tr \bm{\tau}}{ \nabla \tr \bm{\sigma} } \, \dd \surf 
    \notag \\
    & \qquad = \int_{\partial \surf} \psi \con{\bm{\tau}}{\Di \bm{\sigma}\otimes \vb{n}} + \chi \tr \bm{\tau} \con{\nabla \tr \bm{\sigma}}{\vb{n}} \, \dd \curv
    \\
    & \qquad \qquad + \int_\surf \psi \con{\bm{\tau}}{\sym \D \vb{f}} + \con{\tr \bm{\tau}}{\di \vb{f}} \, \dd \surf \qquad \forall \, \bm{\tau} \in \C^\infty(\overline{\surf})\otimes\Sym(2) \, .
    \notag 
\end{align}
We split the boundary between Dirichlet and Neumann $\partial \surf = \curv_D \cup \curv_N$ with $\curv_D \cap \curv_N = \emptyset$ to finally derive the new equations and boundary conditions, where we apply $\Di \bm{\sigma} = -\vb{f}$ for the Neumann boundary.
\begin{definition}[Pure stress planar boundary value problem I]
\label{def:planarBVP}
    Let $\chi = 1/(1+\nu)$ or $\chi = 1-\nu$ for plane stress or plane strain, respectively,  
    the complete boundary value problem reads
    \begin{subequations}
        \begin{align}
        -\psi \sym(\D \Di \bm{\sigma}) - (\chi \Delta \tr \bm{\sigma})\one &= \psi \sym \D \vb{f} + (\di \vb{f})\one && \text{in} && \surf \, , \\ 
        \chi \con{\nabla \tr \bm{\sigma}}{\vb{n}} \one - \psi (\vb{f} \otimes \vb{n}) &= \bm{\kappa}  && \text{on} && \curv_N \, , \\
        \bm{\sigma} &= \widetilde{\bm{\sigma}}  && \text{on} && \curv_D \, .
    \end{align}
    \end{subequations}
\end{definition}
With the boundary conditions at hand, we can construct the variational form of \cref{def:planarBVP}.
\begin{definition}[Pure stress planar variational form I]
\label{def:var2D}
    The variational form reads 
    \begin{align}
    \int_\surf \psi \con{\Di \bm{\tau}}{ \Di \bm{\sigma} } + \chi \con{\nabla \tr \bm{\tau}}{ \nabla \tr \bm{\sigma} } \, \dd \surf 
     &= \int_\surf \psi \con{\bm{\tau}}{\sym \D \vb{f}} + \con{\tr \bm{\tau}}{\di \vb{f}}  \, \dd \surf \label{eq:weak2D}
     \\& \qquad + \int_{\curv_N} \con{\bm{\tau}}{\bm{\kappa}} \, \dd \curv   \qquad \forall \, \bm{\tau} \in \C_{\curv_D}^\infty(\overline{\surf})\otimes\Sym(2) \, ,  \notag
\end{align}
    with $\psi > 0$ and $\chi = 1/(1+\nu)$ or $\chi = 1-\nu$, for plane-stress or plane stain, respectively.
\end{definition}
The latter can clearly be derived from a variation functional.
\begin{definition}[Planar variation functional in pure stresses I]
\label{def:energy2D}
    The weak formulation of the new boundary value problem can directly constructed as the variation of the functional 
    \begin{align}
    I(\bm{\sigma}) &= \int_\surf \dfrac{\psi}{2} \norm{\Di \bm{\sigma}}^2 + \dfrac{\chi}{2}\norm{\nabla \tr \bm{\sigma}}^2  - \psi \con{\bm{\tau}}{ \sym \D \vb{f}} - \con{\tr \bm{\tau}}{\di \vb{f}} \, \dd \surf - \int_{\curv_N} \con{\bm{\sigma}}{\bm{\tau}} \, \dd \curv  \, , 
\end{align}
    with $\psi > 0$ and $\chi = 1/(1+\nu)$ or $\chi = 1-\nu$, for plane-stress or plane stain, respectively.
\end{definition}

In the planar case we have $\bm{\sigma}:\overline{\surf} \subset \R^2 \to \Sym(2)$, such that the stress tensor is a function of the $x-y$-plane $\bm{\sigma} = \bm{\sigma}(x,y)$. Consequently, we find  
\begin{align}
    \inc \widehat{\bm{\sigma}} = (\Delta \tr \widehat{\bm{\sigma}} - \di \Di \widehat{\bm{\sigma}} ) \vb{e}_3 \otimes \vb{e}_3 \, ,
\end{align}
where $\widehat{\bm{\sigma}} \in \C^\infty(\overline{\surf}) \otimes \Sym(3)$ is the three dimensional insertion of the two-dimensional tensor $\bm{\sigma} \in \C^\infty(\overline{\surf}) \otimes \Sym(2)$ into a three-dimensional tensor $\widehat{\sigma}_{\alpha \beta} = \sigma_{\alpha \beta}$ with zeroes in the out-of-plane direction $\widehat{\sigma}_{i3} = \widehat{\sigma}_{3i} = 0$. Due to the latter \cref{eq:laplacesig} reduces to
\begin{align}
     \Delta \bm{\sigma} = 2\sym (\D \Di \bm{\sigma}) - \hess (\tr \bm{\sigma}) + [ \Delta \tr \bm{\sigma}  - \di \Di \bm{\sigma}] \one \, ,
     \label{eq:laplaceSigPlanar}
\end{align}
being a new identity for the Laplacian of a symmetric two-dimensional tensor field. We can reformulate the identity into 
\begin{align}
    -\sym (\D \Di \bm{\sigma}) - (\Delta \tr \bm{\sigma})\one = \sym (\D \Di \bm{\sigma}) - \Delta \bm{\sigma} - \hess (\tr \bm{\sigma}) -  (\di \Di \bm{\sigma}) \one \, .  
\end{align}
Now, let $\psi = \chi$, then we find for the left-hand side of \cref{eq:planarStrong}
\begin{align}
    -\chi \sym (\D \Di \bm{\sigma}) - (\chi \Delta \tr \bm{\sigma} ) \one &= \chi[\sym (\D \Di \bm{\sigma}) - \Delta \bm{\sigma} - \hess (\tr \bm{\sigma}) -  (\di \Di \bm{\sigma}) \one] 
    \notag \\
    &= \chi[-\sym \D \vb{f} - \Delta \bm{\sigma} - \hess (\tr \bm{\sigma}) -  (\di \Di \bm{\sigma}) \one] \, . 
\end{align}
Thus, we derive a new equivalent form for the planar equations 
\begin{align}
    - \Delta \bm{\sigma} - \hess(\tr \bm{\sigma}) - (\di \Di \bm{\sigma}) \one = 2\sym \D \vb{f} + \dfrac{1}{\chi}(\di \vb{f})\one \, ,
    \label{eq:2deq}
\end{align}
by dividing the entire equation by $\chi$. Clearly, the new form is extremely similar to our symmetrised formulation in \cref{eq:symBeltrami} of the three-dimensional problem. However, it leads to a more involved bilinear form.  

\subsection{Existence and uniqueness in the planar case}
We extract the bilinear and linear forms of the planar problem from \cref{eq:weak2D}
\begin{subequations}
\label{eq:biform2d}
    \begin{align}
        a(\bm{\tau},\bm{\sigma}) &= \int_\surf \psi  \con{\Di \bm{\tau}}{ \Di \bm{\sigma} } + \chi \con{\nabla \tr \bm{\tau}}{ \nabla \tr \bm{\sigma} } \, \dd \surf \, , \\
     l(\bm{\tau}) &=  \int_\surf \psi \con{\bm{\tau}}{\sym \D \vb{f}} + \con{\tr \bm{\tau}}{\di \vb{f}}  \, \dd \surf \, . 
    \end{align}
\end{subequations}
As a preparatory step we introduce a new norm equivalence.
\begin{lemma}[Norm equivalence for $\norm{\D \bm{\sigma}}_{\Le}$]
\label{le:normeq}
    Let $\bm{\sigma} \in \Honez(\surf) \otimes \Sym(2)$, then there holds
    \begin{align}
        \exists \, c_1, c_2 > 0: \quad c_1 (\norm{\Di \bm{\sigma}}^2_{\Le} + \norm{\nabla \tr \bm{\sigma}}^2_{\Le}) \leq \norm{\D \bm{\sigma}}_{\Le}^2 \leq c_2 (\norm{\Di \bm{\sigma}}^2_{\Le} + \norm{\nabla \tr \bm{\sigma}}^2_{\Le}) \quad \forall \, \bm{\sigma} \in \Honez(\surf) \otimes \Sym(2) \, ,
    \end{align}
    with 
    \begin{align}
        &c_1 = \min \left \{ 2 - \epsilon, 1 - \dfrac{1}{\epsilon} \right \} \, , &&
        c_2 = \max \left \{ 2 + \epsilon, 1 + \dfrac{1}{\epsilon} \right \} \, ,
    \end{align}
    where $\epsilon > 0$.
\end{lemma}
\begin{proof}
    Starting with \cref{eq:laplaceSigPlanar}, we test with $-\bm{\sigma} \in \Honez(\surf) \otimes \Sym(2)$ and integrate by parts to find
    \begin{align}
        \norm{\D \bm{\sigma}}^2_{\Le} = 2 \norm{\Di \bm{\sigma}}^2_{\Le} - 2 \con{\Di \bm{\sigma}}{\nabla \tr \bm{\sigma}}_{\Le} + \norm{\nabla \tr \bm{\sigma}}^2_{\Le} \qquad \forall \, \bm{\sigma} \in \Honez(\surf) \otimes \Sym(2) \, . 
        \label{eq:con2Deq}
    \end{align}
    Now, applying Young's inequality to find an upper bound yields
    \begin{align}
        \norm{\D \bm{\sigma}}^2_{\Le} &\overset{Y}{\leq} (2 + \epsilon) \norm{\Di \bm{\sigma}}^2_{\Le} + \left(1 + \dfrac{1}{\epsilon} \right ) \norm{\nabla \tr \bm{\sigma}}^2_{\Le} \notag \\
        &\overset{\phantom{Y}}{\leq} \max \left \{ 2+ \epsilon, 1 + \dfrac{1}{\epsilon} \right \}  (\norm{\Di \bm{\sigma}}^2_{\Le} + \norm{\nabla \tr \bm{\sigma}}^2_{\Le}) \qquad \forall \, \bm{\sigma} \in \Honez(\surf) \otimes \Sym(2) \, ,
    \end{align}
    which is satisfied for any $\epsilon > 0$. Alternatively, we find the lower bound via the negative Young's inequality
    \begin{align}
        \norm{\D \bm{\sigma}}^2_{\Le} &\overset{Y}{\geq} (2 - \epsilon) \norm{\Di \bm{\sigma}}^2_{\Le} + \left(1 - \dfrac{1}{\epsilon} \right ) \norm{\nabla \tr \bm{\sigma}}^2_{\Le} \notag \\
        &\overset{\phantom{Y}}{\geq} \min \left \{ 2 - \epsilon, 1 - \dfrac{1}{\epsilon} \right \}  (\norm{\Di \bm{\sigma}}^2_{\Le} + \norm{\nabla \tr \bm{\sigma}}^2_{\Le}) \qquad \forall \, \bm{\sigma} \in \Honez(\surf) \otimes \Sym(2) \, ,
    \end{align}
    which holds for $1 < \epsilon < 2$.
\end{proof}
\begin{corollary}[Norm equivalence for $\norm{\bm{\sigma}}_{\Hone}$]
    From \cref{le:normeq} we automatically derive
    \begin{align}
        c_1 (\norm{\Di \bm{\sigma}}^2_{\Le} + \norm{\nabla \tr \bm{\sigma}}^2_{\Le}) \leq \norm{\bm{\sigma}}_{\Hone}^2 \leq c_2(1+ c_F^2) (\norm{\Di \bm{\sigma}}^2_{\Le} + \norm{\nabla \tr \bm{\sigma}}^2_{\Le}) \quad \forall \, \bm{\sigma} \in \Honez(\surf) \otimes \Sym(2) \, ,
    \end{align}
    by using $\norm{\D \bm{\sigma}}^2_{\Le} \leq \norm{\bm{\sigma}}^2_{\Le} + \norm{\D \bm{\sigma}}^2_{\Le}$ and the Poincar\'e-Friedrich inequality.
\end{corollary}
With the norm equivalence at hand, we can directly show well-posedness.
\begin{theorem}[Well-posedness of the planar forms] \label{th:planar}
    Let the Poisson ratio be $\nu \in [0,1/2]$ such that $\chi \in [2/3, 1]$ or $\chi \in [1/2, 1]$ and $\psi > 0$, then the problem
    \begin{align}
        a(\bm{\tau},\bm{\sigma}) = l(\bm{\tau}) \qquad \forall \, \bm{\tau} \in \Y_0(\surf) = \Honez(\surf) \otimes \Sym(2) \, ,  
    \end{align}
    where the bilinear and linear forms are given in \cref{eq:biform2d}, has a unique solution $\bm{\sigma} \in \Y_0(\surf)$ for every right-hand side with the stability estimate 
    \begin{align}
        \norm{\bm{\sigma}}_{\Y} \leq \dfrac{1}{\beta} \norm{\vb{f}}_{\Y'} \, ,
    \end{align}
    where $\beta = \beta(\psi,\chi) > 0$ is the coercivity constant and $\norm{\cdot}_{\Y'}$ is the dual norm with respect to $\Y_0(\surf)$.  
\end{theorem}
\begin{proof}
    We assume a vanishing Neumann boundary. The continuity of the linear form is obvious. For the bilinear form we first observe 
    \begin{align}
        \norm{\nabla \tr \bm{\sigma}}_{\Le} &= \norm{\Tr  (\nabla \otimes \bm{\sigma})}_{\Le} \overset{CS}{\leq} \norm{\Tr}_*\norm{\nabla \otimes  \bm{\sigma}}_{\Le} = \sqrt{2}\norm{\D \bm{\sigma}}_{\Le} \, ,  \notag \\ 
        \norm{\Di \bm{\sigma}}_{\Le} &= \norm{\Tr \D \bm{\sigma}}_{\Le} \overset{CS}{\leq} \norm{\Tr}_*\norm{\D \bm{\sigma}}_{\Le} = \sqrt{2}\norm{\D \bm{\sigma}}_{\Le} \, ,
        \label{eq:tr2dap}
    \end{align}
    which follows by the Cauchy-Schwarz inequality with $\norm{\Tr}_* = \norm{\one} = \sqrt{2}$ since $\one \in \R^{2 \times 2}$ in the two-dimensional case.
    Now, the continuity of bilinear form is given due to
    \begin{align}
        a(\bm{\tau},\bm{\sigma}) &=  \psi  \con{\Di \bm{\tau}}{ \Di \bm{\sigma} }_{\Le} + \chi \con{\nabla \tr \bm{\tau}}{ \nabla \tr \bm{\sigma} }_{\Le} \notag \\
        &\overset{CS}{\leq}  \psi \norm{\Di \bm{\tau}}_{\Le}\norm{ \Di \bm{\sigma} }_{\Le} + \chi \norm{\nabla \tr \bm{\tau}}_{\Le}\norm{ \nabla \tr \bm{\sigma} }_{\Le} \notag \\
        &\overset{*}{\leq}  2\psi \norm{\D \bm{\tau}}_{\Le}\norm{ \D \bm{\sigma} }_{\Le} + 2\chi \norm{\D  \bm{\tau}}_{\Le}\norm{ \D \bm{\sigma} }_{\Le} \notag \\
        &\overset{\phantom{Y}}{\leq}  2(\psi  + \chi) \norm{\D \bm{\tau}}_{\Y}\norm{ \D \bm{\sigma} }_{\Y} \qquad \forall \, \bm{\sigma} \in \Y_0(\surf) \, , 
    \end{align}
    where we used the Cauchy-Schwarz inequality and \cref{eq:tr2dap} in $*$. 
    For the coercivity we find
    \begin{align}
        a(\bm{\sigma},\bm{\sigma}) &= \psi  \norm{\Di \bm{\sigma}}_\Le^2 + \chi \norm{\nabla \tr \bm{\sigma}}_{\Le}^2  \notag \\
        &\overset{\phantom{Y}}{\geq} \min\{\psi, \chi\} ( \norm{\Di \bm{\sigma}}_\Le^2 +  \norm{\nabla \tr \bm{\sigma}}_{\Le}^2) \notag \\
        &\overset{*}{\geq} \dfrac{1}{3}\min\{\psi, \chi\}  \norm{\D \bm{\sigma}}_\Le^2 \notag \\
        &\overset{PF}{\geq} \dfrac{1}{3(1+c_F^2)}\min\{\psi, \chi\}  \norm{ \bm{\sigma}}_\Y^2 \qquad \forall \, \bm{\sigma} \in \Y_0(\surf) \, , 
    \end{align}
    where we used the norm equivalence from \cref{le:normeq} with $\epsilon = 1$ in $*$, and then applied the Poincar\'e-Friedrich inequality.
\end{proof}
\begin{remark}[Vanishing $\psi$-constant]
    Note that the proof fails for $\psi = 0$, emphasising that the original planar Beltrami-Michell equation is insufficient in order to fully characterise the stress tensor $\bm{\sigma}$. 
\end{remark}
\begin{corollary}[Well-posedness of the alternative planar form] \label{co:well2d}
    For $\psi = \chi$ the bilinear form given by partial integration of \cref{eq:2deq} 
    \begin{align}
        a(\bm{\tau},\bm{\sigma}) = \int_\surf \con{\D \bm{\tau}}{ \D \bm{\sigma}} + \con{\Di \bm{\tau}}{\nabla \tr \bm{\sigma}} + \con{\nabla \tr \bm{\tau}}{\Di \bm{\sigma}} \, \dd \surf \, ,
    \end{align}
    is coercive. The proof follows the same lines as in the proof of \cref{th:planar}, since for coercivity we have
    \begin{align}
        a(\bm{\sigma},\bm{\sigma}) = \norm{\D \bm{\sigma}}^2_{\Le} + 2 \con{\Di \bm{\sigma}}{ \nabla \tr \bm{\sigma} }_{\Le} = 2\norm{\Di \bm{\sigma}}^2_{\Le} + \norm{\nabla \tr \bm{\sigma}}^2_{\Le} \qquad \forall \, \bm{\sigma} \in \Y_0(\surf) \, ,
    \end{align}
    due to \cref{eq:con2Deq}.
\end{corollary}

\section{Discrete variational formulations}
The well-posedness of the three-dimensional and planar formulations is proven by the Lax-Milgram theorem, such that any conforming discretisation is automatically well-posed.
Consequently, we directly obtain the quasi-best approximation via Cea's lemma.
\begin{theorem}[Quasi-best approximation] \label{th:cea}
    Let $\bm{\sigma} \in {\X(\body)}$ or $\bm{\sigma} \in \Y(\surf)$ be the exact solution and $\bm{\sigma}^h \in \X^h(\body)\subset \X(\body)$ or $\bm{\sigma}^h \in \Y^h(\body)\subset \Y(\body)$ be the approximate solution, respectively, then there holds
    \begin{align}
         &\norm{\bm{\sigma} - \bm{\sigma}^h}_\X \leq \dfrac{\alpha}{\beta} \inf_{\bm{\tau}\in\X^h} \norm{\bm{\sigma}-\bm{\tau}}_{\X} \, , && \norm{\bm{\sigma} - \bm{\sigma}^h}_\Y \leq \dfrac{\alpha}{\beta} \inf_{\bm{\tau}\in\Y^h} \norm{\bm{\sigma}-\bm{\tau}}_{\Y} \, ,
    \end{align}
    where $\alpha$ and $\beta$ are the continuity and coercivity constants, respectively.
\end{theorem}
Since Cea's lemma is satisfied, the a priori error estimates follow from standard interpolation errors.
\begin{theorem}[Convergence estimates]
    Let $\bm{\sigma} \in \X(\body)$ or $\bm{\sigma} \in \Y(\surf)$ be the exact solution and $\bm{\sigma}^h \in \X^h(\body)\subset \X(\body)$  or $\bm{\sigma}^h \in \Y^h(\body)\subset \Y(\body)$ be its approximation, respectively, there holds the a priori error estimate
    \begin{align}
        &\norm{\bm{\sigma}-\bm{\sigma}^h}_{\X} \leq c h^p |\bm{\sigma}|_{\H^{p+1}} \, , && 
        \norm{\bm{\sigma}-\bm{\sigma}^h}_{\Y} \leq c h^p |\bm{\sigma}|_{\H^{p+1}} \, ,
    \end{align}
    where the constant $c$ is independent of the element size and the exact solution $c \neq c(h,\bm{\sigma})$, and $|\cdot|_{\H^{p+1}}$ is the Sobolev semi-norm.
\end{theorem}
\begin{proof}
    The proof follows directly by replacing $\bm{\tau} \in \X^h(\body)$ or $\bm{\tau} \in \Y^h(\body)$ in \cref{th:cea} with the interpolation $\Pi_g^p\bm{\sigma} \in \X^h(\body)$ or $\Pi_g^p\bm{\sigma} \in \Y^h(\body)$, respectively, and classical interpolation error estimates \cite{Zaglmayr2006,Demkowicz2000}.
\end{proof}

We conclude our numerical treatment with an observation with respect to discontinuous body forces under the assertion that $\C^0$-continuous Lagrange-type elements are employed in the discretisation.
\begin{observation}[Discontinuous forces]
\label{obs:forces}
    In both the three-dimensional problem and the planar problems we can understand scalar products with derivatives of $\vb{f}$ in the form of distributions 
    \begin{subequations}
        \begin{align}
        \con{\bm{\tau}}{\sym \D \vb{f}}_{\mathcal{T}} &= \sum_{T \in \mathcal{T}}\con{\bm{\tau}}{\vb{f} \otimes \vb{n}}_{\partial T} - \con{\Di\bm{\tau}}{\vb{f}}_{T} \, ,  \\
        \con{\tr\bm{\tau}}{\di  \vb{f}}_{\mathcal{T}} &= \sum_{T \in \mathcal{T}}(\tr\bm{\tau})\con{\vb{f}}{ \vb{n}}_{\partial T} - \con{\nabla \tr\bm{\tau}}{\vb{f}}_{T} \, ,
    \end{align}
    \end{subequations}
    where $\mathcal{T}$ represents the triangulation, $T$ is an element in the mesh and $\partial T$ is its respective boundary. The term on the boundary of the elements vanishes if the jump in the force field does not occur directly on the interface between two elements, since $\bm{\tau}$ is continuous. Consequently, the distributive scalar products read    
    \begin{align}
        &\con{\bm{\tau}}{\sym \D \vb{f}}_{\mathcal{T}} = \con{\bm{\tau}}{\vb{f} \otimes \vb{n}}_{\partial_N} -\sum_{T \in \mathcal{T}} \con{\Di\bm{\tau}}{\vb{f}}_{T} \, ,  &&
        \con{\tr\bm{\tau}}{\di  \vb{f}}_{\mathcal{T}} = \con{\vb{f}}{ \vb{n}}_{\partial_N} -\sum_{T \in \mathcal{T}} \con{\nabla \tr\bm{\tau}}{\vb{f}}_{T} \, ,
    \end{align}
    \textbf{In practice, these terms are simply Lebesgue integrals over the domain, and over the corresponding Neumann boundary $\partial_N$, since no jump occurs there}.
\end{observation}   
With \cref{obs:forces} we can state the discrete weak forms as
\begin{align}
    \int_\body \con{\D \bm{\tau}}{ \D \bm{\sigma}} + \dfrac{1}{1+\nu} [\con{\Di \bm{\tau}}{\nabla \tr \bm{\sigma}} + \con{\nabla \tr \bm{\tau}}{\Di \bm{\sigma}}]  \, \dd \body &= 2\con{\bm{\tau}}{\sym\D \vb{f}}_{\mathcal{T}} + \dfrac{1+\nu^2}{1-\nu^2} \con{\tr \bm{\tau}}{ \di \vb{f}}_{\mathcal{T}} \label{eq:var3d} \\
    &\qquad + \int_{\surf_{N}} \con{\bm{\tau}}{\bm{\kappa}} \, \dd \surf   \qquad \forall \, \bm{\tau} \in \X^h(\body) \, , \notag 
\end{align}
for the three-dimensional problem, and
\begin{align}
    \int_\surf \psi \con{\Di \bm{\tau}}{ \Di \bm{\sigma} } + \chi \con{\nabla \tr \bm{\tau}}{ \nabla \tr \bm{\sigma} } \, \dd \surf = \psi \con{\bm{\tau}}{\sym \D \vb{f}}_{\mathcal{T}} + \con{\tr \bm{\tau}}{\di \vb{f}}_{\mathcal{T}} + \int_{\curv_N}  \con{\bm{\tau}}{\bm{\kappa}} \, \dd \curv \qquad \forall \, \bm{\tau} \in \Y^h(\surf) \, , \label{eq:var2d}
\end{align}
for the planar problems of plane stress and plane strain.

\section{Numerics of three-dimensional solids}
In the following, relative errors are measured using the Lebesgue norm $\norm{\widetilde{\bm{\sigma}} - \bm{\sigma}^h}_{\Le}/\norm{\widetilde{\bm{\sigma}}}_{\Le}$, where $\widetilde{\bm{\sigma}}$ represents the exact solution and $\bm{\sigma}^h$ is the finite element approximation.

\subsection{Stress convergence} \label{sec:con}
\begin{figure}
    	\centering
    	\begin{subfigure}{0.3\linewidth}
    		\centering
    		\includegraphics[width = 1\linewidth]{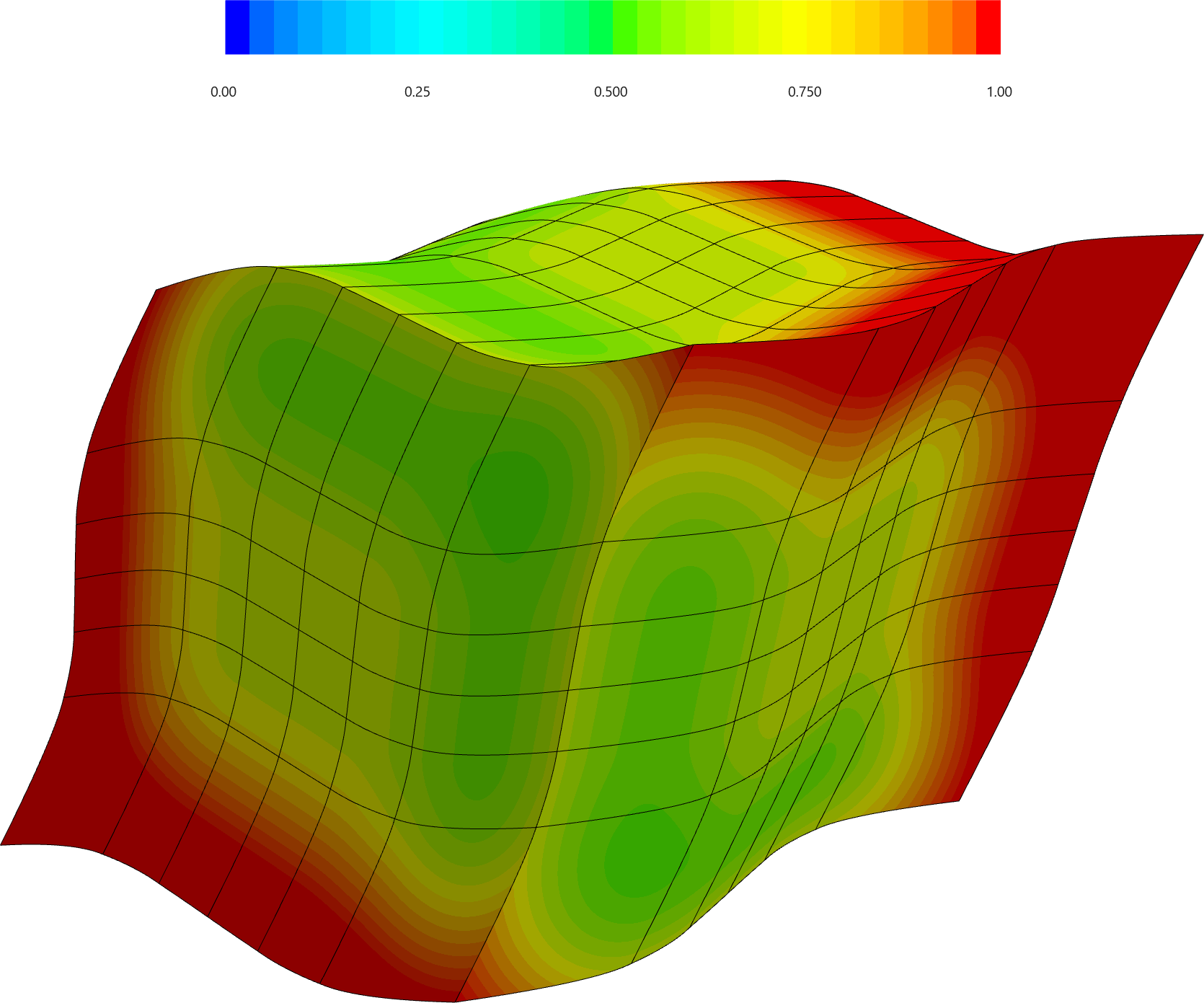}
    		\caption{}
    	\end{subfigure}
    	\begin{subfigure}{0.3\linewidth}
    		\centering
    		\includegraphics[width = 0.7\linewidth]{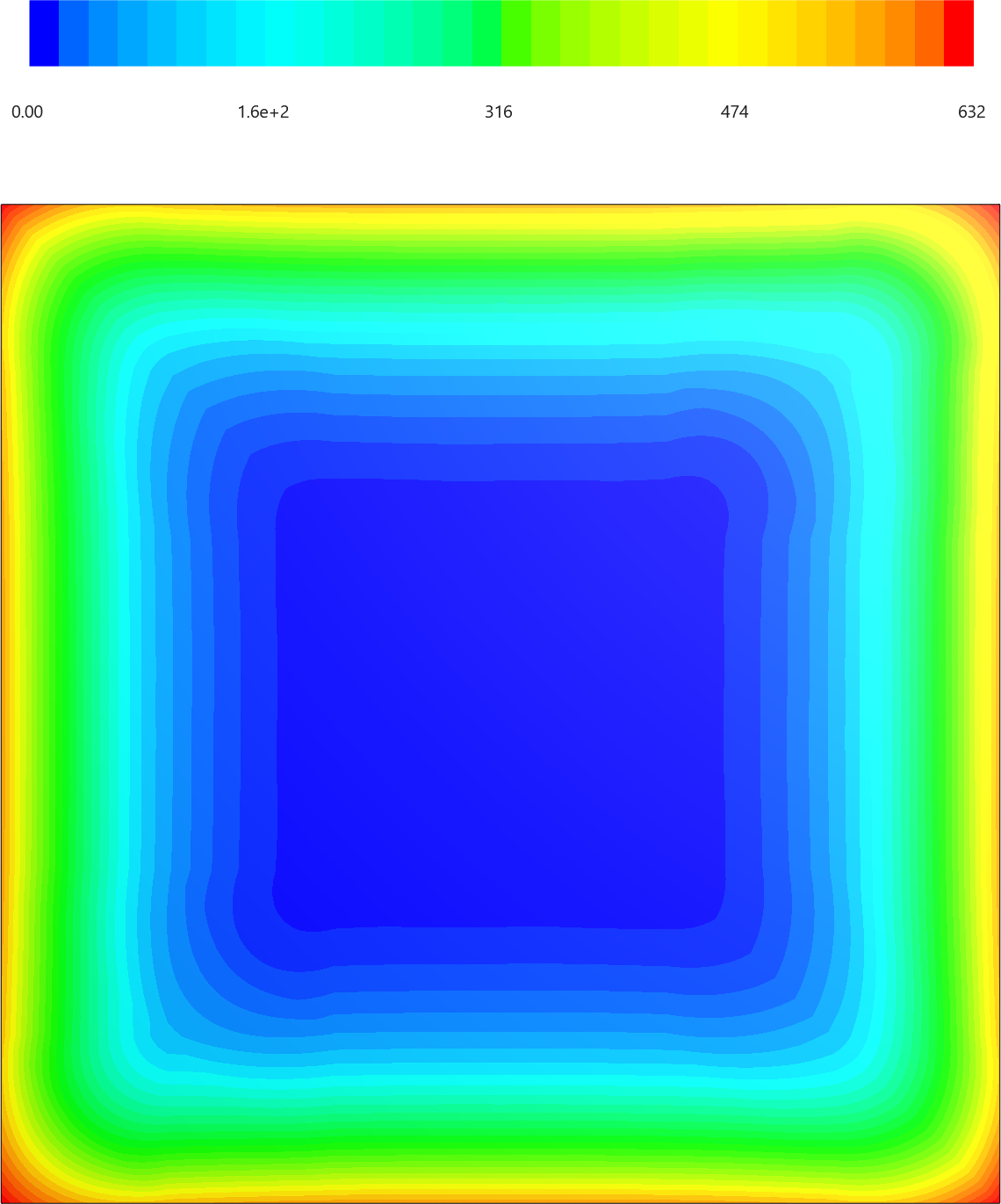}
    		\caption{}
    	\end{subfigure}
        \begin{subfigure}{0.3\linewidth}
    		\centering
    		\includegraphics[width = 0.7\linewidth]{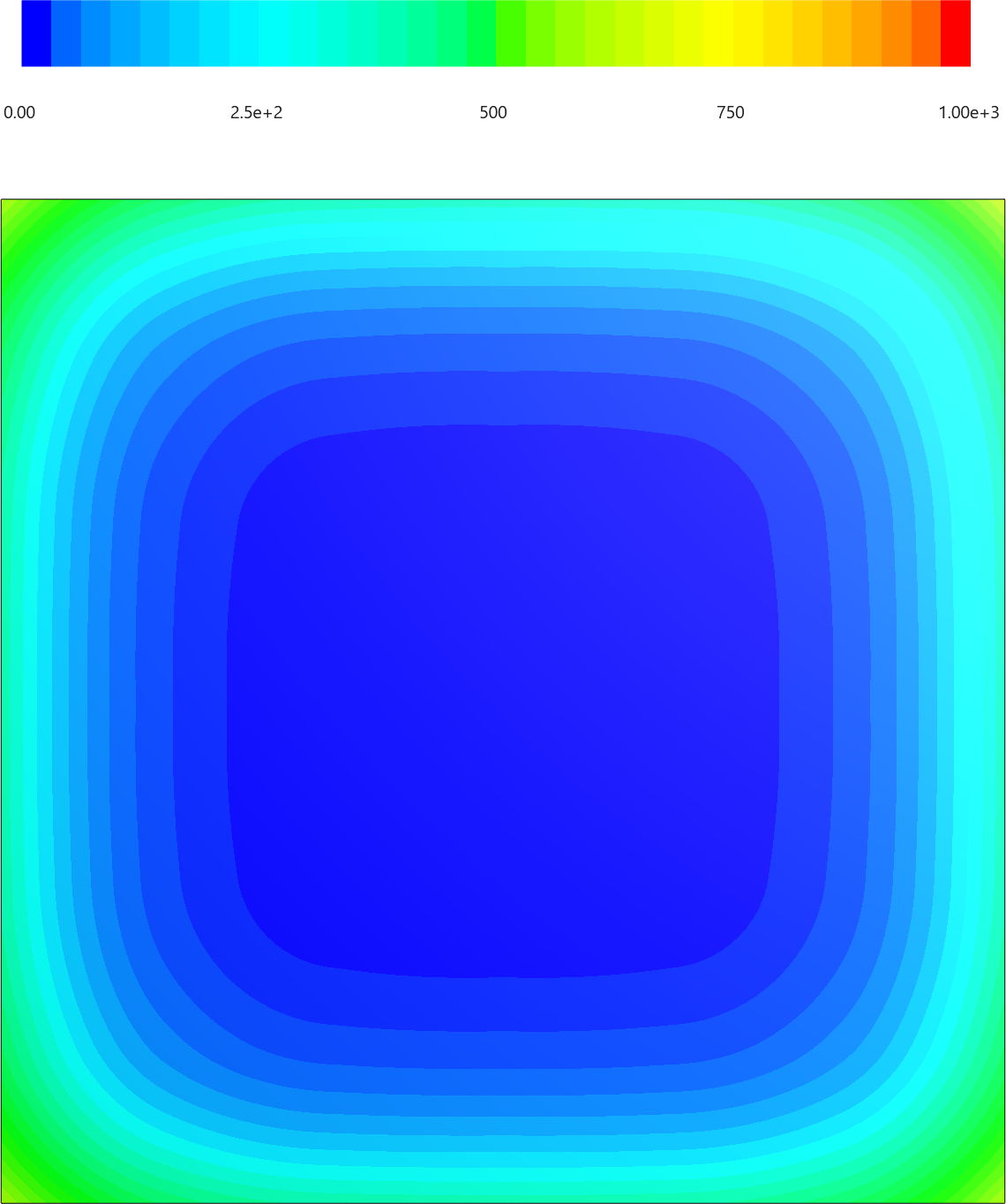}
    		\caption{}
    	\end{subfigure}
    	\caption{Depiction of the displacement field $\widetilde{\vb{u}}$ (a) with the corresponding von Mises $\sigma_v$ (b) and mean stresses $\sigma_m$ (c).
        The stresses are depicted on the middle-plane of the cube.}
    	\label{fig:disp3d}
\end{figure}

We start with the three-dimensional problem. Let the domain be the axis-symmetric cube $\overline{\body} = [-1,1]^3$, with a complete Dirichlet boundary $\surf_D = \partial \body$, we construct an artificial analytical solution by defining the material constants $E = 200$, $\nu = 0.25$ and the displacement field 
\begin{align}
    \widetilde{\vb{u}} &= \dfrac{1}{2} [(x^5 + y^5) \vb{e}_1 + (y^5 + z^5) \vb{e}_2 + (z^5 + x^5)\vb{e}_3] \, , 
\end{align}
from which we retrieve the stress tensor and the forces
\begin{align}
    &\widetilde{\bm{\sigma}} = \A \sym \D \widetilde{\vb{u}} = 200\left[\begin{matrix}3 x^{4} +  y^{4} +  z^{4} &  y^{4} &  x^{4}\\ y^{4} &  x^{4} + 3 y^{4} +  z^{4} & z^{4}\\ x^{4} &  z^{4} &  x^{4} +  y^{4} + 3 z^{4}\end{matrix}\right] \, , &&
    \vb{f} = - \Di \widetilde{\bm{\sigma}} =  800\left[\begin{matrix}- 3 x^{3} -  y^{3}\\- 3 y^{3} -  z^{3}\\-  x^{3} - 3 z^{3}\end{matrix}\right] \, ,
\end{align}
by using the compliance relation and static equilibrium. Further, from the stress tensor we extract the von Mises $\widetilde{\sigma}_v = \sqrt{3/2} \norm{\dev \bm{\sigma}}$ and mean stresses $\widetilde{\sigma}_m = (1/3) \tr \bm{\sigma}$, which are both important invariants of the stress tensor in design.
The displacement field along with the von Mises stress and the mean stress are depicted in \cref{fig:disp3d}.
For comparison we compute the results of our newly introduced stress based formulation $\bm{\sigma}^h$ and the stress field retrieved by post-processing in the displacement based formulation $\bm{\sigma}(\vb{u}^h)$.   
Convergence towards the exact solution is studied on structured meshes of uniform hexahedra of three polynomial orders $p \in \{1,2,3\}$. From the the results in \cref{fig:ex1res} we observe that the stress formulation converges in the optimal rate in all three measurements and yields comparable results to the displacement-based formulation, thus validating the method for a complete Dirichlet boundary. 
\begin{figure}
    	\centering
    	\begin{subfigure}{0.3\linewidth}
    		\centering
    		\begin{tikzpicture}[scale = 0.6]
    			\definecolor{asl}{rgb}{0.4980392156862745,0.,1.}
    			\definecolor{asb}{rgb}{0.,0.4,0.6}
    			\begin{loglogaxis}[
    				/pgf/number format/1000 sep={},
    				axis lines = left,
    				xlabel={Degrees of freedom},
    				ylabel={$\| \widetilde{\bm{\sigma}} - \bm{\sigma}^h \|_{\Le}/\| \widetilde{\bm{\sigma}}\|_{\Le}$},
    				xmin=50, xmax=2e5,
    				ymin=2e-4, ymax=2,
    				xtick={1e2,1e3,1e4,1e5},
    				ytick={1e-3,1e-2,1e-1,1e+0},
    				legend style={at={(0.05,0.05)},anchor= south west},
    				ymajorgrids=true,
    				grid style=dotted,
    				]
    				\addplot[color=asl, mark=triangle] coordinates {
    					(162, 0.39400190125935514)
(384, 0.3130547109627789)
(750, 0.20369691899663336)
(1296, 0.1403734879867633)
(2058, 0.10205078097674194)
(3072, 0.07737381141521538)
(4374, 0.060622745913486056)
(6000, 0.04875503584370207)
(7986, 0.040049683864567816)
(10368, 0.033478716319958884)
(13182, 0.02839896694703428)
(16464, 0.024391896987540664)
(20250, 0.021175862532618862)
(24576, 0.01855580100450658)
    				};
    				\addlegendentry{$\bm{\sigma}^{h.1}$}
    				
    				\addplot[color=violet, mark=pentagon] coordinates {
                        (750, 0.21215172480218902)
(2058, 0.07670543201004319)
(4374, 0.03510794813632964)
(7986, 0.018848688323385685)
(13182, 0.011259947387495263)
(20250, 0.007255823286460624)
(29478, 0.004946784096548262)
    				};
    				\addlegendentry{$\bm{\sigma}^{h.2}$}
    				
    				\addplot[color=purple, mark=square] coordinates {
    					(2058, 0.027477550521765867)
(6000, 0.005597508566867095)
(13182, 0.001795615153676065)
(24576, 0.0007411145864537414)
(41154, 0.00035911632742680977)
    				};
    				\addlegendentry{$\bm{\sigma}^{h.3}$}
    				
    				\addplot[color=blue, mark=10-pointed star] coordinates {
    					(81, 0.8028324304038807)
(192, 0.6401033802454071)
(375, 0.517284258593392)
(648, 0.42750852623380503)
(1029, 0.36199774605317653)
(1536, 0.3129761525039168)
(2187, 0.2752317033093599)
(3000, 0.24540521575874077)
(3993, 0.22130156734462636)
(5184, 0.20144761375454726)
(6591, 0.1848262782903259)
(8232, 0.17071630095916068)
(10125, 0.15859363120264475)
(12288, 0.14806915607560728)
(14739, 0.13884832352407417)
(17496, 0.13070429713618753)
(20577, 0.12345968800127359)
(24000, 0.1169738560191666)
    				};
    				\addlegendentry{$\bm{\sigma}(\vb{u}^{h.1})$}
    				
    				\addplot[color=cyan, mark=oplus] coordinates {
    					(375, 0.3790285947407449)
(1029, 0.20668489695069256)
(2187, 0.12256363243967744)
(3993, 0.08000916680606102)
(6591, 0.05606061912798095)
(10125, 0.0413738402438495)
(14739, 0.03175490763976962)
(20577, 0.0251257472636189)
(27783, 0.020368920183875417)
    				};
    				\addlegendentry{$\bm{\sigma}(\vb{u}^{h.2})$}

                    \addplot[color=asb, mark=o] coordinates {
    					(1029, 0.10905780393082387)
(3000, 0.033495832300536675)
(6591, 0.014180073139638727)
(12288, 0.007249960513109572)
(20577, 0.004186724662511231)
(31944, 0.0026313238013881586)
    				};
    				\addlegendentry{$\bm{\sigma}(\vb{u}^{h.3})$}
    				
    				\addplot[dashed,color=black, mark=none]
    				coordinates {
    					(2e+3, 1e-2)
    					(2e+4, 0.00046415888336127795)
    				};
    			
    			    \addplot[dashed,color=black, mark=none]
    			    coordinates {
    			    	(1e+3, 1e-0)
    			    	(2e+4, 0.36840314986403866)
    			    };
    			\end{loglogaxis}

                \draw (3.5,0.9) 
    			node[anchor=south]{\tiny $\mathcal{O}(h^{4})$};

                \draw (4,4.8) 
    			node[anchor=south]{\tiny $\mathcal{O}(h^{1})$};
    		\end{tikzpicture}
    		\caption{}
    	\end{subfigure}
        \begin{subfigure}{0.3\linewidth}
    		\centering
    		\begin{tikzpicture}[scale = 0.6]
    			\definecolor{asl}{rgb}{0.4980392156862745,0.,1.}
    			\definecolor{asb}{rgb}{0.,0.4,0.6}
    			\begin{loglogaxis}[
    				/pgf/number format/1000 sep={},
    				axis lines = left,
    				xlabel={Degrees of freedom},
    				ylabel={$\| \widetilde{\sigma}_v - \sigma_v^h \|_{\Le}/\| \widetilde{\sigma}_v\|_{\Le}$},
    				xmin=50, xmax=2e5,
    				ymin=2e-4, ymax=2,
    				xtick={1e2,1e3,1e4,1e5},
    				ytick={1e-3,1e-2,1e-1,1e+0},
    				legend style={at={(0.05,0.05)},anchor= south west},
    				ymajorgrids=true,
    				grid style=dotted,
    				]
    				\addplot[color=asl, mark=triangle] coordinates {
    					(162, 0.47387535673109865)
(384, 0.3344728354527508)
(750, 0.21029210011668495)
(1296, 0.14142092343411988)
(2058, 0.10086414010110073)
(3072, 0.07531189235968248)
(4374, 0.05827726848795807)
(6000, 0.046390088415222235)
(7986, 0.03778144696305634)
(10368, 0.03135412359323018)
(13182, 0.026431954630158787)
(16464, 0.02258077398564055)
(20250, 0.019511862571063795)
(24576, 0.017027353856465526)
    				};
    				\addlegendentry{$\bm{\sigma}^{h.1}$}
    				
    				\addplot[color=violet, mark=pentagon] coordinates {
                        (750, 0.1874974373092282)
(2058, 0.0713127573217927)
(4374, 0.033238865744987155)
(7986, 0.018037092272794958)
(13182, 0.010847870163693735)
(20250, 0.007021692079819059)
(29478, 0.004802255051833247)
    				};
    				\addlegendentry{$\bm{\sigma}^{h.2}$}
    				
    				\addplot[color=purple, mark=square] coordinates {
    					(2058, 0.0276148902077691)
(6000, 0.005845444628853821)
(13182, 0.0019022662704213467)
(24576, 0.0007903816727537579)
(41154, 0.0003843309970781168)
    				};
    				\addlegendentry{$\bm{\sigma}^{h.3}$}
    				
    				\addplot[color=blue, mark=10-pointed star] coordinates {
    					(81, 0.9337304311125407)
(192, 0.7835906346176945)
(375, 0.6495149062604442)
(648, 0.5445211595983752)
(1029, 0.46523702222507507)
(1536, 0.4046698329655985)
(2187, 0.35738821933219456)
(3000, 0.3196574973242956)
(3993, 0.28894395736905165)
(5184, 0.26350446168447583)
(6591, 0.24211390919256695)
(8232, 0.2238917784777356)
(10125, 0.20819155913600496)
(12288, 0.19452915678126362)
(14739, 0.18253558986351298)
(17496, 0.17192507291570938)
(20577, 0.16247303813862435)
(24000, 0.15400070757843137)
    				};
    				\addlegendentry{$\bm{\sigma}(\vb{u}^{h.1})$}
    				
    				\addplot[color=cyan, mark=oplus] coordinates {
    					(375, 0.4705886649998098)
(1029, 0.2658561613245903)
(2187, 0.1593592683538338)
(3993, 0.10457783880114971)
(6591, 0.07350282304472532)
(10125, 0.05435523713738639)
(14739, 0.04177541551196226)
(20577, 0.033086601655947526)
(27783, 0.026841839935292374)
    				};
    				\addlegendentry{$\bm{\sigma}(\vb{u}^{h.2})$}

                    \addplot[color=asb, mark=o] coordinates {
    					(1029, 0.14392543634130664)
(3000, 0.04432426517647767)
(6591, 0.0187939766492537)
(12288, 0.009615484703843737)
(20577, 0.005554066132560782)
(31944, 0.003490751915954955)
    				};
    				\addlegendentry{$\bm{\sigma}(\vb{u}^{h.3})$}
    				
    				\addplot[dashed,color=black, mark=none]
    				coordinates {
    					(2e+3, 1e-2)
    					(2e+4, 0.00046415888336127795)
    				};
    			
    			    \addplot[dashed,color=black, mark=none]
    			    coordinates {
    			    	(1e+3, 1e-0)
    			    	(2e+4, 0.36840314986403866)
    			    };
    			\end{loglogaxis}

                \draw (3.5,0.9) 
    			node[anchor=south]{\tiny $\mathcal{O}(h^{4})$};

                \draw (4,4.8) 
    			node[anchor=south]{\tiny $\mathcal{O}(h^{1})$};
    		\end{tikzpicture}
    		\caption{}
    	\end{subfigure}
    	\begin{subfigure}{0.3\linewidth}
    		\centering
    		\begin{tikzpicture}[scale = 0.6]
    			\definecolor{asl}{rgb}{0.4980392156862745,0.,1.}
    			\definecolor{asb}{rgb}{0.,0.4,0.6}
    			\begin{loglogaxis}[
    				/pgf/number format/1000 sep={},
    				axis lines = left,
    				xlabel={Degrees of freedom},
    				ylabel={$\| \widetilde{\sigma}_m - \sigma_m^h \|_{\Le}/\| \widetilde{\sigma}_m\|_{\Le}$},
    				xmin=50, xmax=2e5,
    				ymin=2e-4, ymax=2,
    				xtick={1e2,1e3,1e4,1e5},
    				ytick={1e-3,1e-2,1e-1,1e+0},
    				legend style={at={(0.05,0.05)},anchor= south west},
    				ymajorgrids=true,
    				grid style=dotted,
    				]
    				\addplot[color=asl, mark=triangle] coordinates {
    					(162, 0.3705548971987655)
(384, 0.30727302069936385)
(750, 0.20195632723133197)
(1296, 0.14010134190778598)
(2058, 0.10235541300884424)
(3072, 0.07789815175976066)
(4374, 0.061214673803675415)
(6000, 0.04934818116807011)
(7986, 0.040615650908645606)
(10368, 0.03400655121004629)
(13182, 0.028885855300936153)
(16464, 0.024838779877549323)
(20250, 0.021585314910011765)
(24576, 0.018930992439076352)
    				};
    				\addlegendentry{$\bm{\sigma}^{h.1}$}
    				
    				\addplot[color=violet, mark=pentagon] coordinates {
                        (750, 0.2180743058258449)
(2058, 0.0780393819691882)
(4374, 0.035575324133432236)
(7986, 0.019052945490339465)
(13182, 0.011364080986315899)
(20250, 0.007315151281401416)
(29478, 0.004983477886607028)
    				};
    				\addlegendentry{$\bm{\sigma}^{h.2}$}
    				
    				\addplot[color=purple, mark=square] coordinates {
    					(2058, 0.027441921980808973)
(6000, 0.005531582913918415)
(13182, 0.0017669873951233971)
(24576, 0.0007278313674462518)
(41154, 0.0003523021303731759)
    				};
    				\addlegendentry{$\bm{\sigma}^{h.3}$}
    				
    				\addplot[color=blue, mark=10-pointed star] coordinates {
    					(81, 0.7653468585899376)
(192, 0.5974144928013093)
(375, 0.47716220451770075)
(648, 0.39159917108010034)
(1029, 0.33008593491335375)
(1536, 0.28449420834071965)
(2187, 0.24962333918375923)
(3000, 0.22220124655487553)
(3993, 0.20012207306984978)
(5184, 0.18198764262207165)
(6591, 0.16684034058667663)
(8232, 0.1540052975814652)
(10125, 0.14299461364318272)
(12288, 0.13344749630470362)
(14739, 0.12509175903409525)
(17496, 0.1177183894167982)
(20577, 0.11116433610638583)
(24000, 0.10530059785749252)
    				};
    				\addlegendentry{$\bm{\sigma}(\vb{u}^{h.1})$}
    				
    				\addplot[color=cyan, mark=oplus] coordinates {
    					(375, 0.3514904793534479)
(1029, 0.18838216550177095)
(2187, 0.11108195292530376)
(3993, 0.07230961060301339)
(6591, 0.050580427166363504)
(10125, 0.03728843512581097)
(14739, 0.028597750969608787)
(20577, 0.022615494675739256)
(27783, 0.01832663460970209)
    				};
    				\addlegendentry{$\bm{\sigma}(\vb{u}^{h.2})$}

                    \addplot[color=asb, mark=o] coordinates {
    					(1029, 0.09804321467951609)
(3000, 0.030067465426817503)
(6591, 0.012717334606247256)
(12288, 0.006499594979551056)
(20577, 0.003752907054490273)
(31944, 0.0023586480694802437)
    				};
    				\addlegendentry{$\bm{\sigma}(\vb{u}^{h.3})$}
    				
    				\addplot[dashed,color=black, mark=none]
    				coordinates {
    					(2e+3, 1e-2)
    					(2e+4, 0.00046415888336127795)
    				};
    			
    			    \addplot[dashed,color=black, mark=none]
    			    coordinates {
    			    	(1e+3, 1e-0)
    			    	(2e+4, 0.36840314986403866)
    			    };
    			\end{loglogaxis}

                \draw (3.5,0.9) 
    			node[anchor=south]{\tiny $\mathcal{O}(h^{4})$};

                \draw (4,4.8) 
    			node[anchor=south]{\tiny $\mathcal{O}(h^{1})$};
    		\end{tikzpicture}
    		\caption{}
    	\end{subfigure}
    	\caption{Convergence slopes for $\norm{\bm{\sigma}}_\Le$ (a), $\norm{\sigma_v}_\Le$ (b), and $\norm{\sigma_m}_\Le$ (c), for polynomial orders $p \in \{1,2,3\}$.}
    	\label{fig:ex1res}
\end{figure}

\subsection{Operator spectrum} \label{sec:opnorm}
\begin{figure}
    	\centering
    	\begin{subfigure}{0.3\linewidth}
    		\centering
    		\begin{tikzpicture}[scale = 0.6]
    			\definecolor{asl}{rgb}{0.4980392156862745,0.,1.}
    			\definecolor{asb}{rgb}{0.,0.4,0.6}
    			\begin{axis}[
    				/pgf/number format/1000 sep={},
    				axis lines = left,
    				xlabel={$\nu$},
    				ylabel={Negative eigenvalues $[\%]$},
    				xmin=-0.1, xmax=0.6,
    				ymin=-0.1, ymax=1.1,
    				xtick={0,0.125, 0.25, 0.375, 0.5},
    				ytick={0,1},
    				legend style={at={(0.95,1)},anchor= north east},
    				ymajorgrids=true,
    				grid style=dotted,
    				]
    				\addplot[color=asl, mark=triangle] coordinates {
    					( 0 , 0.9833333333333333 )
( 0.125 , 0.31666666666666665 )
( 0.25 , 0.2 )
( 0.375 , 0.06666666666666667 )
( 0.5 , 0.0 )
    				};
    				\addlegendentry{$\partial \body = \surf_N$}
    				
    				\addplot[color=blue, mark=pentagon] coordinates {
                        ( 0 , 0.21666666666666667 )
( 0.125 , 0.03333333333333333 )
( 0.25 , 0.016666666666666666 )
( 0.375 , 0.0 )
( 0.5 , 0.0 )
    				};
    				\addlegendentry{$\partial \body = \surf_N \cup \surf_D$}
    				
    				\addplot[color=asb, mark=diamond] coordinates {
    					( 0 , 0.0 )
( 0.125 , 0.0 )
( 0.25 , 0.0 )
( 0.375 , 0.0 )
( 0.5 , 0.0 )

    				};
    				\addlegendentry{$\partial \body = \surf_D$}
    			    				    				
    			\end{axis}
    		\end{tikzpicture}
      \caption{}
    	\end{subfigure}
    	\begin{subfigure}{0.3\linewidth}
    		\centering
    		\definecolor{zzttqq}{rgb}{0.,0.4,0.6}
\definecolor{qqqqff}{rgb}{0.,0.4,0.6}
\begin{tikzpicture}[scale=0.7,line cap=round,line join=round,>=triangle 45,x=1cm,y=1cm]
\clip(3,2) rectangle (8,7);
\fill[line width=2pt,color=zzttqq,fill=zzttqq,fill opacity=0.2] (3,6) -- (5,7) -- (5,4) -- (3,3) -- cycle;
\fill[line width=2pt,color=zzttqq,fill=zzttqq,fill opacity=0.2] (5,7) -- (8,6) -- (8,3) -- (5,4) -- cycle;
\fill[line width=2pt,color=zzttqq,fill=zzttqq,fill opacity=0.2] (3,6) -- (5,7) -- (8,6) -- (6,5) -- cycle;
\draw [line width=0.7pt,color=black] (3,6)-- (5,7);
\draw [line width=0.7pt,color=black] (5,7)-- (5,4);
\draw [line width=0.7pt,color=black] (5,4)-- (3,3);
\draw [line width=0.7pt,color=black] (3,3)-- (3,6);
\draw [line width=0.7pt,color=black] (3,3)-- (6,2);
\draw [line width=0.7pt,color=black] (6,2)-- (6,5);
\draw [line width=0.7pt,color=black] (6,5)-- (3,6);
\draw [line width=0.7pt,color=black] (3,6)-- (3,3);
\draw [line width=0.7pt,color=black] (5,7)-- (8,6);
\draw [line width=0.7pt,color=black] (8,6)-- (8,3);
\draw [line width=0.7pt,color=black] (8,3)-- (5,4);
\draw [line width=0.7pt,color=black] (5,4)-- (5,7);
\draw [line width=0.7pt,color=black] (6,2)-- (8,3);
\draw [line width=0.7pt,color=black] (8,3)-- (8,6);
\draw [line width=0.7pt,color=black] (8,6)-- (6,5);
\draw [line width=0.7pt,color=black] (6,5)-- (6,2);
\draw [line width=0.7pt,color=black] (3,6)-- (5,7);
\draw [line width=0.7pt,color=black] (5,7)-- (8,6);
\draw [line width=0.7pt,color=black] (8,6)-- (6,5);
\draw [line width=0.7pt,color=black] (6,5)-- (3,6);
\draw [line width=0.7pt,color=black] (3,3)-- (6,2);
\draw [line width=0.7pt,color=black] (6,2)-- (8,3);
\draw [line width=0.7pt,color=black] (8,3)-- (5,4);
\draw [line width=0.7pt,color=black] (5,4)-- (3,3);
\draw [-to,line width=0.7pt] (5.5,4.5) -- (7,4);
\draw [-to,line width=0.7pt] (5.5,4.5) -- (6.5,5);
\draw [-to,line width=0.7pt] (5.5,4.5) -- (5.5,6);

\draw (7,4.2) node[color=black,anchor=north west] {$x$};
\draw (6.5,4.8) node[color=black,anchor=south west] {$y$};
\draw (5.5,6) node[color=black,anchor=south] {$z$};

\begin{scriptsize}
\draw [fill=black] (5.5,4.5) circle (2pt);
\end{scriptsize}
\end{tikzpicture}
      \caption{}
    	\end{subfigure}
     \begin{subfigure}{0.3\linewidth}
    		\centering
      \includegraphics[width = 0.7\linewidth]{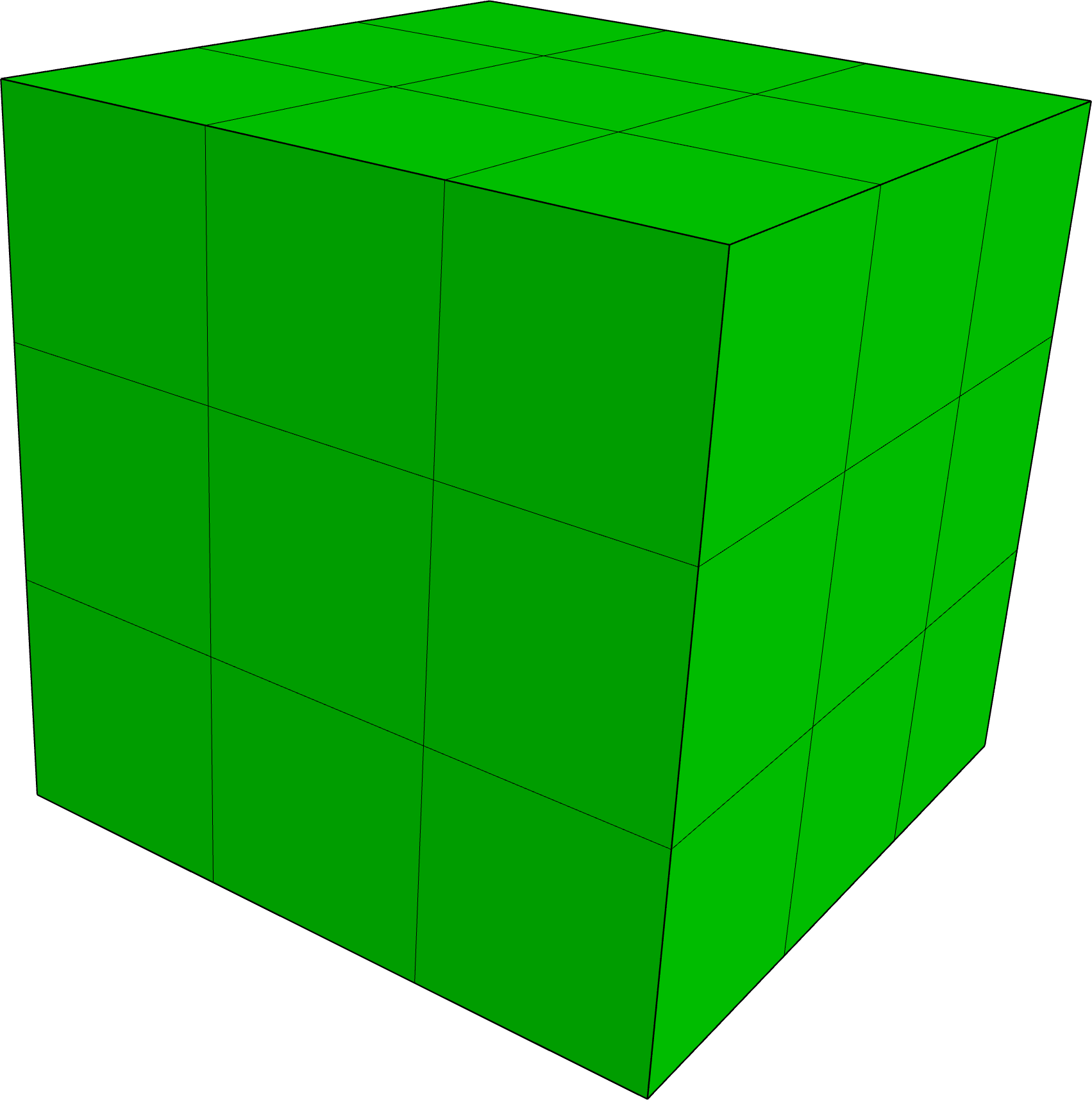}
    		\caption{}
    	\end{subfigure}
    	\caption{Percent of negative eigenvalues out of the total eigenvalues for varying Poisson ratios (a). Depiction of the Neumann boundary in the mixed domain (b) and the $3^3 = 27$-hexahedral mesh (c).}
    	\label{fig:opspec}
\end{figure}

The discrete three-dimensional variational problem in \cref{eq:var3d} can be written in operator form as
\begin{align}
    &\text{find} \quad \bm{\sigma}^h \in \X^h(\body) && \text{s.t.} &&  A^h \bm{\sigma}^h = \vb{f} \qquad \text{in} \qquad [\X^h(\body)]' \, , 
\end{align}
giving rise to the discrete operator $A^h$, which maps an element of $\X^h(\body)$ to the dual space $[\X^h(\body)]'$. By our proof, the operator is coercive under the very sharp restriction of a complete Dirichlet boundary $\surf_D = \partial \body$ if a conforming finite element space is employed $\X^h_0(\body) \subset \X_0(\body)$. However, often the boundary is not completely Dirichlet, such that Neumann boundary conditions are imposed. Thus, it is necessary to determine the behaviour of the operator for a mixed boundary. In the finite case this measure is given by the spectrum of the operator, which is characterised by its eigenvalues. Zero eigenvalues correspond with a positive semi-definite operator, whereas negative eigenvalues imply an indefinite operator, and correspondingly a non-coercive problem. In the following we count the number of negative and zero eigenvalues in relation to the full spectrum (all eigenvalues), which corresponds to the total number of degrees of freedom. We employ the cubical domain $\overline{\body} = [-1,1]^3$ and consider a uniform hexahedral mesh with $3^3 = 27$-elements of cubic polynomial order corresponding to $6000$-degrees of freedom. We examine the three boundary definitions 
\begin{align}
    &\partial \body = \surf_N \, , && \partial \body = \surf_N \cup \surf_D \, , && \partial \body =  \surf_D \, ,    
\end{align}
where the mixed boundary is defined such that $\surf_N$ is composed of the left $x = -1$, back, $y = -1$ and top $z = 1$ surfaces of the cube, and the Dirichlet boundary is $\surf_D =  \partial \body \setminus \surf_N$, see \cref{fig:opspec}. Further, we vary the Poisson ratio within its admissible bounds $\nu \in \{0, 0.125, 0.25, 0.375, 0.5\}$.   
From the results of the spectral analysis in \cref{fig:opspec} we observe that a domain with a Neumann boundary may introduce negative eigenvalues, depending on the Poisson ratio. Interestingly, for incompressible materials $\nu = 0.5$, no negative eigenvalues arise. For a total Neumann boundary we observe exactly $6$-zero eigenvalues irrespective of the Poisson ratio, which is consistent with the displacement formulation for which (in 3D) a total of $6$ rigid body motions have to be eliminated for the positive definiteness of the operator \cite{GEORGIYEVSKII2004941}. 
\begin{observation}[Discrete coercivity of form I]
    Given the numerical results of the eigenvalues of the discrete operator $A^h$ we surmise that for a complete Dirichlet boundary $\surf_D = \partial \body$ or an incompressible material $\nu = 0.5$ with a prescription of $\widetilde{\bm{\sigma}}$ on at least one node $|\surf_D| > 0$ the discrete variational problem is coercive.
\end{observation}

\subsection{Stabilised formulation}
In order to introduce a coercive form of the equations also for a mixed boundary conditions we add the zero-sum term derived from static equilibrium  
\begin{align}
    &-\omega \sym \D \Di \bm{\sigma} = \omega \sym \D \vb{f} \, , && \omega \geq 0 \, .
\end{align}
to the strong form in \cref{def:bvp3d}. The corresponding variational form and boundary term are given by testing with $\bm{\tau} \in \C^\infty(\overline{\body}) \otimes \Sym(3)$ and partial integration
\begin{align}
    - \omega \int_\body \con{\bm{\tau}}{\sym \D \Di \bm{\sigma}} \, \dd \body = \omega \int_\body \con{\Di \bm{\tau}}{\Di \bm{\sigma}} \, \dd \body - \omega \int_{\partial\body} \con{\bm{\tau}}{ \Di \bm{\sigma} \otimes \vb{n} } \, .
\end{align}
Thus, we recover the Neumann boundary term by splitting the boundary between Dirichlet and Neumann $\partial \body = \surf_D \cup \surf_N$, allowing us to state the modified boundary value problem.
\begin{definition}[The pure stress boundary value problem for solids II]
\label{def:bvp3dII}
    The field equations and boundary conditions of the stabilised symmetric three-dimensional stress boundary value problem read
    \begin{subequations}
        \begin{align}
            -\Delta \bm{\sigma} - \dfrac{1}{1+\nu} [\hess( \tr \bm{\sigma} ) + (\di \Di \bm{\sigma}) \one ] -\omega \sym \D \Di \bm{\sigma}  &= (2 + \omega) \sym \D \vb{f}   + \dfrac{1+\nu^2}{1-\nu^2} (\di \vb{f}) \one && \text{in} && \body \, ,  \\
            (\D \bm{\sigma})\vb{n} + \dfrac{1}{1+\nu} [\nabla \tr \bm{\sigma} \otimes \vb{n} - \con{\vb{f}}{\vb{n}}\one]  -\omega (\vb{f} \otimes \vb{n})  &= \bm{\kappa} && \text{on} && \surf_N \, , \\
         \bm{\sigma} &= \widetilde{\bm{\sigma}} && \text{on} && \surf_D  \, ,
        \end{align}
    \end{subequations}
\end{definition}
Further, we retrieve the corresponding variational form.
\begin{definition}[Variational form of the pure stress problem II]
    \label{def:var3DII}
    The weak formulation of the stabilised boundary value problem of linear elasticity written purely in stresses reads
    \begin{align}
    &\int_\body \con{\D \bm{\tau}}{ \D \bm{\sigma}} + \dfrac{1}{1+\nu} [\con{\Di \bm{\tau}}{\nabla \tr \bm{\sigma}} + \con{\nabla \tr \bm{\tau}}{\Di \bm{\sigma}}] + \omega \con{\Di \bm{\tau}}{\Di \bm{\sigma}}  \, \dd \body = \int_{\surf_{N}} \con{\bm{\tau}}{\bm{\kappa}} \, \dd \surf \notag \\
    &\qquad + \int_\body (2 + \omega)\con{\bm{\tau}}{\sym\D \vb{f}} + \dfrac{1+\nu^2}{1-\nu^2} \con{\tr \bm{\tau}}{ \di \vb{f}} \, \dd \body \qquad \forall \, \bm{\tau} \in \C_{\surf_D}^\infty(\body) \otimes \Sym(3) \, , \label{eq:weak3DII}
\end{align}
and yields a fully symmetric left-hand side. Here $\bm{\tau}$ is assumed to be compatible with the Dirichlet boundary.
\end{definition}
Lastly, the stabilised variational form can be derived from a variation functional.
\begin{definition}[Variation functional in pure stresses II]
\label{def:energy3DII}
    The weak formulation of the new boundary value problem can directly constructed as the variation of the stabilised functional
    \begin{align}
    I(\bm{\sigma}) &= \int_\body \dfrac{1}{2} \norm{\D \bm{\sigma}}^2 + \dfrac{1}{1+\nu} \con{\Di \bm{\sigma}}{\nabla \tr \bm{\sigma}} + \dfrac{\omega}{2} \norm{\Di \bm{\sigma}}^2 -  (2 + \omega) \con{\bm{\sigma}}{\sym \D \vb{f}} 
    \notag \\ & \qquad - \dfrac{1+\nu^2}{1-\nu^2} \con{\tr \bm{\sigma}}{ \di \vb{f}}  \, \dd \body - \int_{\surf_{N}} \con{\bm{\sigma}}{\bm{\kappa}} \, \dd \surf  \, , 
\end{align}
    with respect to the stress tensor $\bm{\sigma}$.
\end{definition}
We can now state the well-posedness of the stabilised problem.
\begin{theorem}[Stabilised bilinear form]
\label{th:stab}
    Let the bilinear and linear forms read
    \begin{subequations}
\label{eq:biforms}
    \begin{align}
        a(\bm{\tau},\bm{\sigma}) &= \int_\body \con{\D \bm{\tau}}{ \D \bm{\sigma}} + \chi [\con{\Di \bm{\tau}}{\nabla \tr \bm{\sigma}} + \con{\nabla \tr \bm{\tau}}{\Di \bm{\sigma}}] + \omega \con{\Di \bm{\tau}}{\Di \bm{\sigma}}  \, \dd \body \, , \\
        l(\bm{\tau}) &= \int_\body (2 + \omega)\con{\bm{\tau}}{\sym\D \vb{f}} + \dfrac{1+\nu^2}{1-\nu^2} \con{\tr \bm{\tau}}{ \di \vb{f}}  \, \dd \body \, ,
\end{align}
\end{subequations}
with $\omega \geq 3 \chi^2$, then the variational problem with a non-vanishing Dirichlet boundary $|\surf_D| > 0$, reading 
\begin{align}
    a(\bm{\tau}, \bm{\sigma}) &= l(\bm{\tau}) \qquad \forall \, \bm{\tau} \in \X(\body) \, ,  
\end{align}
has a unique solution $\bm{\sigma} \in \X(\body)$ for every right-hand side with the stability estimate
    \begin{align}
        \norm{\bm{\sigma}}_{\X} \leq \dfrac{1}{\beta} \norm{l}_{\X'} \, ,
    \end{align}
    where $\beta = \beta(\nu,\omega) > 0$. 
\end{theorem}
\begin{proof}
    Continuity follows analogously to the proof of \cref{th:wellposed3D} with the continuity constant $\alpha = 1 + 9 \chi$. For coercivity we observe
    \begin{align}
        a(\bm{\sigma},\bm{\sigma}) &= \norm{\D \bm{\sigma}}^2_{\Le} + 2 \chi \con{\Di \bm{\sigma}}{\nabla \tr \bm{\sigma}}_{\Le} + \omega \norm{\Di \bm{\sigma}}^2
        \notag \\
        & \overset{Y}{\geq} \norm{\D \bm{\sigma}}^2_{\Le} + (\omega - \epsilon  \chi) \norm{\Di \bm{\sigma}}_\Le^2 - \dfrac{\chi}{\epsilon} \norm{\nabla \tr \bm{\sigma}}_{\Le}^2 
        \notag \\
        & \overset{*}{\geq} \left(1 - \dfrac{3\chi}{\epsilon} \right )\norm{\D \bm{\sigma}}^2_{\Le} + (\omega - \epsilon  \chi) \norm{\Di \bm{\sigma}}_\Le^2 
        \notag \\
        & \overset{PF}{\geq} \dfrac{1}{1+c_F^2} \left(1 - \dfrac{3\chi}{\epsilon} \right )\norm{ \bm{\sigma}}^2_{\X} \qquad \forall \, \bm{\sigma} \in \X(\body) \, ,
    \end{align}
    where we used Young's inequality, the bound of $\norm{\nabla \tr \bm{\sigma}}_\Le^2 \leq 3 \norm{\D \bm{\sigma}}_\Le^2$, and the Poincar\'e-Friedrich inequality.
    The proof holds if $\epsilon > 3 \chi$, which is possible for the choice $\omega > 3\chi^2$.
\end{proof}
For the discrete variational problem the coercivity condition can be relaxed to $a(\bm{\sigma}^h, \bm{\sigma}^h) > 0$, implying a positive-definite matrix. Thus, we can inspect $\omega$ in the discrete case by investigating the eigenvalues of the newly induced operator $A^h$ of the stabilised bilinear form, characterising its spectrum. We apply the same test as in the previous benchmark with a complete Neumann boundary $\surf_N = \partial \body$ to find that for $\omega = \chi$ no negative eigenvalue arise. However, specifically and only for $\nu = 0$ we find 10 zero eigenvalues instead of $6$, which indicates a positive semi-definite operator. For $\omega = 1.01 \chi$ all negative eigenvalues vanish and we always retrieve exactly $6$ zero eigenvalues for all possible Poisson ratios $\nu \in [0, 0.5]$. We thus conclude that the discrete variational problem is coercive for $\omega > \chi$.  
\begin{definition}[Stabilised discrete variational problem for solids]
    The stabilised discrete variational problem reads
    \begin{align}
    &\int_\body \con{\D \bm{\tau}}{ \D \bm{\sigma}} + \dfrac{1}{1+\nu} [\con{\Di \bm{\tau}}{\nabla \tr \bm{\sigma}} + \con{\nabla \tr \bm{\tau}}{\Di \bm{\sigma}}] + \omega \con{\Di \bm{\tau}}{\Di \bm{\sigma}}  \, \dd \body \notag   \\
    &\qquad = \int_{\surf_{N}} \con{\bm{\tau}}{\bm{\kappa}} \, \dd \surf + (2+\omega) \con{\bm{\tau}}{\sym\D \vb{f}}_{\mathcal{T}} + \dfrac{1+\nu^2}{1-\nu^2} \con{\tr \bm{\tau}}{ \di \vb{f}}_{\mathcal{T}} \qquad \forall \, \bm{\tau} \in \X^h(\body) \, ,  \label{eq:var3dII}
\end{align}
and is coercive for any choice $\omega > \chi$.
\end{definition}

In order to compare the formulations we repeat the benchmark from \cref{sec:con} using the three-sided Neumann boundary definition from \cref{fig:opspec} (b) for the stabilised formulation with $\omega = 1.01 \chi$ and for the non-stabilised formulation over $\nu \in \{0, 0.25, 0.5\}$.  
\begin{figure}
    	\centering
    	\begin{subfigure}{0.3\linewidth}
    		\centering
    		\begin{tikzpicture}[scale = 0.6]
    			\definecolor{asl}{rgb}{0.4980392156862745,0.,1.}
    			\definecolor{asb}{rgb}{0.,0.4,0.6}
    			\begin{loglogaxis}[
    				/pgf/number format/1000 sep={},
    				axis lines = left,
    				xlabel={Degrees of freedom},
    				ylabel={$\| \widetilde{\bm{\sigma}} - \bm{\sigma}^h \|_{\Le}/\| \widetilde{\bm{\sigma}}\|_{\Le}$},
    				xmin=50, xmax=2e5,
    				ymin=2e-4, ymax=2,
    				xtick={1e2,1e3,1e4,1e5},
    				ytick={1e-3,1e-2,1e-1,1e+0},
    				legend style={at={(0.05,0.05)},anchor= south west},
    				ymajorgrids=true,
    				grid style=dotted,
    				]
    				\addplot[color=asl, mark=triangle] coordinates {
    					(162, 13.662906833095875)
(384, 0.6541701521054096)
(750, 0.38236077295392623)
(1296, 0.5265410186262903)
(2058, 0.28987021158695697)
(3072, 0.2961273457716303)
(4374, 0.961906334840834)
(6000, 1.1923125535852903)
(7986, 0.23560001066121647)
(10368, 1.09816855594621)
(13182, 0.11791119193395207)
(16464, 0.796343865222108)
(20250, 0.15915206387734668)
(24576, 0.1250408449561178)
    				};
    				\addlegendentry{$\bm{\sigma}^{h.1}$}
    				
    				\addplot[color=violet, mark=pentagon] coordinates {
                        (750, 0.6329201261045625)
(2058, 0.6994807180065562)
(4374, 0.07371569706741217)
(7986, 0.13377494567849255)
(13182, 0.07591477020527945)
(20250, 0.014663392078988983)
(29478, 0.021351038092826312)
    				};
    				\addlegendentry{$\bm{\sigma}^{h.2}$}
    				
    				\addplot[color=purple, mark=square] coordinates {
    					(2058, 0.042026730990191344)
(6000, 0.011511478449639073)
(13182, 0.003426863651922426)
(24576, 0.0010266253725154575)
(41154, 0.0016338071069316108)
    				};
    				\addlegendentry{$\bm{\sigma}^{h.3}$}
    				
    				\addplot[color=blue, mark=10-pointed star] coordinates {
    					(162, 1.1924422891749038)
(384, 0.7681899854927514)
(750, 0.507058628219259)
(1296, 0.3560354694203482)
(2058, 0.26281849705453586)
(3072, 0.20163432409752893)
(4374, 0.15943340961965421)
(6000, 0.12914316259769226)
(7986, 0.10668706044419808)
(10368, 0.08958843310286244)
(13182, 0.0762746241955097)
(16464, 0.06570904800806388)
(20250, 0.05718640038274006)
(24576, 0.050213488324374375)
    				};
    				\addlegendentry{$\bm{\sigma}_s^{h.1}$}
    				
    				\addplot[color=cyan, mark=oplus] coordinates {
    					(750, 0.34218934274653673)
(2058, 0.11449864816426047)
(4374, 0.050486824401255646)
(7986, 0.02642958347330865)
(13182, 0.015493550677092589)
(20250, 0.009837013965886569)
(29478, 0.006626362490625309)
    				};
    				\addlegendentry{$\bm{\sigma}_s^{h.2}$}

                    \addplot[color=asb, mark=o] coordinates {
    					(2058, 0.037083770545440715)
(6000, 0.007203979101023384)
(13182, 0.00225808343192521)
(24576, 0.0009191879016225698)
(41154, 0.00044130875359409366)
    				};
    				\addlegendentry{$\bm{\sigma}_s^{h.3}$}
    				
    				\addplot[dashed,color=black, mark=none]
    				coordinates {
    					(2e+3, 1e-2)
    					(2e+4, 0.00046415888336127795)
    				};
    			
    			\end{loglogaxis}

                \draw (3.5,0.9) 
    			node[anchor=south]{\tiny $\mathcal{O}(h^{4})$};

    		\end{tikzpicture}
    		\caption{}
    	\end{subfigure}
        \begin{subfigure}{0.3\linewidth}
    		\centering
    		\begin{tikzpicture}[scale = 0.6]
    			\definecolor{asl}{rgb}{0.4980392156862745,0.,1.}
    			\definecolor{asb}{rgb}{0.,0.4,0.6}
    			\begin{loglogaxis}[
    				/pgf/number format/1000 sep={},
    				axis lines = left,
    				xlabel={Degrees of freedom},
    				ylabel={$\| \widetilde{\bm{\sigma}} - \bm{\sigma}^h \|_{\Le}/\| \widetilde{\bm{\sigma}}\|_{\Le}$},
    				xmin=50, xmax=2e5,
    				ymin=2e-4, ymax=2,
    				xtick={1e2,1e3,1e4,1e5},
    				ytick={1e-3,1e-2,1e-1,1e+0},
    				legend style={at={(0.05,0.05)},anchor= south west},
    				ymajorgrids=true,
    				grid style=dotted,
    				]
    				\addplot[color=asl, mark=triangle] coordinates {
    					(162, 1.1823257668299443)
(384, 0.8875490095910311)
(750, 0.759101033759975)
(1296, 0.7694846507161439)
(2058, 0.9577755848125048)
(3072, 1.800187777127411)
(4374, 6.569077761841614)
(6000, 0.9381666606560242)
(7986, 0.45959758648512555)
(10368, 0.2877213348441578)
(13182, 0.2015542593326232)
(16464, 0.15079604401218932)
(20250, 0.11785426502448744)
(24576, 0.095034656980718)
    				};
    				\addlegendentry{$\bm{\sigma}^{h.1}$}
    				
    				\addplot[color=violet, mark=pentagon] coordinates {
                        (750, 9.059337113943162)
(2058, 0.12413069767437075)
(4374, 0.04890081270732632)
(7986, 0.025053949126465774)
(13182, 0.014728664774798275)
(20250, 0.009529130262781623)
(29478, 0.006677734276205932)
    				};
    				\addlegendentry{$\bm{\sigma}^{h.2}$}
    				
    				\addplot[color=purple, mark=square] coordinates {
    					(2058, 0.03246606721834037)
(6000, 0.006295342082263691)
(13182, 0.0019992565262421944)
(24576, 0.0008748663224401589)
(41154, 0.003334172460197108)
    				};
    				\addlegendentry{$\bm{\sigma}^{h.3}$}
    				
    				\addplot[color=blue, mark=10-pointed star] coordinates {
    					(162, 0.9178098354453458)
(384, 0.5715104179589847)
(750, 0.36786312237257185)
(1296, 0.2529919833593938)
(2058, 0.18368963168419947)
(3072, 0.13908304651306705)
(4374, 0.10881802672154747)
(6000, 0.08739200802833715)
(7986, 0.07169026685663764)
(10368, 0.0598509041406936)
(13182, 0.05070855294280969)
(16464, 0.04350483293165991)
(20250, 0.03772954174242912)
(24576, 0.03302947984251954)
    				};
    				\addlegendentry{$\bm{\sigma}_s^{h.1}$}
    				
    				\addplot[color=cyan, mark=oplus] coordinates {
    					(750, 0.2998878898104958)
(2058, 0.1002296014633521)
(4374, 0.04422752147203387)
(7986, 0.02317174524851714)
(13182, 0.013593978363514708)
(20250, 0.008636795522973518)
(29478, 0.005821373773178679)
    				};
    				\addlegendentry{$\bm{\sigma}_s^{h.2}$}

                    \addplot[color=asb, mark=o] coordinates {
    					(2058, 0.03282679030562198)
(6000, 0.006361319805813556)
(13182, 0.001990607022235663)
(24576, 0.0008094244908311059)
(41154, 0.0003883198449502816)
    				};
    				\addlegendentry{$\bm{\sigma}_s^{h.3}$}
    				
    				\addplot[dashed,color=black, mark=none]
    				coordinates {
    					(2e+3, 1e-2)
    					(2e+4, 0.00046415888336127795)
    				};
    			
    			\end{loglogaxis}

                \draw (3.5,0.9) 
    			node[anchor=south]{\tiny $\mathcal{O}(h^{4})$};

    		\end{tikzpicture}
    		\caption{}
    	\end{subfigure}
    	\begin{subfigure}{0.3\linewidth}
    		\centering
    		\begin{tikzpicture}[scale = 0.6]
    			\definecolor{asl}{rgb}{0.4980392156862745,0.,1.}
    			\definecolor{asb}{rgb}{0.,0.4,0.6}
    			\begin{loglogaxis}[
    				/pgf/number format/1000 sep={},
    				axis lines = left,
    				xlabel={Degrees of freedom},
    				ylabel={$\| \widetilde{\bm{\sigma}} - \bm{\sigma}^h \|_{\Le}/\| \widetilde{\bm{\sigma}}\|_{\Le}$},
    				xmin=50, xmax=2e5,
    				ymin=2e-4, ymax=2,
    				xtick={1e2,1e3,1e4,1e5},
    				ytick={1e-3,1e-2,1e-1,1e+0},
    				legend style={at={(0.05,0.05)},anchor= south west},
    				ymajorgrids=true,
    				grid style=dotted,
    				]
    				\addplot[color=asl, mark=triangle] coordinates {
    					(162, 0.760106407450243)
(384, 0.47794575178509197)
(750, 0.30900321743395526)
(1296, 0.21348226216786184)
(2058, 0.15559791818446145)
(3072, 0.11817592061130937)
(4374, 0.09268477981917983)
(6000, 0.07457774712161143)
(7986, 0.0612710893800588)
(10368, 0.05121432892010668)
(13182, 0.04343349906613143)
(16464, 0.037292710720807865)
(20250, 0.03236293269099819)
(24576, 0.028346403619301543)
    				};
    				\addlegendentry{$\bm{\sigma}^{h.1}$}
    				
    				\addplot[color=violet, mark=pentagon] coordinates {
                        (750, 0.240154020084217)
(2058, 0.08169182041524733)
(4374, 0.036328235729999075)
(7986, 0.019121254357813273)
(13182, 0.011254838593150612)
(20250, 0.0071692279646932205)
(29478, 0.00484255810031293)
    				};
    				\addlegendentry{$\bm{\sigma}^{h.2}$}
    				
    				\addplot[color=purple, mark=square] coordinates {
    					(2058, 0.027375966428977536)
(6000, 0.0053676451342808944)
(13182, 0.001689640265286692)
(24576, 0.0006895279711806555)
(41154, 0.00033160586429339446)
    				};
    				\addlegendentry{$\bm{\sigma}^{h.3}$}
    				
    				\addplot[color=blue, mark=10-pointed star] coordinates {
    					(162, 0.671993088347779)
(384, 0.42847187217109206)
(750, 0.2751726698036707)
(1296, 0.18868944247673827)
(2058, 0.13671140555936154)
(3072, 0.10336235530745735)
(4374, 0.08078882700071309)
(6000, 0.06483563864019332)
(7986, 0.05315957577776485)
(10368, 0.04436412503347633)
(13182, 0.03757726738552251)
(16464, 0.03223258469667089)
(20250, 0.02794958779889347)
(24576, 0.024465192486382365)

    				};
    				\addlegendentry{$\bm{\sigma}_s^{h.1}$}
    				
    				\addplot[color=cyan, mark=oplus] coordinates {
    					(750, 0.2412487327330682)
(2058, 0.08191775674104969)
(4374, 0.036377408232469675)
(7986, 0.019135435875593815)
(13182, 0.011260096643658997)
(20250, 0.007171622603118277)
(29478, 0.004843835964704118)
    				};
    				\addlegendentry{$\bm{\sigma}_s^{h.2}$}

                    \addplot[color=asb, mark=o] coordinates {
    					(2058, 0.027441574640851295)
(6000, 0.005375284618359145)
(13182, 0.0016914349525652495)
(24576, 0.0006901387228663628)
(41154, 0.00033186618565301574)
    				};
    				\addlegendentry{$\bm{\sigma}_s^{h.3}$}
    				
    				\addplot[dashed,color=black, mark=none]
    				coordinates {
    					(2e+3, 1e-2)
    					(2e+4, 0.00046415888336127795)
    				};
    			\end{loglogaxis}

                \draw (3.5,0.9) 
    			node[anchor=south]{\tiny $\mathcal{O}(h^{4})$};
    		\end{tikzpicture}
    		\caption{}
    	\end{subfigure}
    	\caption{Convergence slopes of the non-stabilised $\bm{\sigma}$ and stabilised $\bm{\sigma}_s$ forms with $p \in \{1,2,3\}$ for $\nu = 0$ (a), $\nu = 0.25$ (b) and $\nu = 0.499 \approx 0.5$ (c).}
    	\label{fig:exst}
\end{figure}
The convergence curves depicted in \cref{fig:exst} clearly demonstrate that the initial formulation without the stabilisation term leads to fluctuating results unless an incompressible material is used $\nu = 0.499 \approx 0.5$. In contrast, the stabilised formulation is robust regardless of the Poisson ratio. We note that the latter results are reproduced when using the Cholesky, PARDISO, or UMFPACK solvers in NGSolve.     

\section{Numerics of the planar problems}
\subsection{Characterisation of $\psi$} \label{sec:chara}
\begin{figure}
    	\centering
    	\begin{subfigure}{0.3\linewidth}
    		\centering
    		\includegraphics[width = 1\linewidth]{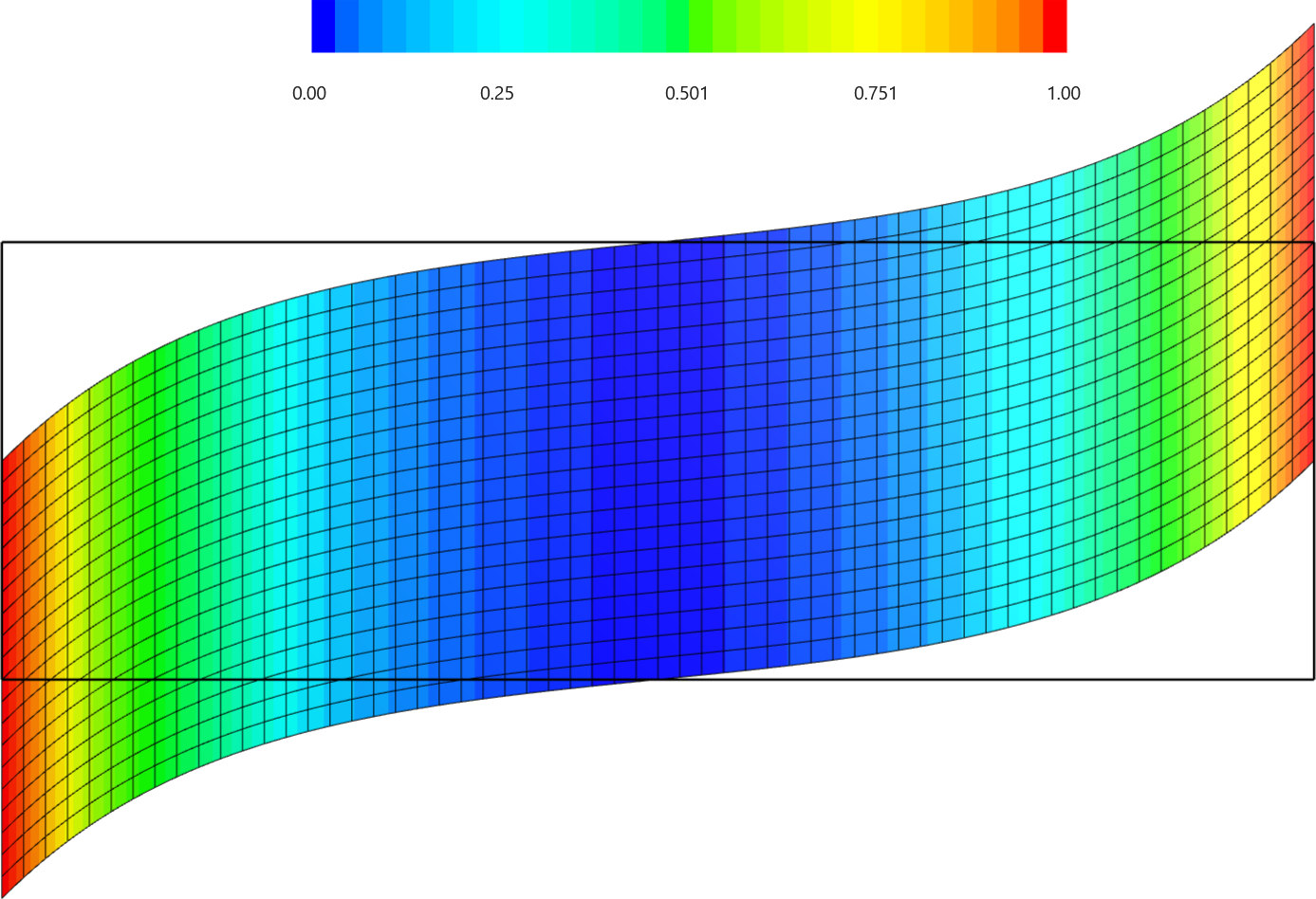}
      \includegraphics[width = 1\linewidth]{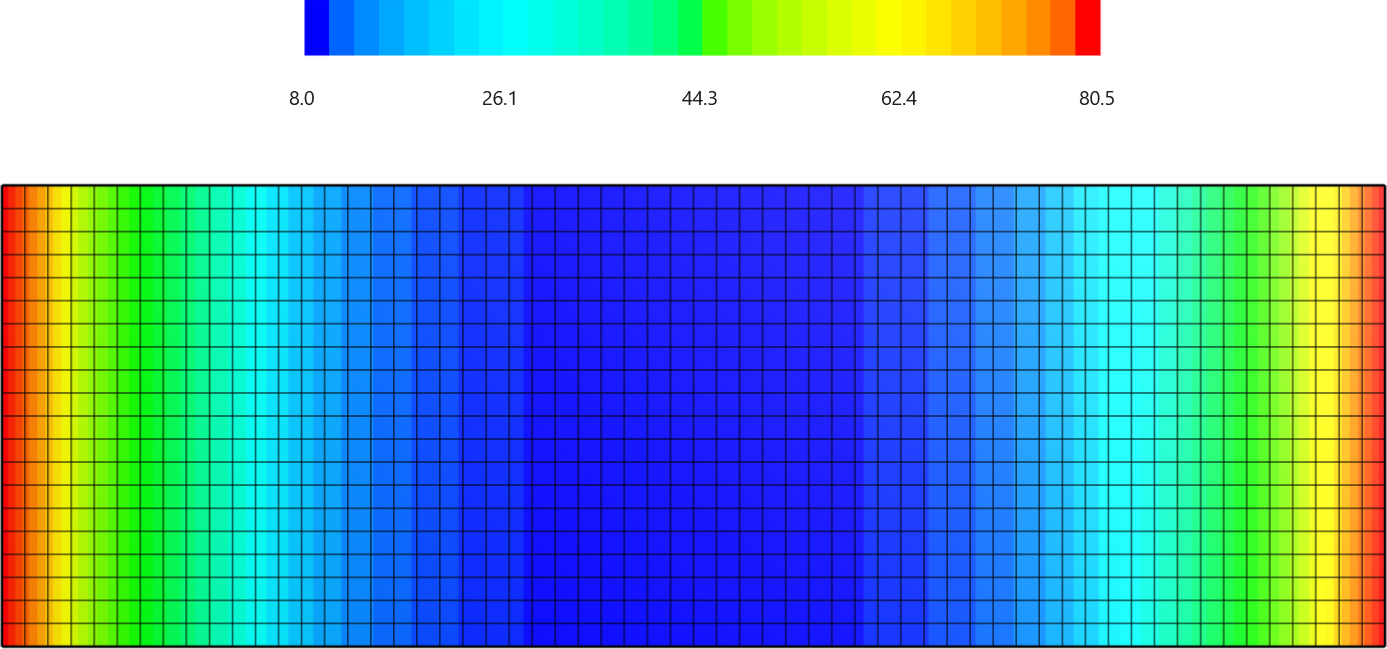}
    		\caption{}
    	\end{subfigure}
    	\begin{subfigure}{0.3\linewidth}
    		\centering
    		\includegraphics[width = 1\linewidth]{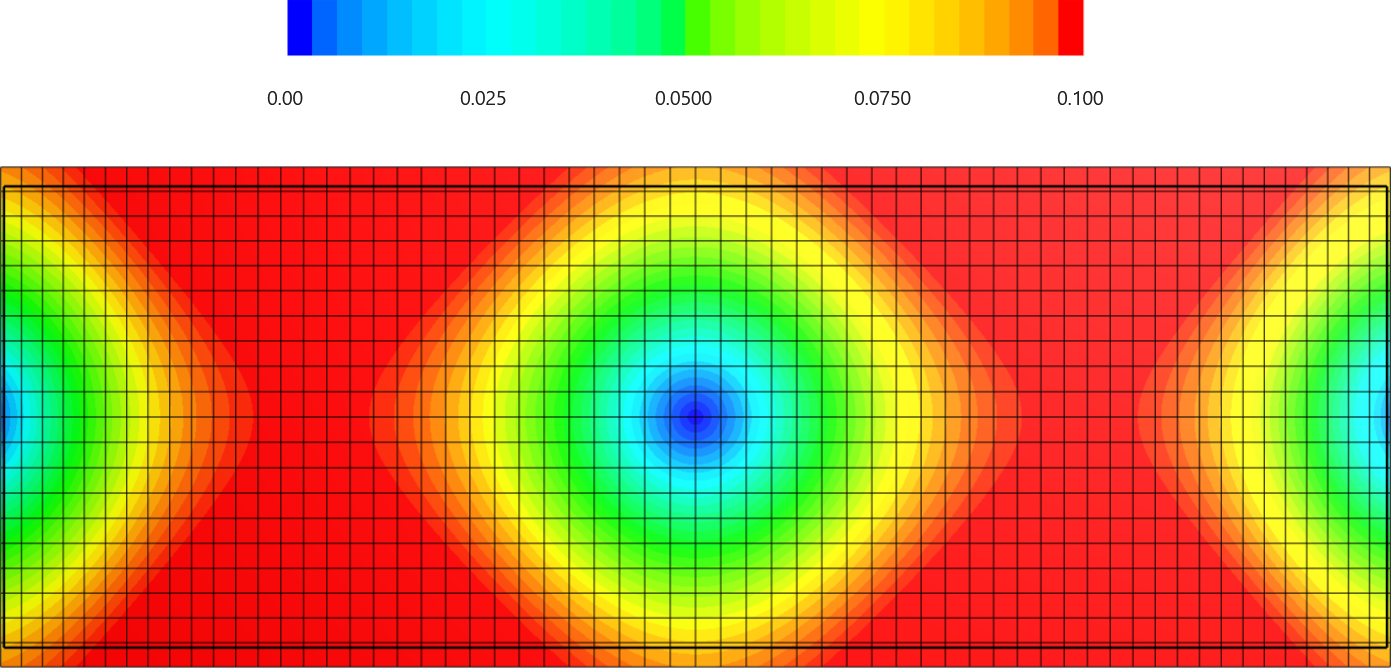}
      \includegraphics[width = 1\linewidth]{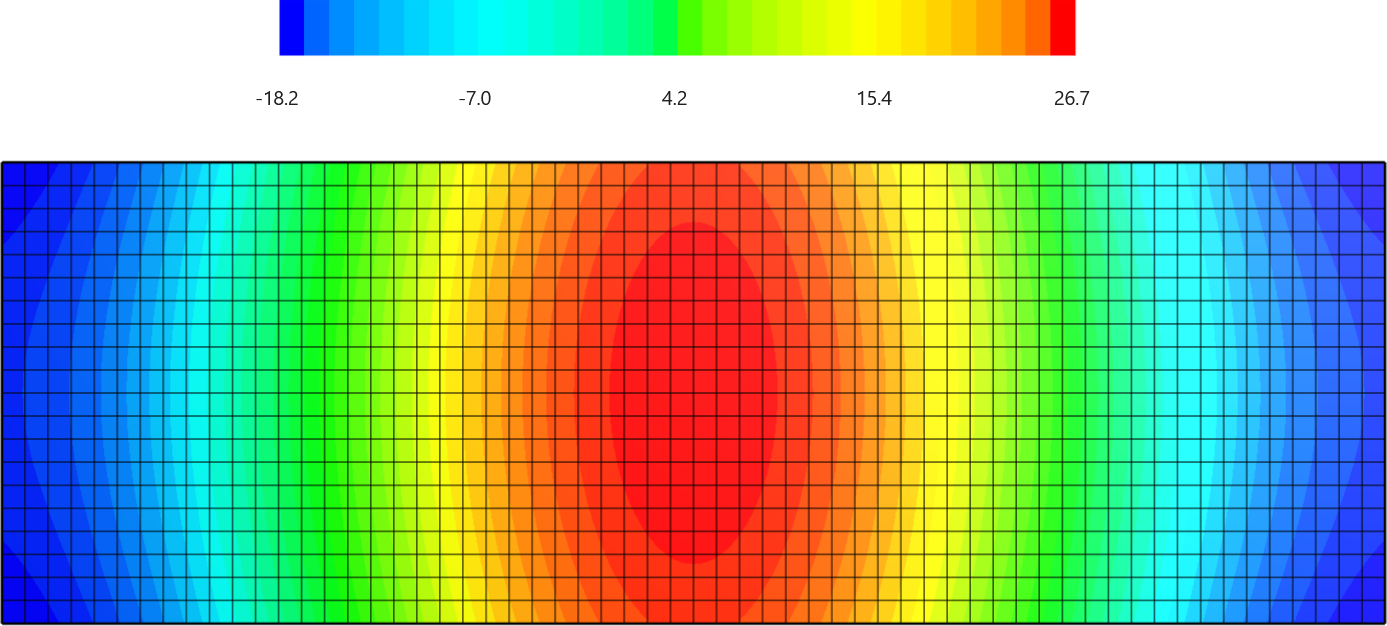}
      \includegraphics[width = 1\linewidth]{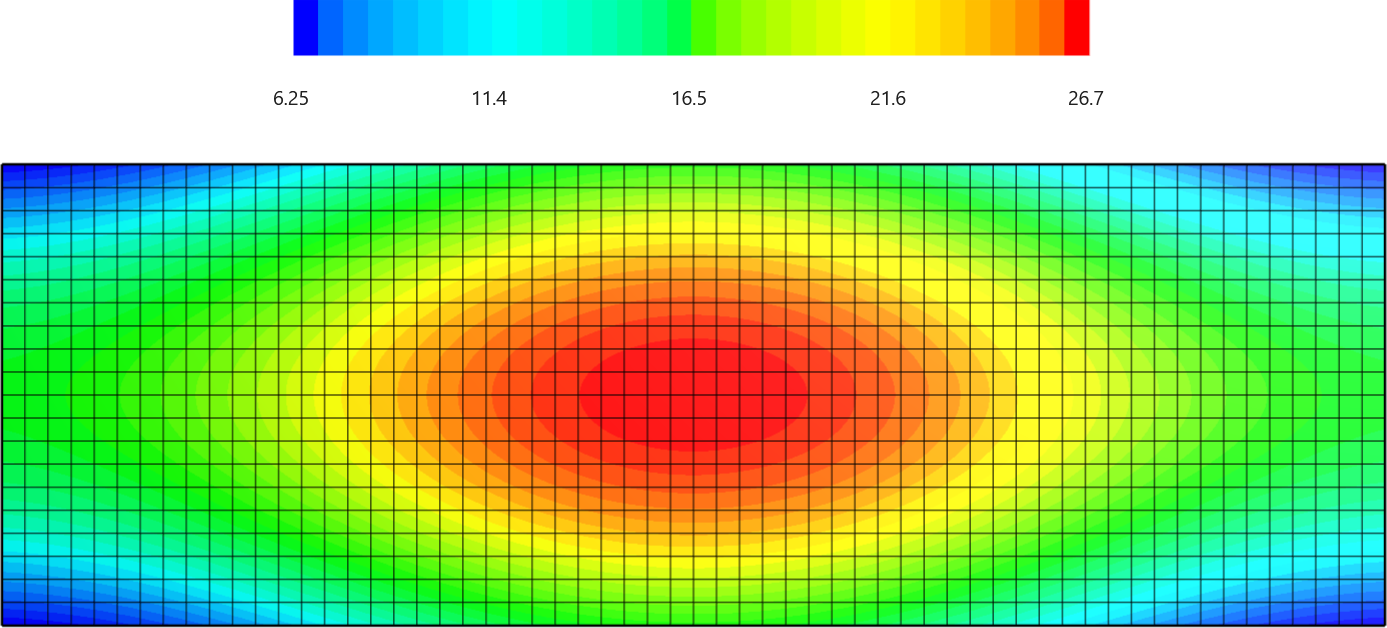}
    		\caption{}
    	\end{subfigure}
        \begin{subfigure}{0.3\linewidth}
    		\centering
    		\includegraphics[width = 1\linewidth]{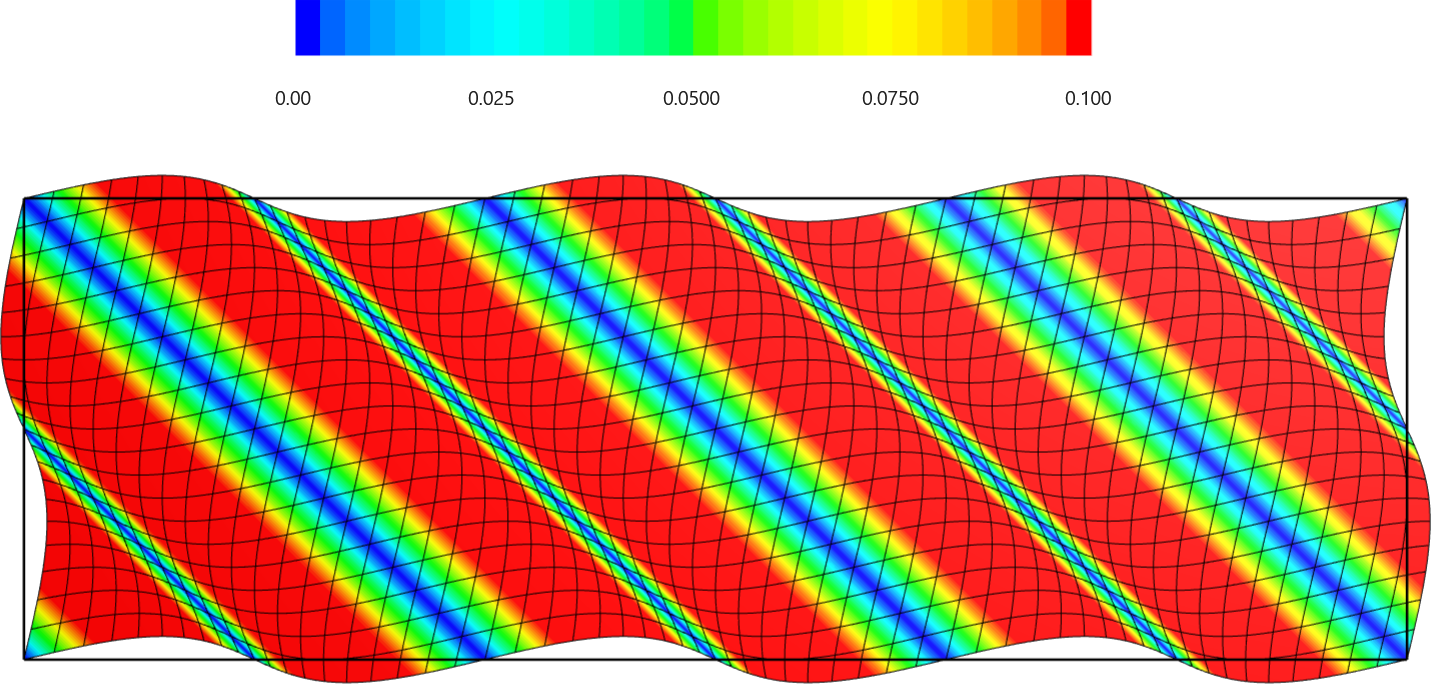}
      \includegraphics[width = 1\linewidth]{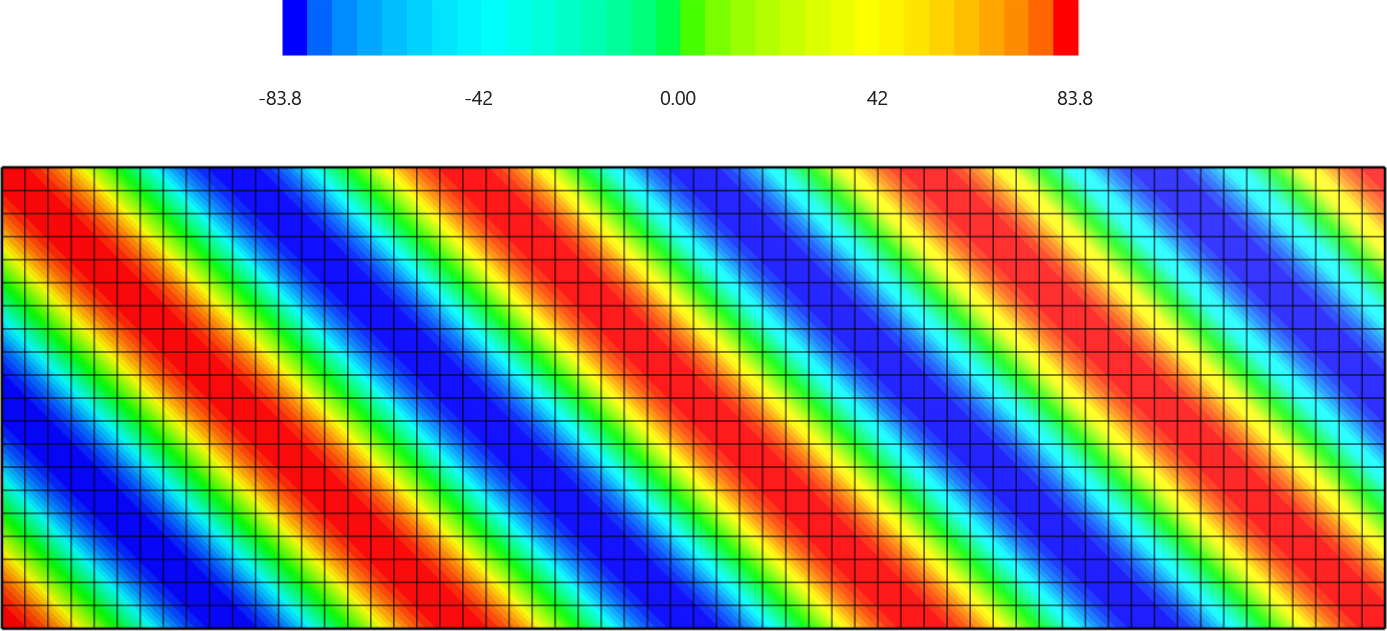}
      \includegraphics[width = 1\linewidth]{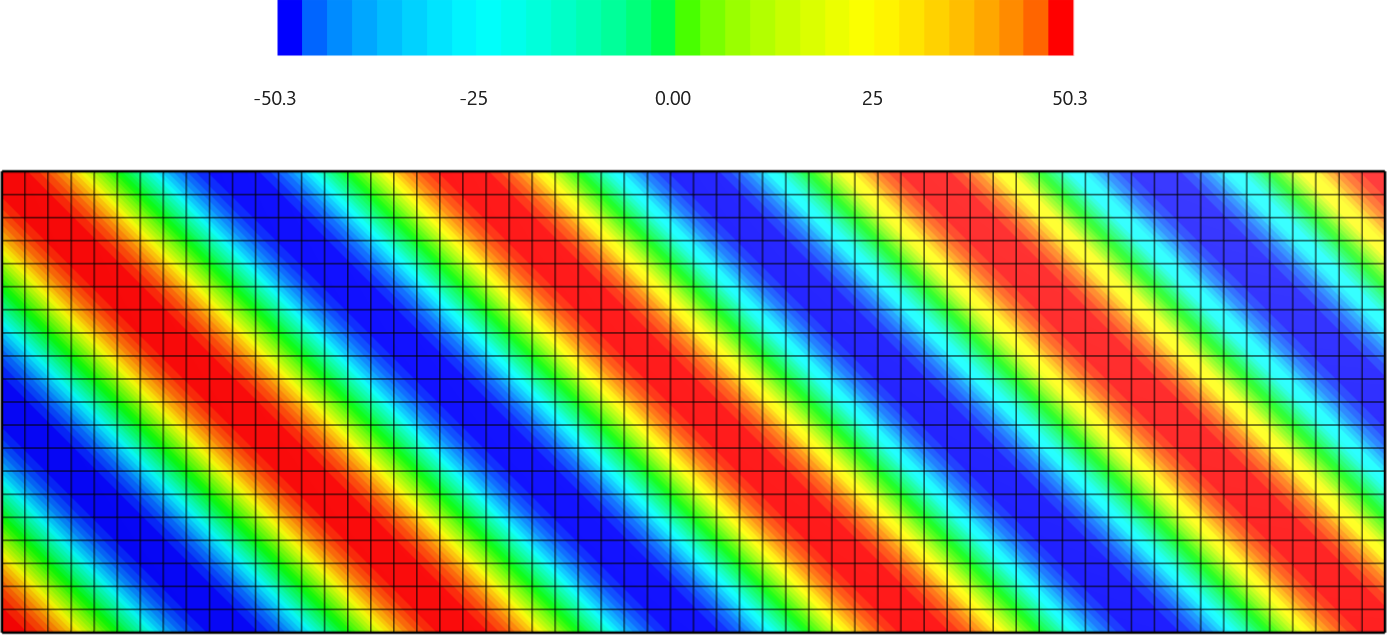}
    		\caption{}
    	\end{subfigure}
    	\caption{The displacement field $\widetilde{\vb{u}}_s$ and the corresponding stress field $\sigma_{12}$, where the remaining stress components vanish (a). The bi-axial tension field $\widetilde{\vb{u}}_b$ along with the stresses $\sigma_{11}$ and $\sigma_{22}$, such that the shear stress vanishes (b). The periodic field $\widetilde{\vb{u}}_p$ along with $\sigma_{11} = \sigma_{22}$ and $\sigma_{12}$ (c). For all the visualisations we set $\nu = 0.25$.}
    	\label{fig:psitest}
\end{figure}

In order to characterise appropriate values for the dimensionless parametre $\psi$ we examine plane stress and plane strain with varying Poisson ratios and Young's modulus set to $E = 200$. The domain is defined as $\overline{\surf} = [-3,3]\times[-1,1]$ with a complete Dirichlet boundary $\curv_D = \partial\surf$. We guarantee highly accurate approximations by relying on a structured mesh with $1200$ uniform quadrilateral elements of cubic polynomial order, yielding $33123$ degrees of freedom.   
We note beforehand that standard problems such as simple shear, pure extension and bending, lead to constant stress fields for the first two and a linear stress field in the case of bending, and as a result, the formulation captures these exactly, making them invalid for the characterisation of $\psi$. Consequently, we construct the following artificial analytical solutions 
\begin{align}
    &\widetilde{\vb{u}}_s = \dfrac{1}{10} \sinh (x) \vb{e}_2 \, , &&
    &\widetilde{\vb{u}}_b = \dfrac{1}{10} [\sin (x) \vb{e}_1 + \sin (y) \vb{e}_2]  \, , &&
    &\widetilde{\vb{u}}_p = \dfrac{1}{10} (\sin[\pi(x+y)] \vb{e}_1 + \sin[\pi(x+y)] \vb{e}_2) \, ,
\end{align}
which lead to a pure shear field, bi-axial tension, and periodicity. 
We retrieve the stresses and forces via the constitutive equation $\widetilde{\bm{\sigma}} = \Cm \sym \D \widetilde{\vb{u}}$ and static equilibrium $\vb{f} = - \Di \widetilde{\bm{\sigma}}$. We use $\widetilde{\bm{\sigma}}$ to enforce the Dirichlet boundary conditions and compute error estimates.
The fields and the resulting stresses for plane stress with $\nu = 0.25$ are depicted in \cref{fig:psitest}. 
The errors with respect to the exact solution are explored in \cref{fig:error} for the cases of plane stress and plane strain. The value of $\psi$ is defined with respect to $\chi$ using multiples of $k \in \{1e-3, 1e-2, \dots, 1e+2, 1e+3\}$ via $\psi = k \chi$, and the Poisson ration is varied as $\nu \in \{0,0.125,0.25, 0.375, 0.5\}$. In the case of plane strain we use $\nu = 0.499 \approx 0.5$ to derive the analytical solution for a quasi-incompressible material.
\begin{figure}
    	\centering
    	\begin{subfigure}{0.3\linewidth}
    		\centering
    		\begin{tikzpicture}[scale = 0.6]
    			\definecolor{asl}{rgb}{0.4980392156862745,0.,1.}
    			\definecolor{asb}{rgb}{0.,0.4,0.6}
    			\begin{loglogaxis}[
    				/pgf/number format/1000 sep={},
    				axis lines = left,
    				xlabel={$k$-value},
    				ylabel={$\| \widetilde{\bm{\sigma}} - \bm{\sigma}^h \|_{\Le}/\| \widetilde{\bm{\sigma}}\|_{\Le}$},
    				xmin=0.5*1e-3, xmax=2e3,
    				ymin=2e-8, ymax=5e-5,
    				xtick={1e-3,1e-2,1e-1,1,1e+1,1e+2,1e+3},
    				ytick={1e-8,1e-7,1e-6,1e-5,1e-4},
    				legend style={at={(0.95,1)},anchor= north east},
    				ymajorgrids=true,
    				grid style=dotted,
    				]
    				\addplot[color=asl, mark=triangle] coordinates {
    					( 0.001 , 4.922348774893941e-08 )
( 0.01 , 4.922501772899755e-08 )
( 0.1 , 4.924075114148115e-08 )
( 1 , 4.9435878461234236e-08 )
( 1.0 , 4.9435878461234236e-08 )
( 10.0 , 5.3168707225336725e-08 )
( 100.0 , 8.105872357788166e-08 )
( 1000.0 , 1.0131001243793102e-07 )
    				};
    				\addlegendentry{$\nu = 0$}
    				
    				\addplot[color=asb, mark=pentagon] coordinates {
                        ( 0.001 , 4.922348767571194e-08 )
( 0.01 , 4.922501772868709e-08 )
( 0.1 , 4.924074950629988e-08 )
( 1 , 4.9435881201875234e-08 )
( 1.125 , 4.94677957298107e-08 )
( 10.0 , 5.316870665066468e-08 )
( 100.0 , 8.1058723044589e-08 )
( 1000.0 , 1.0131000825300518e-07 )
    				};
    				\addlegendentry{$\nu = 0.125$}
    				
    				\addplot[color=violet, mark=square] coordinates {
    					( 0.001 , 4.922348681834029e-08 )
( 0.01 , 4.922501857800661e-08 )
( 0.1 , 4.924074901829029e-08 )
( 1 , 4.94358812414918e-08 )
( 1.25 , 4.950075393721495e-08 )
( 10.0 , 5.3168705943189246e-08 )
( 100.0 , 8.10587234194185e-08 )
( 1000.0 , 1.0131000579498556e-07 )
    				};
    				\addlegendentry{$\nu = 0.25$}
    				
    				\addplot[color=blue, mark=diamond] coordinates {
    					( 0.001 , 4.9223487365999876e-08 )
( 0.01 , 4.9225018081512354e-08 )
( 0.1 , 4.9240749955970404e-08 )
( 1 , 4.9435880696748644e-08 )
( 1.375 , 4.953471551787137e-08 )
( 10.0 , 5.316870643890856e-08 )
( 100.0 , 8.105872447556235e-08 )
( 1000.0 , 1.013100144848154e-07 )
    				};
    				\addlegendentry{$\nu = 0.375$}
    				
    				\addplot[color=cyan, mark=o] coordinates {
    					( 0.001 , 4.9223488952971577e-08 )
( 0.01 , 4.922501931905519e-08 )
( 0.1 , 4.9240750946017665e-08 )
( 1 , 4.943587998593981e-08 )
( 1.5 , 4.956964552051497e-08 )
( 10.0 , 5.3168709252722e-08 )
( 100.0 , 8.105872422873048e-08 )
( 1000.0 , 1.0131001056155572e-07 )
    				};
    				\addlegendentry{$\nu = 0.5$}
    				    				
    			\end{loglogaxis}
    		\end{tikzpicture}
      \caption{}
    	\end{subfigure}
    	\begin{subfigure}{0.3\linewidth}
    		\centering
    		\begin{tikzpicture}[scale = 0.6]
    			\definecolor{asl}{rgb}{0.4980392156862745,0.,1.}
    			\definecolor{asb}{rgb}{0.,0.4,0.6}
    			\begin{loglogaxis}[
    				/pgf/number format/1000 sep={},
    				axis lines = left,
    				xlabel={$k$-value},
    				ylabel={$\| \widetilde{\bm{\sigma}} - \bm{\sigma}^h \|_{\Le}/\| \widetilde{\bm{\sigma}}\|_{\Le}$},
    				xmin=0.5*1e-3, xmax=2e3,
    				ymin=2e-8, ymax=5e-5,
    				xtick={1e-3,1e-2,1e-1,1,1e+1,1e+2,1e+3},
    				ytick={1e-8,1e-7,1e-6,1e-5,1e-4},
    				legend style={at={(0.95,1)},anchor= north east},
    				ymajorgrids=true,
    				grid style=dotted,
    				]
    				\addplot[color=asl, mark=triangle] coordinates {
    					( 0.001 , 4.902400548876349e-08 )
( 0.01 , 4.9027530841703406e-08 )
( 0.1 , 4.906516010789614e-08 )
( 1 , 4.947321628675091e-08 )
( 1.0 , 4.947321628675091e-08 )
( 10.0 , 5.387758701371452e-08 )
( 100.0 , 6.55498421541749e-08 )
( 1000.0 , 7.009239756023263e-08 )
    				};
    				\addlegendentry{$\nu = 0$}
    				
    				\addplot[color=asb, mark=pentagon] coordinates {
                        ( 0.001 , 4.813870137946055e-08 )
( 0.01 , 4.8144234037215305e-08 )
( 0.1 , 4.817182711170255e-08 )
( 1 , 4.8490716364667506e-08 )
( 1.125 , 4.8540056958728927e-08 )
( 10.0 , 5.261034263798017e-08 )
( 100.0 , 6.471316654419156e-08 )
( 1000.0 , 6.957539824216099e-08 )
    				};
    				\addlegendentry{$\nu = 0.125$}
    				
    				\addplot[color=violet, mark=square] coordinates {
    					( 0.001 , 4.7003552608571386e-08 )
( 0.01 , 4.700650298566663e-08 )
( 0.1 , 4.702234485217927e-08 )
( 1 , 4.7234470337659774e-08 )
( 1.25 , 4.730807170091947e-08 )
( 10.0 , 5.092880463102422e-08 )
( 100.0 , 6.318882211758357e-08 )
( 1000.0 , 6.827319731627459e-08 )
    				};
    				\addlegendentry{$\nu = 0.25$}
    				
    				\addplot[color=blue, mark=diamond] coordinates {
    					( 0.001 , 4.572735840343113e-08 )
( 0.01 , 4.572811611292132e-08 )
( 0.1 , 4.573112085674382e-08 )
( 1 , 4.5827296924558226e-08 )
( 1.375 , 4.5896896851311784e-08 )
( 10.0 , 4.899435448792455e-08 )
( 100.0 , 6.116342368384329e-08 )
( 1000.0 , 6.637546158960342e-08 )
    				};
    				\addlegendentry{$\nu = 0.375$}
    				
    				\addplot[color=cyan, mark=o] coordinates {
    					( 0.001 , 4.441818809152503e-08 )
( 0.01 , 4.4416875609837705e-08 )
( 0.1 , 4.440714408103999e-08 )
( 1 , 4.4386097352567934e-08 )
( 1.5 , 4.442408181082424e-08 )
( 10.0 , 4.696668893994335e-08 )
( 100.0 , 5.8842296206510385e-08 )
( 1000.0 , 6.409947803563048e-08 )
    				};
    				\addlegendentry{$\nu = 0.5$}
    				    				
    			\end{loglogaxis}
    		\end{tikzpicture}
      \caption{}
    	\end{subfigure}
     \begin{subfigure}{0.3\linewidth}
    		\centering
    		\begin{tikzpicture}[scale = 0.6]
    			\definecolor{asl}{rgb}{0.4980392156862745,0.,1.}
    			\definecolor{asb}{rgb}{0.,0.4,0.6}
    			\begin{loglogaxis}[
    				/pgf/number format/1000 sep={},
    				axis lines = left,
    				xlabel={$k$-value},
    				ylabel={$\| \widetilde{\bm{\sigma}} - \bm{\sigma}^h \|_{\Le}/\| \widetilde{\bm{\sigma}}\|_{\Le}$},
    				xmin=0.5*1e-3, xmax=2e3,
    				ymin=2e-8, ymax=5e-5,
    				xtick={1e-3,1e-2,1e-1,1,1e+1,1e+2,1e+3},
    				ytick={1e-8,1e-7,1e-6,1e-5,1e-4},
    				legend style={at={(0.05,0.05)},anchor= south west},
    				ymajorgrids=true,
    				grid style=dotted,
    				]
    				\addplot[color=asl, mark=triangle] coordinates {
    					( 0.001 , 4.831118393992531e-06 )
( 0.01 , 4.83082646729399e-06 )
( 0.1 , 4.828090960120863e-06 )
( 1 , 4.815647153054629e-06 )
( 1.0 , 4.815647153054665e-06 )
( 10.0 , 5.077195954992945e-06 )
( 100.0 , 9.834781707678435e-06 )
( 1000.0 , 4.964618170994774e-05 )
    				};
    				\addlegendentry{$\nu = 0$}
    				
    				\addplot[color=asb, mark=pentagon] coordinates {
                        ( 0.001 , 4.905362924946483e-06 )
( 0.01 , 4.905087725676867e-06 )
( 0.1 , 4.902486931158953e-06 )
( 1 , 4.888767104353023e-06 )
( 1.125 , 4.888253812345791e-06 )
( 10.0 , 5.075513721650014e-06 )
( 100.0 , 8.934897675375146e-06 )
( 1000.0 , 4.3422722137319145e-05 )
    				};
    				\addlegendentry{$\nu = 0.125$}
    				
    				\addplot[color=violet, mark=square] coordinates {
    					( 0.001 , 4.971721584795066e-06 )
( 0.01 , 4.971462674559404e-06 )
( 0.1 , 4.968997193495019e-06 )
( 1 , 4.954356797513652e-06 )
( 1.25 , 4.952774552013529e-06 )
( 10.0 , 5.075967345831185e-06 )
( 100.0 , 8.029973296887549e-06 )
( 1000.0 , 3.6820255612819834e-05 )
    				};
    				\addlegendentry{$\nu = 0.25$}
    				
    				\addplot[color=blue, mark=diamond] coordinates {
    					( 0.001 , 5.025976427711357e-06 )
( 0.01 , 5.025732278746974e-06 )
( 0.1 , 5.0233921415444146e-06 )
( 1 , 5.008203680198447e-06 )
( 1.375 , 5.005179431274809e-06 )
( 10.0 , 5.0782557791652315e-06 )
( 100.0 , 7.19300116314077e-06 )
( 1000.0 , 3.021245523148524e-05 )
    				};
    				\addlegendentry{$\nu = 0.375$}
    				
    				\addplot[color=cyan, mark=o] coordinates {
    					( 0.001 , 5.066592919257201e-06 )
( 0.01 , 5.066361206462402e-06 )
( 0.1 , 5.064129800570042e-06 )
( 1 , 5.0487260579776755e-06 )
( 1.5 , 5.044132485007212e-06 )
( 10.0 , 5.081861368588671e-06 )
( 100.0 , 6.48651375301588e-06 )
( 1000.0 , 2.3950993987714765e-05 )
    				};
    				\addlegendentry{$\nu = 0.5$}
    				    				
    			\end{loglogaxis}
    		\end{tikzpicture}
      \caption{}
    	\end{subfigure}
        \begin{subfigure}{0.3\linewidth}
    		\centering
    		\begin{tikzpicture}[scale = 0.6]
    			\definecolor{asl}{rgb}{0.4980392156862745,0.,1.}
    			\definecolor{asb}{rgb}{0.,0.4,0.6}
    			\begin{loglogaxis}[
    				/pgf/number format/1000 sep={},
    				axis lines = left,
    				xlabel={$k$-value},
    				ylabel={$\| \widetilde{\bm{\sigma}} - \bm{\sigma}^h \|_{\Le}/\| \widetilde{\bm{\sigma}}\|_{\Le}$},
    				xmin=0.5*1e-3, xmax=2e3,
    				ymin=2e-8, ymax=5e-5,
    				xtick={1e-3,1e-2,1e-1,1,1e+1,1e+2,1e+3},
    				ytick={1e-8,1e-7,1e-6,1e-5,1e-4},
    				legend style={at={(0.95,1)},anchor= north east},
    				ymajorgrids=true,
    				grid style=dotted,
    				]
    				\addplot[color=asl, mark=triangle] coordinates {
    					( 0.001 , 4.9223487687391645e-08 )
( 0.01 , 4.922501753999173e-08 )
( 0.1 , 4.924075093646369e-08 )
( 1 , 4.943587834246358e-08 )
( 1.0 , 4.943587834246358e-08 )
( 10.0 , 5.316870677336018e-08 )
( 100.0 , 8.10587237457314e-08 )
( 1000.0 , 1.0131001283826343e-07 )
    				};
    				\addlegendentry{$\nu = 0$}
    				
    				\addplot[color=asb, mark=pentagon] coordinates {
                        ( 0.001 , 4.922348819557342e-08 )
( 0.01 , 4.922501893897112e-08 )
( 0.1 , 4.924074888548089e-08 )
( 1 , 4.943587934178958e-08 )
( 1.1428571428571428 , 4.947244154713961e-08 )
( 10.0 , 5.31687072364598e-08 )
( 100.0 , 8.105872449081354e-08 )
( 1000.0 , 1.0131001689424773e-07 )
    				};
    				\addlegendentry{$\nu = 0.125$}
    				
    				\addplot[color=violet, mark=square] coordinates {
    					( 0.001 , 4.9223488902849825e-08 )
( 0.01 , 4.92250192538705e-08 )
( 0.1 , 4.924075124132937e-08 )
( 1 , 4.9435880181891764e-08 )
( 1.3333333333333333 , 4.952328578846165e-08 )
( 10.0 , 5.316870779374647e-08 )
( 100.0 , 8.105872426890697e-08 )
( 1000.0 , 1.0131001490482918e-07 )
    				};
    				\addlegendentry{$\nu = 0.25$}
    				
    				\addplot[color=blue, mark=diamond] coordinates {
    					( 0.001 , 4.9223486735707373e-08 )
( 0.01 , 4.9225019859061075e-08 )
( 0.1 , 4.924074833705193e-08 )
( 1 , 4.943587978762838e-08 )
( 1.6 , 4.959826356362338e-08 )
( 10.0 , 5.316870702014554e-08 )
( 100.0 , 8.105872381760701e-08 )
( 1000.0 , 1.0131001243013988e-07 )
    				};
    				\addlegendentry{$\nu = 0.375$}
    				
    				\addplot[color=cyan, mark=o] coordinates {
    					( 0.001 , 4.922348760328797e-08 )
( 0.01 , 4.922501824694724e-08 )
( 0.1 , 4.9240749183846165e-08 )
( 1 , 4.943587951290344e-08 )
( 1.996007984031936 , 4.971715742203098e-08 )
( 10.0 , 5.316870676873016e-08 )
( 100.0 , 8.105872514199951e-08 )
( 1000.0 , 1.013100159819449e-07 )
    				};
    				\addlegendentry{$\nu \approx 0.5$}  
    			\end{loglogaxis}
    		\end{tikzpicture}
      \caption{}
    	\end{subfigure}
    	\begin{subfigure}{0.3\linewidth}
    		\centering
    		\begin{tikzpicture}[scale = 0.6]
    			\definecolor{asl}{rgb}{0.4980392156862745,0.,1.}
    			\definecolor{asb}{rgb}{0.,0.4,0.6}
    			\begin{loglogaxis}[
    				/pgf/number format/1000 sep={},
    				axis lines = left,
    				xlabel={$k$-value},
    				ylabel={$\| \widetilde{\bm{\sigma}} - \bm{\sigma}^h \|_{\Le}/\| \widetilde{\bm{\sigma}}\|_{\Le}$},
    				xmin=0.5*1e-3, xmax=2e3,
    				ymin=2e-8, ymax=5e-5,
    				xtick={1e-3,1e-2,1e-1,1,1e+1,1e+2,1e+3},
    				ytick={1e-8,1e-7,1e-6,1e-5,1e-4},
    				legend style={at={(0.95,1)},anchor= north east},
    				ymajorgrids=true,
    				grid style=dotted,
    				]
    				\addplot[color=asl, mark=triangle] coordinates {
    					( 0.001 , 4.9023690872677925e-08 )
( 0.01 , 4.9027488573559145e-08 )
( 0.1 , 4.9065158835939403e-08 )
( 1 , 4.9473215785617045e-08 )
( 1.0 , 4.9473215785617045e-08 )
( 10.0 , 5.387758712930894e-08 )
( 100.0 , 6.554984332791442e-08 )
( 1000.0 , 7.0092396374782e-08 )
    				};
    				\addlegendentry{$\nu = 0$}
    				
    				\addplot[color=asb, mark=pentagon] coordinates {
                        ( 0.001 , 4.79935772499081e-08 )
( 0.01 , 4.799454615247358e-08 )
( 0.1 , 4.80207453751088e-08 )
( 1 , 4.832527014540147e-08 )
( 1.1428571428571428 , 4.8379731480155155e-08 )
( 10.0 , 5.23919549248566e-08 )
( 100.0 , 6.453368845295137e-08 )
( 1000.0 , 6.94336872843975e-08 )
    				};
    				\addlegendentry{$\nu = 0.125$}
    				
    				\addplot[color=violet, mark=square] coordinates {
    					( 0.001 , 4.6161517357578606e-08 )
( 0.01 , 4.616295253024028e-08 )
( 0.1 , 4.617014712613163e-08 )
( 1 , 4.630549320978711e-08 )
( 1.3333333333333333 , 4.63795965586083e-08 )
( 10.0 , 4.9656985464102995e-08 )
( 100.0 , 6.188138571610672e-08 )
( 1000.0 , 6.706075952279486e-08 )
    				};
    				\addlegendentry{$\nu = 0.25$}
    				
    				\addplot[color=blue, mark=diamond] coordinates {
    					( 0.001 , 4.340267532153589e-08 )
( 0.01 , 4.34011561879436e-08 )
( 0.1 , 4.338177566491489e-08 )
( 1 , 4.327031399916604e-08 )
( 1.6 , 4.326536440051456e-08 )
( 10.0 , 4.536546666179684e-08 )
( 100.0 , 5.6895977435746595e-08 )
( 1000.0 , 6.214071235169597e-08 )
    				};
    				\addlegendentry{$\nu = 0.375$}
    				
    				\addplot[color=cyan, mark=o] coordinates {
    					( 0.001 , 4.0104799303394754e-08 )
( 0.01 , 4.009959404557521e-08 )
( 0.1 , 4.004912413585133e-08 )
( 1 , 3.964301038240013e-08 )
( 1.996007984031936 , 3.936489568004807e-08 )
( 10.0 , 3.9946789249207816e-08 )
( 100.0 , 4.962481333104947e-08 )
( 1000.0 , 5.4560377199008364e-08 )
    				};
    				\addlegendentry{$\nu \approx 0.5$}
    				    				
    			\end{loglogaxis}
    		\end{tikzpicture}
    		\caption{}
    	\end{subfigure}
     \begin{subfigure}{0.3\linewidth}
    		\centering
    		\begin{tikzpicture}[scale = 0.6]
    			\definecolor{asl}{rgb}{0.4980392156862745,0.,1.}
    			\definecolor{asb}{rgb}{0.,0.4,0.6}
    			\begin{loglogaxis}[
    				/pgf/number format/1000 sep={},
    				axis lines = left,
    				xlabel={$k$-value},
    				ylabel={$\| \widetilde{\bm{\sigma}} - \bm{\sigma}^h \|_{\Le}/\| \widetilde{\bm{\sigma}}\|_{\Le}$},
    				xmin=0.5*1e-3, xmax=2e3,
    				ymin=2e-8, ymax=5e-5,
    				xtick={1e-3,1e-2,1e-1,1,1e+1,1e+2,1e+3},
    				ytick={1e-8,1e-7,1e-6,1e-5,1e-4},
    				legend style={at={(0.05,0.05)},anchor= south west},
    				ymajorgrids=true,
    				grid style=dotted,
    				]
    				\addplot[color=asl, mark=triangle] coordinates {
    					( 0.001 , 4.8311184144745615e-06 )
( 0.01 , 4.8308264697635525e-06 )
( 0.1 , 4.828090960045397e-06 )
( 1 , 4.8156471528785104e-06 )
( 1.0 , 4.815647152878504e-06 )
( 10.0 , 5.077195955003104e-06 )
( 100.0 , 9.834781708095675e-06 )
( 1000.0 , 4.96461817095758e-05 )
    				};
    				\addlegendentry{$\nu = 0$}
    				
    				\addplot[color=asb, mark=pentagon] coordinates {
                        ( 0.001 , 4.915452228473412e-06 )
( 0.01 , 4.915179338961468e-06 )
( 0.1 , 4.9125981151791155e-06 )
( 1 , 4.898723876407559e-06 )
( 1.1428571428571428 , 4.898088891810889e-06 )
( 10.0 , 5.075451330523683e-06 )
( 100.0 , 8.80420642251835e-06 )
( 1000.0 , 4.249408229593084e-05 )
    				};
    				\addlegendentry{$\nu = 0.125$}
    				
    				\addplot[color=violet, mark=square] coordinates {
    					( 0.001 , 5.009384546849451e-06 )
( 0.01 , 5.009135708020126e-06 )
( 0.1 , 5.006755392635723e-06 )
( 1 , 4.991709781011849e-06 )
( 1.3333333333333333 , 4.989192798582966e-06 )
( 10.0 , 5.077318651683569e-06 )
( 100.0 , 7.460215983134531e-06 )
( 1000.0 , 3.239305068145383e-05 )
    				};
    				\addlegendentry{$\nu = 0.25$}
    				
    				\addplot[color=blue, mark=diamond] coordinates {
    					( 0.001 , 5.089674205295996e-06 )
( 0.01 , 5.089450676012784e-06 )
( 0.1 , 5.087292417037793e-06 )
( 1 , 5.0719091493486664e-06 )
( 1.6 , 5.066125414399832e-06 )
( 10.0 , 5.085318060888366e-06 )
( 100.0 , 6.043320789309587e-06 )
( 1000.0 , 1.9401855821078755e-05 )
    				};
    				\addlegendentry{$\nu = 0.375$}
    				
    				\addplot[color=cyan, mark=o] coordinates {
    					( 0.001 , 5.120750163512182e-06 )
( 0.01 , 5.120545504714897e-06 )
( 0.1 , 5.118567293133773e-06 )
( 1 , 5.104201093847943e-06 )
( 1.996007984031936 , 5.095786244654546e-06 )
( 10.0 , 5.099735191013673e-06 )
( 100.0 , 5.368462495515704e-06 )
( 1000.0 , 9.30481136233513e-06 )
    				};
    				\addlegendentry{$\nu \approx 0.5$}
    				    				
    			\end{loglogaxis}
    		\end{tikzpicture}
      \caption{}
    	\end{subfigure}
    	\caption{Respective errors for plane stress (a)-(c) and plane strain (d)-(e) relative to $\psi = k \chi$ for the cases of shear $\widetilde{\vb{u}}_s$, bi-axial tension $\widetilde{\vb{u}}_b$, and periodic displacements $\widetilde{\vb{u}}_p$.}
    	\label{fig:error}
\end{figure}
We observe that up to $\psi = \chi$ and even for $\psi = 1$ with a negligibly worse result, the error remains practically the same. At the same time, the solver fails for $\psi = 0$. Consequently, it appears that any not too small or too large value leads to accurate results. Clearly, a simple and valid choice is therefore $\psi = \chi$, which allows to derive the simplest form of the equations by diving by $\chi$. 

\subsection{Operator spectrum and stabilisation}
Analogously to \cref{sec:opnorm} for the three-dimensional case, we introduce the operator form of the two-dimensional variational problem in \cref{eq:var2d} as
\begin{align}
    &\text{find} \quad \bm{\sigma} \in \Y^h(\surf) && \text{s.t.} &&  A^h \bm{\sigma} = \vb{f} \qquad \text{in} \qquad [\Y^h(\surf)]' \, , 
\end{align}
and investigate its eigenvalues on the domain $\overline{\surf} = [-1,1]^2$ with a complete Neumann boundary $\curv_N = \partial \surf$. The finite element mesh is given by $3^2 = 9$-quadrilateral elements of cubic polynomial order, yielding $300$-degrees of freedom.

From our computations we observe that irrespective of the value the Poisson ratio $\nu$ or $\psi > 0$ take, we always retrieve $9$-zero eigenvalues and no negative eigenvalues for both plane stress and plane strain, implying a positive semi-definite operator. However, from the displacement-based formulation we expect only $3$-zero eigenvalues, accounting for rigid body motions in the plane. 
Observe that for $\psi = \chi$ the resulting field equation is equivalent to \cref{eq:2deq}, allowing us to reformulate the boundary value problem.   
\begin{definition}[Pure stress planar boundary value problem II]
\label{def:planarBVPII}
    Let $\chi = 1/(1+\nu)$ or $\chi = 1-\nu$ for plane stress or plane strain, respectively,  
    the complete boundary value problem reads
    \begin{subequations}
        \begin{align}
            - \Delta \bm{\sigma} - \hess(\tr \bm{\sigma}) - (\di \Di \bm{\sigma}) \one &= 2\sym \D \vb{f} + \dfrac{1}{\chi}(\di \vb{f})\one \, , && \text{in} && \surf \, , \\ 
            (\D \bm{\sigma})\vb{n} + \nabla \tr \bm{\sigma} \otimes \vb{n} - \con{\vb{f}}{\vb{n}}\one  &= \bm{\kappa} && \text{on} && \curv_N \, , \\
        \bm{\sigma} &= \widetilde{\bm{\sigma}}  && \text{on} && \curv_D \, .
    \end{align}
    \end{subequations}
\end{definition}
The variational form of \cref{def:planarBVPII} follows by testing and partial integration, and is well-posed due to \cref{co:well2d}.
\begin{definition}[Pure stress planar variational form II]
\label{def:var2DII}
    The variational form reads 
    \begin{align}
    &\int_\surf \con{\D \bm{\tau}}{ \D \bm{\sigma}} + \con{\Di \bm{\tau}}{\nabla \tr \bm{\sigma}} + \con{\nabla \tr \bm{\tau}}{\Di \bm{\sigma}} \, \dd \surf = \int_{\curv_N} \con{\bm{\tau}}{\bm{\kappa}} \, \dd \curv
     \notag \\
     &\qquad +  \int_\surf 2\con{\bm{\tau}}{\sym \D \vb{f}} +\dfrac{1}{\chi} \con{\tr \bm{\tau}}{\di \vb{f}}  \, \dd \surf  \qquad \forall \, \bm{\tau} \in \C_{\curv_D}^\infty(\overline{\surf})\otimes\Sym(2) \, , \label{eq:weak2DII} 
\end{align}
    with $\chi = 1/(1+\nu)$ or $\chi = 1-\nu$, for plane-stress or plane stain, respectively.
\end{definition}
The latter can clearly be derived from a functional.
\begin{definition}[Planar variational functional in pure stresses II]
\label{def:energy2DII}
    The weak formulation of the new boundary value problem can directly constructed as the variation of the functional 
    \begin{align}
    I(\bm{\sigma}) &= \int_\surf \dfrac{1}{2} \norm{\D \bm{\sigma}}^2 +  \con{\Di \bm{\sigma}}{\nabla \tr \bm{\sigma}}   - 2 \con{\bm{\sigma}}{ \sym \D \vb{f}} - \dfrac{1}{\chi} \con{\tr \bm{\sigma}}{\di \vb{f}} \, \dd \surf- \int_{\curv_N} \con{\bm{\sigma}}{\bm{\kappa}} \, \dd \curv \, , \notag
\end{align}
    with $\chi = 1/(1+\nu)$ or $\chi = 1-\nu$, for plane-stress or plane stain, respectively.
\end{definition}
Lastly, the discrete form can be derived by using distributions of the derivatives of the body forces.
\begin{definition}[Planar discrete variational form]
\label{def:var2DIId}
    The discrete variational form reads 
    \begin{align}
    &\int_\surf \con{\D \bm{\tau}}{ \D \bm{\sigma}} + \con{\Di \bm{\tau}}{\nabla \tr \bm{\sigma}} + \con{\nabla \tr \bm{\tau}}{\Di \bm{\sigma}} \, \dd \surf \\ 
    &\qquad =  2 \con{\bm{\tau}}{\sym\D \vb{f}}_{\mathcal{T}} + \dfrac{1}{\chi} \con{\tr \bm{\tau}}{ \di \vb{f}}_{\mathcal{T}} + \int_{\curv_N} \con{\bm{\tau}}{\bm{\kappa}} \, \dd \curv \qquad \forall \, \bm{\tau} \in \Y^h(\body)     \, , \notag  
\end{align}
    with $\chi = 1/(1+\nu)$ or $\chi = 1-\nu$, for plane-stress or plane stain, respectively.
\end{definition}
Now, repeating the benchmark with the new formulation reveals that the associated discrete operator $\widehat{A}^h$ always yields $3$-zero eigenvalues, being a sound result considering the analogy to the displacement-based formulation.

We apply to both formulations the benchmark from \cref{sec:chara} for a complete Dirichlet $\curv_D = \partial \surf$ and mixed Dirichlet Neumann boundary $\partial \surf = \curv_D \cup \curv_N$, using the periodic displacement field $\widetilde{\vb{u}}_p$ for the generation of appropriate boundary conditions. We test for both plane stress $\bm{\sigma}_\sigma$ and plane strain $\bm{\sigma}_\varepsilon$ with $\nu = 0.25$ using cubic quadrilateral finite elements.   
\begin{figure}
    	\centering
    	\begin{subfigure}{0.3\linewidth}
    		\centering
    		\begin{tikzpicture}[scale = 0.6]
    			\definecolor{asl}{rgb}{0.4980392156862745,0.,1.}
    			\definecolor{asb}{rgb}{0.,0.4,0.6}
    			\begin{loglogaxis}[
    				/pgf/number format/1000 sep={},
    				axis lines = left,
    				xlabel={Degrees of freedom},
    				ylabel={$\| \widetilde{\bm{\sigma}} - \bm{\sigma}^h \|_{\Le}/\| \widetilde{\bm{\sigma}}\|_{\Le}$},
    				xmin=0.9e+3, xmax=1.1e5,
    				ymin=1e-6, ymax=2e-1,
    				xtick={1e2,1e3,1e4,1e5},
    				ytick={1e-5,1e-5,1e-4,1e-3,1e-2,1e-1},
    				legend style={at={(0.95,0.95)},anchor= north east},
    				ymajorgrids=true,
    				grid style=dotted,
    				]
    				\addplot[color=asl, mark=triangle] coordinates {
    					( 2208 , 0.0011811342624151748 )
( 8463 , 7.781968257553548e-05 )
( 18768 , 1.557530491523361e-05 )
( 33123 , 4.9527867145745e-06 )
( 51528 , 2.0335889370543256e-06 )
    				};
    				\addlegendentry{$\bm{\sigma}_{\sigma}^{h}$}
    				
    				\addplot[color=violet, mark=square] coordinates {
                        ( 2208 , 0.0011931792077948446 )
( 8463 , 7.842910132821979e-05 )
( 18768 , 1.5691571335432072e-05 )
( 33123 , 4.989205040593749e-06 )
( 51528 , 2.0484500176206513e-06 )
    				};
    				\addlegendentry{$\bm{\sigma}_{\varepsilon}^{h}$}

    				\addplot[color=blue, mark=10-pointed star] coordinates {
    					( 2208 , 0.001195044259229365 )
( 8463 , 7.79942720589485e-05 )
( 18768 , 1.558294832889305e-05 )
( 33123 , 4.952424629739216e-06 )
( 51528 , 2.032957078008511e-06 )
    				};
    				\addlegendentry{$\bm{\sigma}_{\sigma s}^{h}$}
    				
    				\addplot[color=asb, mark=oplus] coordinates {
    					( 2208 , 0.0012043893672952432 )
( 8463 , 7.854993278918275e-05 )
( 18768 , 1.5693416304225657e-05 )
( 33123 , 4.987523579631606e-06 )
( 51528 , 2.0473695887190555e-06 )
    				};
    				\addlegendentry{$\bm{\sigma}_{\varepsilon s}^{h}$}
    				
    				\addplot[dashed,color=black, mark=none]
    				coordinates {
    					(5000, 0.5e-4)
    					(20000, 3.125e-06)
    				};
    			
    			\end{loglogaxis}

                \draw (3.1,0.5) 
    			node[anchor=south]{\tiny $\mathcal{O}(h^{4})$};

    		\end{tikzpicture}
    		\caption{}
    	\end{subfigure}
        \begin{subfigure}{0.3\linewidth}
    		\centering
    		\begin{tikzpicture}[scale = 0.6]
    			\definecolor{asl}{rgb}{0.4980392156862745,0.,1.}
    			\definecolor{asb}{rgb}{0.,0.4,0.6}
    			\begin{loglogaxis}[
    				/pgf/number format/1000 sep={},
    				axis lines = left,
    				xlabel={Degrees of freedom},
    				ylabel={$\| \widetilde{\bm{\sigma}} - \bm{\sigma}^h \|_{\Le}/\| \widetilde{\bm{\sigma}}\|_{\Le}$},
    				xmin=0.9e+3, xmax=1.1e5,
    				ymin=1e-6, ymax=2e-1,
    				xtick={1e2,1e3,1e4,1e5},
    				ytick={1e-5,1e-5,1e-4,1e-3,1e-2,1e-1},
    				legend style={at={(0.05,0.05)},anchor= south west},
    				ymajorgrids=true,
    				grid style=dotted,
    				]
    				\addplot[color=asl, mark=triangle] coordinates {
    					( 2208 , 0.08232824219246421 )
( 8463 , 0.08904808332353946 )
( 18768 , 0.08024552780623966 )
( 33123 , 0.07238331200927114 )
( 51528 , 0.0647776149245417 )
    				};
    				\addlegendentry{$\bm{\sigma}_{\sigma}^{h}$}
    				
    				\addplot[color=violet, mark=square] coordinates {
                        ( 2208 , 0.07435522854853867 )
( 8463 , 0.07479670821797677 )
( 18768 , 0.06838110609811651 )
( 33123 , 0.06226122731799847 )
( 51528 , 0.05588358357205787 )
    				};
    				\addlegendentry{$\bm{\sigma}_{\varepsilon}^{h}$}

    				\addplot[color=blue, mark=10-pointed star] coordinates {
    					( 2208 , 0.0010972962894630622 )
( 8463 , 7.150883377925193e-05 )
( 18768 , 1.4279184112098849e-05 )
( 33123 , 4.537441536683593e-06 )
( 51528 , 1.8625063694028735e-06 )
    				};
    				\addlegendentry{$\bm{\sigma}_{\sigma s}^{h}$}
    				
    				\addplot[color=asb, mark=oplus] coordinates {
    					( 2208 , 0.0011053238762182476 )
( 8463 , 7.193870905517421e-05 )
( 18768 , 1.4362732834402402e-05 )
( 33123 , 4.563733157142203e-06 )
( 51528 , 1.873244858916143e-06 )
    				};
    				\addlegendentry{$\bm{\sigma}_{\varepsilon s}^{h}$}
    				
    				\addplot[dashed,color=black, mark=none]
    				coordinates {
    					(5000, 0.5e-4)
    					(20000, 3.125e-06)
    				};

                    \addplot[dashed,color=black, mark=none]
    				coordinates {
    					(15000, 3e-2)
    					(50000, 0.02220248413476856)
    				};
    			
    			\end{loglogaxis}

                \draw (3.1,0.5) 
    			node[anchor=south]{\tiny $\mathcal{O}(h^{4})$};

                \draw (4.75,4.05) 
    			node[anchor=south]{\tiny $\mathcal{O}(\sqrt{h})$};

    		\end{tikzpicture}
    		\caption{}
    	\end{subfigure}
    	\begin{subfigure}{0.3\linewidth}
    		\centering
    		\definecolor{xfqqff}{rgb}{0.4980392156862745,0,1}
\definecolor{zzttqq}{rgb}{0.,0.4,0.6}
\begin{tikzpicture}[scale =0.7, line cap=round,line join=round,>=triangle 45,x=1cm,y=1cm]
\clip(5.5, 5) rectangle (13,9.5);
\fill[line width=0.7pt,color=zzttqq,fill=zzttqq,fill opacity=0.10000000149011612] (6,7) -- (6,9) -- (12,9) -- (12,7) -- cycle;
\draw [line width=0.7pt,color=xfqqff] (6,7)-- (6,9);
\draw [line width=0.7pt,color=xfqqff] (6,9)-- (12,9);
\draw [line width=0.7pt,color=zzttqq, dashed] (12,9)-- (12,7);
\draw [line width=0.7pt,color=zzttqq, dashed] (12,7)-- (6,7);
\draw [-to,line width=0.7pt] (9,8) -- (10.5,8);
\draw [-to,line width=0.7pt] (9,8) -- (9,9.5);
\draw (6.5,9) node[anchor=south,color=xfqqff] {$\curv_D$};
\draw (12,7.5) node[anchor=west,color=zzttqq] {$\curv_N$};
\draw (7,7.5) node[anchor=south west] {$\surf$};
\begin{scriptsize}
\draw [fill=black] (9,8) circle (2.5pt);
\end{scriptsize}
\end{tikzpicture}
    		\caption{}
    	\end{subfigure}
    	\caption{Convergence of non-stabilised forms for plane stress $\bm{\sigma}^h_\sigma$ and plane strain $\bm{\sigma}^h_\varepsilon$ versus the respective stabilised forms denoted with $\bm{\sigma}^h_{\sigma s}$ and $\bm{\sigma}^h_{\varepsilon s}$. In (a) the complete boundary is Dirichlet $\curv_D = \partial \surf$, whereas in (b) the boundary is mixed between Dirichlet and Neumann $\partial\surf = \curv_D \cup \curv_N$ as per the depiction in (c).}
    	\label{fig:st2d}
\end{figure}
From the convergence rates in \cref{fig:st2d} we observe that for a complete Dirichlet boundary both formulations are stable and converge optimally. However, for a mixed Dirichlet-Neumann boundary the initial formulation barely converges at all at the suboptimal rate $\mathcal{O}(\sqrt{h})$. In contrast, the stabilised formulation is robust and yields optimal convergence irrespective of the boundary condition in use.   

\section{Conclusions and outlook}
This work introduces new formulations of isotropic linear elasticity in terms of the stress tensor. The formulations are based on the Beltrami-Michell equations, which in the three-dimensional case we symmetrised and stabilised. In the two-dimensional cases of plane stress and plane strain, the Beltrami-Michell equations characterise only the mean stress and are thus insufficient for a variational formulation of the full stress tensor. Consequently, we introduced a new boundary value problem in two-dimensions by adding to the Beltrami-Michell equations a weaker form of static equilibrium scaled with a constant, which we later characterised numerically. In addition, for all the presented variational formulations we proposed the treatment of constant and piece-wise constant body forces via distributions.     

Both, the modified Beltrami-Michell equations in three dimensions as well as the newly introduced equations for two-dimensions, were proven to be well-posed on a domain with a complete Dirichlet boundary, which was verified by numerical benchmarks yielding optimal convergence. In the case of mixed Dirichlet-Neumann boundary conditions stabilisation is required in both three- and two dimensions. In three-dimensions we introduced stabilisation by adding a weighted zero-sum term of a weaker form of static equilibrium into the boundary value problem. We investigated the constant weight $\omega$ of the stabilised form both analytically and numerically, concluding that the result $\omega > \chi$ guarantees stability, which we verified by spectral analysis and convergence tests. In the two-dimensional case we derived an alternative boundary value problem via a novel identity of the Laplacian of a symmetric second order tensor. The formulation was shown to be robust via spectral and convergence analysis irrespective of the split of the boundary between Dirichlet and Neumann. The newly introduced equations are summarised in \cref{ap:eqs}. Considering our observations on the numerical behaviour of the variational forms, we recommend to only use the stabilised versions in order to guarantee robust computations.

Despite being classical instruments for the (analytical) determination of stresses in mechanics, the Beltrami-Michell equations were not considered a suitable starting point for approximate solutions to the stress state, e.g. using finite element computations. As shown in this work, this is likely due to the inadequacy of the classical form for a variational framework. Thus, the newly introduced equations and variational forms in this paper represent the foundation for future numerical application of pure stress formulations of linear elasticity, for example for thermal-hydraulic-mechanical processes \cite{WANG2022109752,10.2118/182595-MS,10.2118/193931-MS}. A clear advantage of the numerical framework is its flexibility with respect to the geometry of the domain and variation of the material constants within it.
Specifically in material characterisation \cite{Alshaya2021} and in the theory of mixtures \cite{MUTI2015140} we expect the equations to be beneficial, considering the intensive research effort currently being invested into the reconstruction of computational representations of such materials for numerical simulations \cite{Seibert,Seibert2024,DURETH2023105608}.

The equations presented in this manuscript are restricted to the case of isotropic linear elasticity. Investigation of a more general constitutive setting is a subject for future works. 

\bibliographystyle{spmpsci}   

\footnotesize
\bibliography{ref}   

\normalsize
\appendix

\section{Glossary of stress equations} \label{ap:eqs}
\begin{definition}[Pure stress boundary value problem for solids I]
    The field equations and boundary conditions of the novel three-dimensional stress boundary value problem read
    \begin{subequations}
        \begin{align}
            -\Delta \bm{\sigma} - \dfrac{1}{1+\nu} [\hess( \tr \bm{\sigma} ) + (\di \Di \bm{\sigma}) \one ]  &= 2 \sym \D \vb{f}   + \dfrac{1+\nu^2}{1-\nu^2} (\di \vb{f}) \one && \text{in} && \body \, ,  \\
            (\D \bm{\sigma})\vb{n} + \dfrac{1}{1+\nu} [\nabla \tr \bm{\sigma} \otimes \vb{n} - \con{\vb{f}}{\vb{n}}\one]  &= \bm{\kappa} && \text{on} && \surf_N \, , \\
         \bm{\sigma} &= \widetilde{\bm{\sigma}} && \text{on} && \surf_D  \, .
        \end{align}
    \end{subequations}
    where $\bm{\kappa}$ is a mixed measure of equilibrium and compatibility on the Neumann boundary. The problem yields a stable variational form the case of a complete Dirichlet boundary condition, or a mixed Dirichlet-Neumann boundary provided the material is incompressible $\nu = 0.5$.   
\end{definition}
\begin{definition}[Discrete variational problem for solids I]
    The corresponding discrete variational problem reads
    \begin{align}
    &\int_\body \con{\D \bm{\tau}}{ \D \bm{\sigma}} + \dfrac{1}{1+\nu} [\con{\Di \bm{\tau}}{\nabla \tr \bm{\sigma}} + \con{\nabla \tr \bm{\tau}}{\Di \bm{\sigma}}]   \, \dd \body \notag   \\
    &\qquad = \int_{\surf_{N}} \con{\bm{\tau}}{\bm{\kappa}} \, \dd \surf + 2 \con{\bm{\tau}}{\sym\D \vb{f}}_{\mathcal{T}} + \dfrac{1+\nu^2}{1-\nu^2} \con{\tr \bm{\tau}}{ \di \vb{f}}_{\mathcal{T}} \qquad \forall \, \bm{\tau} \in \X^h(\body) \, , 
\end{align}
and is stable either for a complete Dirichlet boundary $\surf_D = \partial \body$ or a mixed boundary assuming incomprehensibility $\nu = 0.5$.
\end{definition}
\begin{definition}[Pure stress boundary value problem for solids II]
    The field equations and boundary conditions of the stabilised symmetric three-dimensional stress boundary value problem read
    \begin{subequations}
        \begin{align}
            -\Delta \bm{\sigma} - \dfrac{1}{1+\nu} [\hess( \tr \bm{\sigma} ) + (\di \Di \bm{\sigma}) \one ] -\omega \sym \D \Di \bm{\sigma}  &= (2 + \omega) \sym \D \vb{f}   + \dfrac{1+\nu^2}{1-\nu^2} (\di \vb{f}) \one && \text{in} && \body \, ,  \\
            (\D \bm{\sigma})\vb{n} + \dfrac{1}{1+\nu} [\nabla \tr \bm{\sigma} \otimes \vb{n} - \con{\vb{f}}{\vb{n}}\one]  -\omega (\vb{f} \otimes \vb{n})  &= \bm{\kappa} && \text{on} && \surf_N \, , \\
         \bm{\sigma} &= \widetilde{\bm{\sigma}} && \text{on} && \surf_D  \, .
        \end{align}
    \end{subequations}
\end{definition}
\begin{definition}[Stabilised discrete variational problem for solids II]
    The stabilised discrete variational problem reads
    \begin{align}
    &\int_\body \con{\D \bm{\tau}}{ \D \bm{\sigma}} + \dfrac{1}{1+\nu} [\con{\Di \bm{\tau}}{\nabla \tr \bm{\sigma}} + \con{\nabla \tr \bm{\tau}}{\Di \bm{\sigma}}] + \omega \con{\Di \bm{\tau}}{\Di \bm{\sigma}}  \, \dd \body \notag   \\
    &\qquad = \int_{\surf_{N}} \con{\bm{\tau}}{\bm{\kappa}} \, \dd \surf + (2+\omega) \con{\bm{\tau}}{\sym\D \vb{f}}_{\mathcal{T}} + \dfrac{1+\nu^2}{1-\nu^2} \con{\tr \bm{\tau}}{ \di \vb{f}}_{\mathcal{T}} \qquad \forall \, \bm{\tau} \in \X^h(\body) \, ,  \label{eq:var3dII}
\end{align}
and is coercive for any choice $\omega > \chi$ even for a mixed Dirichlet and Neumann boundary irrespective of the value of the Poisson ratio.
\end{definition}
\begin{definition}[Pure stress planar boundary value problem I]
    Let $\chi = 1/(1+\nu)$ or $\chi = 1-\nu$ for plane stress or plane strain, respectively, and $\psi > 0$, the complete boundary value problem reads
    \begin{subequations}
        \begin{align}
        -\psi \sym(\D \Di \bm{\sigma}) - (\chi \Delta \tr \bm{\sigma})\one &= \psi \sym \D \vb{f} + (\di \vb{f})\one && \text{in} && \surf \, , \\ 
        \chi \con{\nabla \tr \bm{\sigma}}{\vb{n}} \one - \psi (\vb{f} \otimes \vb{n}) &= \bm{\kappa}  && \text{on} && \curv_N \, , \\
        \bm{\sigma} &= \widetilde{\bm{\sigma}}  && \text{on} && \curv_D \, .
    \end{align}
    \end{subequations}
    The boundary value problem leads to an optimally convergent discrete formulation only for a complete Dirichlet boundary.
\end{definition}
\begin{definition}[Pure stress planar discrete variational form I]
    The variational form reads 
    \begin{align}
    \int_\surf \psi \con{\Di \bm{\tau}}{ \Di \bm{\sigma} } + \chi \con{\nabla \tr \bm{\tau}}{ \nabla \tr \bm{\sigma} } \, \dd \surf 
     &=  \psi \con{\bm{\tau}}{\sym\D \vb{f}}_{\mathcal{T}} + \con{\tr \bm{\tau}}{ \di \vb{f}}_{\mathcal{T}}   
     \\& \qquad + \int_{\curv_N} \con{\bm{\tau}}{\bm{\kappa}} \, \dd \curv   \qquad \forall \, \bm{\tau} \in \Y^h(\surf) \, ,  \notag
\end{align}
    with $\psi > 0$ and $\chi = 1/(1+\nu)$ or $\chi = 1-\nu$, for plane-stress or plane stain, respectively. This form converges optimally only for a complete Dirichlet boundary $\curv_D = \partial \surf$.
\end{definition}
\begin{definition}[Pure stress planar boundary value problem II]
    Let $\chi = 1/(1+\nu)$ or $\chi = 1-\nu$ for plane stress or plane strain, respectively,  
    the complete boundary value problem reads
    \begin{subequations}
        \begin{align}
            - \Delta \bm{\sigma} - \hess(\tr \bm{\sigma}) - (\di \Di \bm{\sigma}) \one &= 2\sym \D \vb{f} + \dfrac{1}{\chi}(\di \vb{f})\one \, , && \text{in} && \surf \, , \\ 
            (\D \bm{\sigma})\vb{n} + \nabla \tr \bm{\sigma} \otimes \vb{n} - \con{\vb{f}}{\vb{n}}\one  &= \bm{\kappa} && \text{on} && \curv_N \, , \\
        \bm{\sigma} &= \widetilde{\bm{\sigma}}  && \text{on} && \curv_D \, .
    \end{align}
    \end{subequations}
\end{definition}
\begin{definition}[Pure stress planar discrete variational form II]
    The discrete variational form reads 
    \begin{align}
    &\int_\surf \con{\D \bm{\tau}}{ \D \bm{\sigma}} + \con{\Di \bm{\tau}}{\nabla \tr \bm{\sigma}} + \con{\nabla \tr \bm{\tau}}{\Di \bm{\sigma}} \, \dd \surf \\ 
    &\qquad =  2 \con{\bm{\tau}}{\sym\D \vb{f}}_{\mathcal{T}} + \dfrac{1}{\chi} \con{\tr \bm{\tau}}{ \di \vb{f}}_{\mathcal{T}} + \int_{\curv_N} \con{\bm{\tau}}{\bm{\kappa}} \, \dd \curv \qquad \forall \, \bm{\tau} \in \Y^h(\body)     \, , \notag  
\end{align}
    with $\chi = 1/(1+\nu)$ or $\chi = 1-\nu$, for plane-stress or plane stain, respectively. The form yields optimal convergence also for a mixed Dirichlet and Neumann boundary.
\end{definition}

\section{Mathematical formulae}
In this work we introduced the following original three-dimensional and two-dimensional identities 
\begin{subequations}
    \begin{align}
        \skw\Curl\bm{\sigma} &= \dfrac{1}{2} \Anti( \Di \bm{\sigma} - \nabla \tr \bm{\sigma} ) \qquad \forall \,  \bm{\sigma} \in \C^\infty(\overline{\body}) \otimes \Sym(3) \, ,  \\
        \Delta \bm{\sigma} &= 2\sym (\D \Di \bm{\sigma}) - \hess (\tr \bm{\sigma}) + [ \Delta \tr \bm{\sigma}  - \di \Di \bm{\sigma}] \one \qquad \forall \,  \bm{\sigma} \in \C^\infty(\overline{\surf}) \otimes \Sym(2) \, ,
    \end{align}
\end{subequations}
which were used in subsequent lemmas. Further, we made use of the following classical three-dimensional results 
\begin{subequations}
\begin{align}
        \inc \bm{\sigma} &= 2\sym (\D \Di \bm{\sigma}) - \Delta \bm{\sigma} - \hess (\tr \bm{\sigma}) + [ \Delta \tr \bm{\sigma}  - \di \Di \bm{\sigma}] \one \, , \\
        \tr (\inc \bm{\sigma}) &= \Delta \tr \bm{\sigma} - \di \Di \bm{\sigma} \, , \\
    \tr (\inc [ (\tr \bm{\sigma}) \one ])  & = 2 \Delta \tr \bm{\sigma} \, , \\
    \inc [ (\tr \bm{\sigma}) \one ] &= (\Delta \tr \bm{\sigma}) \one -  \hess (\tr \bm{\sigma})  \, , \\
    \Curl( \skw \bm{\sigma}) &= (\di \axl \skw \bm{\sigma}) \one - (\D \axl \skw \bm{\sigma} )^T \, ,
    \end{align}
\end{subequations}
whose thorough derivation can be found in \cite{Lewintan2021}.

\section{An observation on planar linear elasticity} \label{ap:a}
In following we characterise the curl of the stress tensor in planar linear elasticity. We start by introducing the planar Navier-Cauchy equations.

\subsection{The planar Navier-Cauchy equations}
Starting with the boundary value problem of linear elasticity, we can derive the planar Navier-Cauchy equations by inserting the constitutive law for plane stress or plane stain, respectively.
In the case of plane stress $\sigma_{i3} = \sigma_{3i} = 0$, the constitutive equation reads
\begin{align}
    \bm{\sigma} = \Cm \bm{\varepsilon} = \dfrac{E}{1-\nu^2}[(1-\nu)\bm{\varepsilon} + (\nu \tr \bm{\varepsilon}) \one] = \dfrac{E}{1-\nu^2}[(1-\nu)\sym \D\vb{u} + (\nu \tr \sym \D\vb{u}) \one] \, .
\end{align}
Now, recalling the identities
\begin{align}
    &\Di \sym\D \vb{u} = \dfrac{1}{2}(\Delta \vb{u} + \nabla \di \vb{u}) \, , && \tr \sym \D \vb{u} = \tr \D \vb{u} = \di \vb{u} \, , &&
    \Di[(\di \vb{u})\one] = \nabla \di \vb{u} \, ,
\end{align}
we find the equations
\begin{subequations}
    \begin{align}
    -\Di \Cm \bm{\varepsilon} =  -\dfrac{E}{2(1+\nu)} \Delta \vb{u}  - \dfrac{E}{2(1-\nu)} \nabla \di \vb{u}  &=\vb{f}  && \text{in} && \surf \, , \\
        \left [ \dfrac{E}{1+\nu} \sym \D \vb{u} + \dfrac{E \nu}{1-\nu^2} (\di \vb{u}) \one  \right ] \vb{n} &= \bm{\mu}  && \text{on} && \curv_N \,  , \\
        \vb{u} &= \widetilde{\vb{u}} && \text{on} && \curv_D \, ,
\end{align}
\end{subequations}
being the Navier-Cauchy equations of plane stress.

In the case of plane strain the constitutive equation reads
\begin{align}
    \bm{\sigma} = \Cm \bm{\varepsilon} = \dfrac{E}{1+\nu}\left[\bm{\varepsilon} + \dfrac{\nu}{1-2\nu} ( \tr \bm{\varepsilon}) \one \right] = \dfrac{E}{1+\nu}\left[\sym \D\vb{u} + \dfrac{\nu}{1-2\nu} ( \tr \sym \D\vb{u}) \one \right] \, ,
\end{align}
leading to the equations
\begin{subequations}
    \begin{align}
    -\Di \Cm \bm{\varepsilon} =  -\dfrac{E}{2(1+\nu)} \Delta \vb{u}  - \dfrac{E}{2(1+\nu)(1-2\nu)} \nabla \di \vb{u}  &=\vb{f}  && \text{in} && \surf \, , \\
        \left [ \dfrac{E}{1+\nu} \sym \D \vb{u} + \dfrac{E \nu}{(1+\nu)(1-2\nu)} (\di \vb{u}) \one  \right ]  \vb{n} &= \bm{\mu}  && \text{on} && \curv_N \,  , \\
        \vb{u} &= \widetilde{\vb{u}} && \text{on} && \curv_D \, ,
\end{align}
\end{subequations}
being the Navier-Cauchy equations of plane strain.
\subsection{The curl of the stress tensor}
We note that the skew-symmetric part of the gradient of the stress tensor is proportional to the curl of the stress tensor $\Rot \bm{\sigma}$. In order to find an identity concerning $\Rot \bm{\sigma}$, we observe
\begin{align}
    \Rot \sym \Di \vb{u} &= \Rot (\D \vb{u})^T = \dfrac{1}{2} \bm{R} \Delta \vb{u} - \dfrac{1}{2} \bm{R} \nabla \di \vb{u} =  \dfrac{1}{2} (\bm{R} \Delta \vb{u} -  \nabla^\perp \di \vb{u})  \, ,  \\
    \Rot[(\tr\sym \D \vb{u})\one] &= \Rot[(\di \vb{u})\one] = -\bm{R}\nabla \di \vb{u} = - \nabla^\perp \di \vb{u} \, . 
\end{align}
Now, the curl of the the stress tensor given by the constitutive equation of plane stress reads
\begin{align}
    \Rot \bm{\sigma} = \Rot \Cm \sym \D \vb{u} &= 
    \dfrac{E}{2(1+\nu)} (\bm{R} \Delta \vb{u} - \nabla^\perp \di \vb{u}) - \dfrac{E \nu}{1-\nu^2}\nabla^\perp \di \vb{u}
    \notag \\
    &= \dfrac{E}{2(1+\nu)} \bm{R} \Delta \vb{u}  - \dfrac{E}{2(1-\nu)}\nabla^\perp \di \vb{u} \, .
\end{align}
By rotating the Navier-Cauchy equation of plane stress we derive
\begin{align}
      \dfrac{E}{2(1+\nu)} \bm{R}\Delta \vb{u} =  - \dfrac{E}{2(1-\nu)} \nabla^\perp \di \vb{u}- \bm{R}\vb{f}  \, ,
\end{align}
leading to
\begin{align}
    \Rot \bm{\sigma} = - \dfrac{E}{1-\nu} \nabla^\perp \di \vb{u} - \bm{R} \vb{f} \, .
\end{align}
We can express
\begin{align}
    \nabla^\perp \di \vb{u} = \bm{R} \nabla \di \vb{u} = \bm{R}\Di[(\tr \sym \D \vb{u})\one] = \bm{R}\Di[(\tr [\A \bm{\sigma}])\one] = \bm{R} \nabla \tr(\A \bm{\sigma}) = \nabla^\perp \tr (\A \bm{\sigma})  \, , 
\end{align}
which by using the compliance relation of plane stress reads
\begin{align}
    \nabla^\perp \di \vb{u} = \nabla^\perp \tr(\A \bm{\sigma}) = \dfrac{1}{E} \nabla^\perp [(1+\nu) \tr \bm{\sigma} - 2 \nu \tr \bm{\sigma}] = \dfrac{1-\nu}{E} \nabla^\perp \tr \bm{\sigma} \, .    
\end{align}
Thus, we find for the curl of the stress tensor 
\begin{align}
    \Rot\bm{\sigma} + \nabla^\perp \tr \bm{\sigma} = - \bm{R} \vb{f} \, , \label{eq:rotS}
\end{align} 
being an alternative form of equilibrium for planar problems.

Following the same procedure for plane strain we observe
\begin{align}
    \Rot \bm{\sigma} = \dfrac{E}{2(1+\nu)} \bm{R} \Delta \vb{u} - \dfrac{E}{2(1+\nu)(1-2\nu)}\nabla^\perp \di \vb{u} \, .
\end{align}
The rotated Navier-Cauchy equation of plane strain yields
\begin{align}
    \dfrac{E}{2(1+\nu)} \bm{R}\Delta \vb{u} =-\dfrac{E}{2(1+\nu)(1-2\nu)} \nabla^\perp \di \vb{u} - \bm{R}\vb{f} \, , 
\end{align}
leading to
\begin{align}
    \Rot\bm{\sigma} = -\dfrac{E}{(1+\nu)(1-2\nu)} \nabla^\perp \di \vb{u} - \bm{R}\vb{f} \, .
\end{align}
We compute
\begin{align}
    \nabla^\perp \di \vb{u} = \nabla^\perp \tr(\A \bm{\sigma}) = \dfrac{1+\nu}{E} \nabla^\perp [ \tr \bm{\sigma} - 2 \nu \tr \bm{\sigma}] = \dfrac{(1+\nu)(1-2\nu)}{E} \nabla^\perp \tr \bm{\sigma} \, ,    
\end{align}
where we relied on the plane strain constitutive equation. Thus, we obtain the same relation as for plane stress in \cref{eq:rotS}. 

Consequently, we find the new differential identities 
\begin{align}
    &\Rot \bm{\sigma} + \nabla^\perp \tr \bm{\sigma} = \bm{R} \Di \bm{\sigma} \, , && \Di \bm{\sigma} = \nabla \tr \bm{\sigma} - \bm{R}\Rot \bm{\sigma} \, .
\end{align}
relating the divergence of a symmetric two-dimensional second order tensor $\bm{\sigma} = \bm{\sigma}^T$ to its curl.  

\end{document}